\documentclass[reqno, twoside, a4paper]{amsart}

\usepackage{amsmath}
\usepackage{amssymb}
\usepackage{amsthm}
\usepackage{bookmark}
\usepackage{indentfirst}                                                                         
\usepackage{multicol}
\usepackage{float}
\usepackage{enumerate}
\usepackage{xcolor}
\usepackage{cleveref}
\usepackage{booktabs}
\usepackage{mathrsfs}
\usepackage{tikz-cd}  
\usepackage{MnSymbol}
\usepackage{pgflibraryplothandlers}
\usepackage{amsrefs}

\newtheorem{theorem}{Theorem}[section]
\newtheorem{lemma}[theorem]{Lemma}
\newtheorem{proposition}[theorem]{Proposition}
\newtheorem{corollary}[theorem]{Corollary}

\theoremstyle{definition}

\newtheorem{definition}[theorem]{Definition}
\newtheorem{remark}[theorem]{Remark}

\newtheorem{example}[theorem]{Example}

\numberwithin{equation}{section}

\makeatletter
\renewcommand*\env@matrix[1][*\c@MaxMatrixCols c]{%
  \hskip -\arraycolsep
  \let\@ifnextchar\new@ifnextchar
  \array{#1}}
\makeatother

\usepackage{tikz}

\usetikzlibrary{patterns,decorations.pathreplacing, calc, shapes}
\usetikzlibrary{arrows,decorations.pathmorphing,backgrounds,positioning,fit,petri}
\tikzstyle{edge} = [fill,opacity=.2,fill opacity=.5,line cap=round, line join=round, line width=50pt]

\tikzset{bullet/.style={
shape = circle,fill = black, inner sep = 0pt, outer sep = 0pt, minimum size = 0.35em, line width = 0pt, draw=black!100}}

\tikzset{circle/.style={
shape = circle,fill = none, inner sep = 0pt, outer sep = 0pt, minimum size = 0.35em, line width = 1pt, draw=black!100}}

\tikzset{rectangle/.style={
shape = rectangle,fill = white, inner sep = 0pt, outer sep = 0pt, minimum size = 0.35em, line width = 0pt, draw=black!100}}

\tikzset{empty/.style={
shape = circle,fill = white, inner sep = 0pt, outer sep = 0pt, minimum size = 0.35em, line width = 0pt, draw=white!100}}

\tikzset{xmark/.style={
shape = x,fill = white, inner sep = 0pt, outer sep = 0pt, minimum size = 0em, line width = 0pt, draw=white!100}}

\tikzset{longrectangle/.style={
inner sep = 1em,
rectangle,
minimum size=1em,
very thick,
draw=black!100, 
}}

\DeclareMathOperator{\Def}{Def}

\tikzset{label distance=-0.15em}

\tikzset{cross/.style={cross out, draw=black, fill=none, minimum size=2*(#1-\pgflinewidth), inner sep=0pt, outer sep=0pt}, cross/.default={2pt}}

\tikzset{font=\scriptsize}

\begin{document}
\title[Deformations of weighted homogeneous surface singularities]{Deformations of weighted homogeneous surface singularities with big central node}

\author[J. Jeon]{Jaekwan Jeon}

\address{Department of Mathematics, Chungnam National University, Daejeon 34134, Korea}

\email{jk-jeon@cnu.ac.kr}

\author[D. Shin]{Dongsoo Shin}

\address{Department of Mathematics, Chungnam National University, Daejeon 34134, Korea}

\email{dsshin@cnu.ac.kr}



\begin{abstract}
We prove Koll\'{a}r conjecture for weighted homogeneous surface singularities with big central node. More precisely, we show that every irreducible component of the deformation space of the singularity is parametrized by a certain partial resolution which is known as a $P$-resolution. 
\end{abstract}

\maketitle
\tableofcontents

\section{Introduction}
J. Koll\'{a}r and N. I. Shepherd-Barron(K-SB \cite{KSB}) proved that each irreducible component of the deformation space of a quotient surface singularity is parametrized by certain partial resolution, known as a \textit{$P$-resolution}. Building on this result, J. Koll\'{a}r(\cite{MR1144527}) introduced a conjecture stating that  every irreducible component of the deformation space of a rational surface singularity is parameterized by a certain partial modification of the singularity, known as a \textit{$P$-modification}. 

A $P$-resolution of a singularity $(X, p)$ is a partial resolution $f: Y \to X$ such that $Y$ has only singularities of class T and the canonical divisor $K_Y$ of $Y$ is $f$-relatively ample. A singularity of class T is a cyclic quotient surface singularity admitting a $\mathbb{Q}$-Gorenstein smoothing. Since every irreducible component of the deformation space of a rational surface singularity contains a smoothing, there is a natural map $\phi_P : \mathscr{P}(X) \to \mathscr{C}(X)$ from the set of all $P$-resolutions to the set of all irreducible components of the deformation space. In this aspect, Koll\'{a}r conjecture means the surjectivity of the map $\phi_P$.

We prove the conjecture for weighted homogeneous surface singularities. The singularity is a normal surface singularity that admits a good $\mathbb{C}^*$-action(ref. Orlik-Wagreich \cite{Orl} for details). The dual resolution graph of the singularity is star-shaped, meaning it consists of a central node of degree $-d$ and $t$ branches. Then we prove:
\begin{theorem}
Let $(X,0)$ be a weighted homogeneous surface singularity with $d \geq t + 3$. We construct all $P$-resolutions of $(X, 0)$ and show that $P$-resolutions parametrize all irreducible components of $\Def(X,0)$, that is, the map $\phi_P$ is surjective. This implies that Koll\'{a}r conjecture holds for the singularity.
\end{theorem}
To prove the theorem, we use a deformation theory of sandwiched surface singularities. 

A sandwiched surface singularity is a normal surface singularity admitting a birational morphism to $\mathbb{C}^2$. It is well known that a sandwiched surface singularity is rational. T. de Jong and D. van Straten(JS, \cite{deJ}) proved that any one-parameter smoothing of a sandwiched surface singularity can be described by a deformation of a plane curve singularity, which is known as a picture deformation. Since every irreducible component of a rational surface singularity is a smoothing component, this work implies that picture deformations parametrize irreducible components of a sandwiched surface singularity. If we denote the set of all picture deformations of a sandwiched surface singularity by $\mathscr{I}(X)$, then there exists a natural map $\phi_I : \mathscr{C}(X) \to \mathscr{I}(X)$.

On the other hand, a picture deformation has combinatorial aspects(roughly speaking, intersection relations between curves) and therefore it induces a matrix, which is called an incidence matrix. Furthermore, there exist combinatorial equations that every incidence matrix satisfies. Matrices satisfying these equations are called combinatorial incidence matrices. Let the set of all combinatorial incidence matrices be denoted by $C\mathscr{I}(X)$. Then we have a natural map $\phi_C : \mathscr{I}(X) \to C\mathscr{I}(X)$. The injectivity and surjectivity of the map are not generally known(JS, \cite{deJ}*{p.485}).

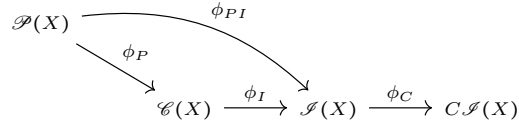
\begin{figure}
\centering
  \begin{tikzcd}
    \mathscr{P}(X) \arrow[dr, "\phi_P"] \arrow[drr, "\phi_{PI}", bend left=25] &                                      &                                     &                 \\
                                                                               & \mathscr{C}(X) \arrow[r, "\phi_I"]   & \mathscr{I}(X) \arrow[r, "\phi_C"]  & C\mathscr{I}(X) 
  \end{tikzcd}
\caption{Correspondence for sandwiched surface singularities}
\label{fig:Correspondence for sandwiched singularities}
\end{figure}

H. Park and D. Shin(\cite{park2022deformations}*{Theorem~6.5}) establish a map $\phi_{PI} : \mathscr{P}(X) \to \mathscr{I}(X)$ by using the minimal model program. The correspondences of the sets that we have discussed so far are illustrated in Figure~\ref{fig:Correspondence for sandwiched singularities}.
If we can find a $P$-resolution $f \in \mathscr{P}(X)$ such that $\phi_{PI}(f) = M$ for each combinatorial incidence matrix $M \in C\mathscr{I}(X)$, then it follows that the map $\phi_{PI}$ is surjective. Furthermore, if we can prove that $\phi_I$ is injective, then the map $\phi_P$ becomes surjective, and as a result, Koll\'{a}r conjecture holds.

Park-Shin(\cite{park2022deformations}*{Theorem~12.6}) proved the injectivity of the map $\phi_I$ for weighted homogeneous surface singularities with $d \geq t + 2$. In this article, we prove the surjectivity of the map $\phi_{PI}$ for the case $d \geq t + 3$. 

The strategy is as follows. We classify combinatorial incidence matrices $M$ of a weighted homogeneous surface singularity into two cases A and B(Theorem~\ref{thm:caseofcombinatorialincidencematrix}). For the classification, we observe that the combinatorial equations that the matrix $M$ satisfies contain the combinatorial equations of cyclic quotient surface singularities. In some sense, the matrix $M$ is a combination of combinatorial incidence matrices of cyclic quotient surface singularities with special restrictions. We prove that the matrix $M$ must be one of the cases $A$ or $B$ because of the restrictions.

We know that, for cyclic quotient surface singularities, every map in Figure~\ref{fig:Correspondence for sandwiched singularities} is bijective. Therefore if a combinatorial incidence matrix of a cyclic quotient surface singularity is given, then we can find the corresponding $P$-resolution. Since we have already observed that the matrix $M$ is a combination of combinatorial incidence matrices of cyclic quotient surface singularities, we construct the corresponding $P$-resolution of the matrix $M$ by combining the $P$-resolutions of the cyclic quotient surface singularities.

Finally, we verify that $\phi_{PI}(f) = M$ by applying MMP method of Park-Shin(\cite{park2022deformations}).

\subsection*{Acknowledgements} 
This article is a revision of Ph.D dissertation of J. Jeon presented at Department of Mathematics, Chungnam National University, Daejeon, Korea in 2023.

\section{Sandwiched surface singularities}
We will briefly review some definitions and theorems based on the work of M. Spivakovsky\cite{S} and de Jong-van Straten (\cite{deJ}).

\subsection{Sandwiched surface singularities}
  A sandwiched surface singularity $(X,p)$ is a normal surface singularity admitting a birational morphism $X \to \mathbb{C}^2$. Since a sandwiched surface singularity is rational, it is characterized by its dual resolution graph:
\begin{definition}[Spivakovsky \cite{S}]
  A weighted graph is \emph{sandwiched} if the graph contracts to a smooth point by properly adding $(-1$)-nodes and contracting them.
\end{definition}
\begin{example}
  Consider the following weighted graph.
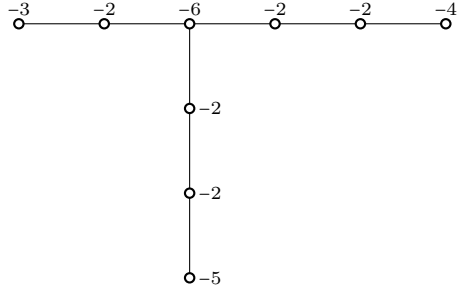
\begin{center}
  \begin{tikzpicture}
    [inner sep=1mm,
  R/.style={circle,draw=black!255,fill=white!20,thick},
  T/.style={rectangle,draw=black!255,fill=white!255,thick}]
      \node[R] (a) [label=above:$-3$]{};
      \node[R] (b) [right=of a, label=above:$-2$] {};
      \node[R] (c) [right=of b, label=above:$-6$] {};
      \node[R] (d) [right=of c, label=above:$-2$] {};
      \node[R] (e) [right=of d, label=above:$-2$] {};
      \node[R] (f) [right=of e, label=above:$-4$] {};
      \node[R] (g) [below=of c, label=right:$-2$] {};
      \node[R] (h) [below=of g, label=right:$-2$] {};
      \node[R] (i) [below=of h, label=right:$-5$] {};
      \draw[-] (a)--(b);
      \draw[-] (b)--(c);
      \draw[-] (c)--(d);
      \draw[-] (d)--(e);
      \draw[-] (e)--(f);
      \draw[-] (c)--(g);
      \draw[-] (g)--(h);
      \draw[-] (h)--(i);
  \end{tikzpicture}  
\end{center}
If we attach two $(-1)$-nodes on the western $(-3)$-node, three $(-1)$-nodes on the eastern $(-4)$-node, four $(-1)$-nodes on the southern $(-5)$-node and two $(-1)$-nodes on the central $(-6)$-node, then the graph contracts to a smooth point.
\end{example}
In \cite{S}*{Proposition~1.11}, Spivakovsky proved that the dual resolution graph of a sandwiched surface singularity is sandwiched. And conversely, for a given sandwiched graph, there exists a sandwiched surface singularity such that its dual resolution graph is the given graph.

In a different aspect, T.de Jong and D.Van Straten show that every sandwiched surface singularity can be obtained from a plane curve singularity with weights assigned to the curves. 

For a plane curve germ $C = \bigcup C_i \subset (\mathbb{C}^2,0)$, we consider the minimal good resolution of $C$. We track the multiplicities of strict transformations of $C_i$ at infinitely near $0$ for each blow-up of $C_i$ to obtain the minimal good resolution except the final one. We denote the sum of the multiplicities by $M(i)$. We then define a decorated curve:
\begin{definition}[de Jong-van Straten {\cite{deJ}*{Definition~1.3}}]
  \emph{A decorated curve} is a pair $(C, l)$ such that
    \begin{enumerate}
      \item $C = \bigcup\limits_{i=1}^{s} C_i \subset (\mathbb{C}^2,0)$ is a plane curve singularity at the origin
      \item a function $l : T=\{1, \cdots, s\} \to \mathbb{Z}$ assigning a number $l(i)$ to $C_i$
      \item $l(i) \geq M(i)$
    \end{enumerate}
\end{definition}
The function $l$ is the information of blow-ups:
\begin{definition}[de Jong-van Straten \cite{deJ}, Definition~(1,4)] 
  Let $(C, l)$ be a decorated curve.
  \begin{enumerate}
    \item The modification $\widetilde{Z}(C,l) \to \mathbb{C}^2$ determined by $(C, l)$ is obtained from the minimal embedded resolution of $C$ by $l(i)-M(i)$ consecutive blow-ups at the $i$-th branch of $C$.
    \item The analytic space $X(C,l)$ is obtained from $\widetilde{Z}(C,l) \backslash \widetilde{C}$ by blowing down all exceptional divisors not intersecting $\widetilde{C} \subset \widetilde{Z}(C,l)$.
  \end{enumerate}
\end{definition}
If $l(i) \geq M(i)+1$, then the exceptional set not intersecting $\widetilde{C_i}$ is connected(\cite{deJ}) and therefore we get one sandwiched surface singularity by blowing down. 
\begin{example}
  Let $C$ be the ordinary cusp given by the equation $y^2 -x^3 = 0$. The following are the modifications $\widetilde{Z}(C,l)$ for $l = 1, 2, 3, 4$.
    \begin{center}
      \begin{tikzpicture}
        \draw (1,1) parabola bend (0,0) (1,-1);
        \node at (1.2,1)[] {C};
        \node at (0.5,-1.2)[] {C};
        \node at (1.5,0)[] {$\leftarrow$};

        \draw[red] (2,1) to (2,-1);
        \draw plot[smooth] coordinates {(3,1)  (2,0)  (3,-1)};
        \node at (3.2,1)[] {C};
        \node at (2.5,-1.2)[] {(C,1)};
        \node at (3.5,0)[] {$\leftarrow$};

        \draw[red] (5,1) to (5,-1);
        \draw (4,1) to (6,-1);
        \draw (6,1) to (4,-1);
        \node at (6.2,1)[] {C};
        \node at (5,-1.2)[] {(C,2)};
        \node at (6.5,0)[] {$\leftarrow$};

        \draw[red] (7.5,1) to (7.5,-1);
        \draw[blue] (7, 0.5) to (8,0.5);
        \draw[blue] (7, 0) to (8,0);
        \draw (7, -0.5) to (8,-0.5);
        \node at (8.2,-0.5)[] {C};
        \node at (7.5,-1.2)[] {(C,3)};
        \node at (8.5,0)[] {$\leftarrow$};

        \draw[blue] (9.5,1) to (9.5,-1);
        \draw[blue] (9, 0.5) to (10,0.5);
        \draw[blue] (9, 0) to (10,0);
        \draw[red] (9, -0.5) to (10.5,-0.5);
        \draw (10, -1) to (10.5,0);
        \node at (10.7,0)[] {C};
        \node at (9.5,-1.2)[] {(C,4)};

      \end{tikzpicture}
    \end{center}
  The red lines are $(-1)$-curves and the blue lines are exceptional curves will be contracted. We see that $X(C, 1)$ and $X(C, 2)$ have no singularity, $X(C,3)$ has two singularities and $X(C,4)$ has a sandwiched surface singularity. 
\end{example}
\begin{proposition}[de Jong-van Straten \cite{deJ}]
\label{pro:X(C,l)}
  Any sandwiched singularity $X$ is isomorphic to $X(C, l)$ for some decorated curves $(C, l)$.
\end{proposition}
\subsection{Picture deformations}
From another point of view, the decoration $l$ can be seen as a subscheme of points on $\widetilde{C}$. Specifically, $l(i)$ is a subscheme of the branch $\widetilde{C_i}$. Similarly, if we consider $m(i)$ as a subscheme of $\widetilde{C_i}$, then the condition $l(i) \geq m(i)$ can be interpreted as $m$ being a subset of $l$.

\begin{definition}[de Jong-van Straten {\cite{deJ}*{4.2}}]
  Let $(\Delta,0)$ be a small open ball. \emph{A one-parameter deformation $(\mathscr{C}, \mathscr{L})$ of a decorated curve $(C, l)$ over $\Delta$} consists of 
    \begin{enumerate}
      \item A $\delta$-constant deformation $\mathscr{C} \to {\Delta}$ of $C$, that is, $\delta(C_{i,t})$ is constant for all $t \in \Delta^*$.
      \item A flat deformation $\mathscr{L} \subset \widetilde{C} \times \Delta$ of the scheme $l$ with the condition $\mathscr{M} \subset \mathscr{L}$ where $\mathscr{M} = \overline{\bigcup\limits_{t \in \Delta \backslash 0}m(C_t)}$.
    \end{enumerate}
\end{definition}
Here, $$\delta(C_{i,t}) = \sum\limits_{Q}\frac{m(C_{i,t},Q)(m(C_{i,t},Q)-1)}{2}$$ where $Q$ ranges over all the points infinitely near $0$(cf.\cite{MR2290112}*{Proposition~3.34}). For convenience, we use the notation $C_i$ instead of $C_{i,t}$ 
\begin{theorem}[de Jong-van Straten \cite{deJ}*{4.4}]
  For any one-parameter deformation $(\mathscr{C}, \mathscr{L})$ of a decorated curve $(C, l)$, there is a flat one parameter deformation $\mathscr{X} \to {\Delta}$ with the property that\\
  (1) $X_0 = X(C, l)$.\\
  (2) $X_t = X(C_t, l_t)$ for all $t \in \Delta^*$.\\
  Moreover, every one parameter deformation of $X(C, l)$ is obtained in this way.
\end{theorem}
We can also describe smoothings of a sandwiched surface singularity $X(C,l)$.
\begin{definition}[de Jong-van Straten {\cite{deJ}*{4.6}}]
  A one-parameter deformation $(\mathscr{C}, \mathscr{L})$ is called a \emph{picture deformation} if for $t \in \Delta^*$, the divisor $l_t$ on $\widetilde{C}_t$ is reduced.
\end{definition}
The definition means that $\mathscr{C}$ has only ordinary $m$-tuple points. For convenience, the ordinary $1$-tuple point is called a free point, a non-singular point in the support of $\mathscr{L}$.
\begin{example}[Continued from~\ref{eg:19/11}]
\label{eg:picture deformation of 19/11}
  We consider the following sandwiched structure.
$$
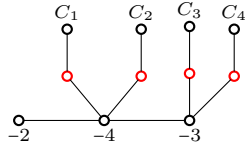
\begin{tikzpicture}
  [inner sep=1mm,
R/.style={circle,draw=black!255,fill=white!20,thick},
T/.style={circle,draw=red!255,fill=white!255,thick}]
    \node[R] (a) [label=below:$-2$]{};
    \node[R] (b) [right=of a, label=below:$-4$] {};
    \node[R] (c) [right=of b, label=below:$-3$] {};

    \node[T] (d) [above right = .5cm and .5cm of c] {};
    \node[T] (e) [above left = .5cm and .4cm of b] {};
    \node[T] (f) [above right = .5cm and .4cm of b] {};
    \node[T] (g) [above = .5cm of c] {};

    \node[R] (h) [above = .5cm of d, label = above:$C_4$] {};
    \node[R] (i) [above = .5cm of e, label = above:$C_1$] {};
    \node[R] (j) [above = .5cm of f, label = above:$C_2$] {};
    \node[R] (k) [above = .5cm of g, label = above:$C_3$] {};

    \draw[-] (a)--(b);
    \draw[-] (b)--(c);

    \draw[-] (c)--(d);
    \draw[-] (b)--(e);
    \draw[-] (b)--(f);
    \draw[-] (c)--(g);

    \draw[-] (d)--(h);
    \draw[-] (e)--(i);
    \draw[-] (f)--(j);
    \draw[-] (g)--(k);
\end{tikzpicture}$$
After the contraction, we obtain the following decorated curve.
$$
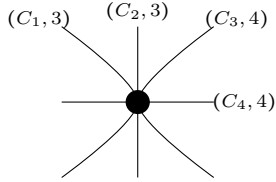
\begin{tikzpicture}
        \draw[black] (0,1) to (0,-1);
        \draw plot[smooth] coordinates {(1,1)  (0,0)  (1,-1)};
        \draw plot[smooth] coordinates {(-1,1)  (0,0)  (-1,-1)};
        \draw[black] (-1,0) to (1,0);
        \node at (-1.3,1.1)[] {$(C_1, 3)$};
        \node at (1.3,1.1)[] {$(C_3, 4)$};
        \node at (0,1.2)[] {$(C_2, 3)$};
        \node at (1.4,0)[] {$(C_4, 4)$};
        \node[circle, fill=black, inner sep=1mm] at (0,0) {}; 
\end{tikzpicture}$$
We obtain three picture deformations:
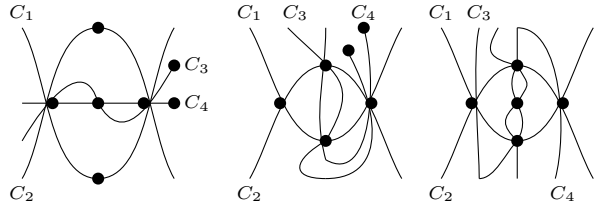
\begin{center}
\begin{tikzpicture}
        \node  (0) at (-1, 0) {};
        \node  (1) at (1, 0) {};
        \node  (2) at (-1, -1) {};
        \node  (3) at (1, -1) {};
        \node  (4) at (-1, 1) {};
        \node  (5) at (0, 1) {};
        \node  (6) at (1, 1) {};
        \node  (7) at (0, 0) {};
        \node  (8) at (0, -1) {};
        \node  (9) at (-1, -0.5) {};
        \node  (10) at (-1, 0.5) {};
        \node  (11) at (1, 0.5) {};
        \node  (12) at (1, 0.5) {};

        \node[circle, fill=black] at (-0.6,0){};
        \node[circle, fill=black] at (0,0){};
        \node[circle, fill=black] at (0.6,0){};
        \node[circle, fill=black] at (0,1){};
        \node[circle, fill=black] at (0,-1){};
        \node[circle, fill=black] at (1,0){};
        \node[circle, fill=black] at (1,0.5){};

        \node at (-1,1.2)[] {$C_1$};
        \node at (-1,-1.2)[] {$C_2$};
        \node at (1.3, 0.5)[] {$C_3$};
        \node at (1.3, 0)[] {$C_4$};

        \draw [in=-180, out=60, looseness=0.75] (2.center) to (5.center);
        \draw [in=120, out=0, looseness=0.75] (5.center) to (3.center);
        \draw [in=180, out=-60, looseness=0.75] (4.center) to (8.center);
        \draw [in=-120, out=0, looseness=0.75] (8.center) to (6.center);
        \draw (0.center) to (1.center);
        \draw [in=105, out=60, looseness=1.50] (9.center) to (7.center);
        \draw [in=-120, out=-60, looseness=1.50] (7.center) to (12.center);
\end{tikzpicture}
\begin{tikzpicture}
        \node  (2) at (-1, -1) {};
        \node  (3) at (1, -1) {};
        \node  (4) at (-1, 1) {};
        \node  (5) at (0, 0.5) {};
        \node  (6) at (1, 1) {};
        \node  (8) at (0, -0.5) {};
        \node  (9) at (0.5, 1) {};
        \node  (10) at (0, -0.75) {};
        \node  (11) at (0, -1) {};
        \node  (12) at (0, 1) {};
        \node  (13) at (0.25, 0.75) {};
        \node  (14) at (-0.5, 1) {};

        \node[circle, fill=black] at (-0.6,0){};
        \node[circle, fill=black] at (0,0.5){};
        \node[circle, fill=black] at (0.6,0){};
        \node[circle, fill=black] at (0,-0.5){};
        \node[circle, fill=black] at (0.5,1){};
        \node[circle, fill=black] at (0.3,0.7){};

        \node at (-1,1.2)[] {$C_1$};
        \node at (-1,-1.2)[] {$C_2$};
        \node at (-0.4, 1.2)[] {$C_3$};
        \node at (0.5, 1.2)[] {$C_4$};

        \draw [in=-180, out=60, looseness=0.75] (2.center) to (5.center);
        \draw [in=120, out=0, looseness=0.75] (5.center) to (3.center);
        \draw [in=180, out=-60, looseness=0.75] (4.center) to (8.center);
        \draw [in=-120, out=0, looseness=0.75] (8.center) to (6.center);
        \draw [in=-30, out=-75] (9.center) to (10.center);
        \draw [in=270, out=105] (10.center) to (12.center);
        \draw [in=0, out=-60, looseness=1.25] (13.center) to (11.center);
        \draw [in=-150, out=180, looseness=2.50] (11.center) to (8.center);
        \draw [in=-45, out=60, looseness=1.50] (8.center) to (14.center);
\end{tikzpicture}
\begin{tikzpicture}
        \node  (2) at (-1, -1) {};
        \node  (3) at (1, -1) {};
        \node  (4) at (-1, 1) {};
        \node  (5) at (0, 0.5) {};
        \node  (6) at (1, 1) {};
        \node  (8) at (0, -0.5) {};
        \node  (9) at (0, -1) {};
        \node  (10) at (0, 1) {};
        \node  (11) at (0.5, -1) {};
        \node  (12) at (-0.5, 1) {};
        \node  (13) at (-0.5, -1) {};
        \node  (14) at (-0.25, 1) {};
        \node  (15) at (0, 0) {};

        \node[circle, fill=black] at (-0.6,0){};
        \node[circle, fill=black] at (0,0.5){};
        \node[circle, fill=black] at (0.6,0){};
        \node[circle, fill=black] at (0,-0.5){};
        \node[circle, fill=black] at (0,0){};

        \node at (-1,1.2)[] {$C_1$};
        \node at (-1,-1.2)[] {$C_2$};
        \node at (-0.5, 1.2)[] {$C_3$};
        \node at (0.6, -1.2)[] {$C_4$};

        \draw [in=-180, out=60, looseness=0.75] (2.center) to (5.center);
        \draw [in=120, out=0, looseness=0.75] (5.center) to (3.center);
        \draw [in=180, out=-60, looseness=0.75] (4.center) to (8.center);
        \draw [in=-120, out=0, looseness=0.75] (8.center) to (6.center);
        \draw [in=75, out=0, looseness=0.75] (10.center) to (11.center);
        \draw (10.center) to (5.center);
        \draw (8.center) to (9.center);
        \draw [in=90, out=-105, looseness=0.50] (12.center) to (13.center);
        \draw [in=225, out=0, looseness=0.75] (13.center) to (8.center);
        \draw [in=-120, out=60, looseness=1.25] (8.center) to (5.center);
        \draw [in=-120, out=165, looseness=1.75] (5.center) to (14.center);
        \draw [bend left, looseness=1.75] (5.center) to (15.center);
        \draw [bend right=45, looseness=1.50] (15.center) to (8.center);
\end{tikzpicture}
\end{center}
\end{example}
\begin{theorem}[de Jong-van Straten \cite{deJ}*{Lemma 4.7}]
  A generic smoothing of $X(C,l)$ is realized by a picture deformation of $(C, l)$.
\end{theorem}
\subsection{Incidence matrices}
A picture deformation has a combinatorial aspect.
\begin{definition}[de Jong-van Straten {\cite{deJ}*{p.483}}] 
  The \emph{incidence matrix} of a picture deformation $(\mathscr{C}, \mathscr{L})$ is the matrix $I(\mathscr{C}, \mathscr{L}) \in M_{s, n}(\mathbb{Z})$ where $I(\mathscr{C}, \mathscr{L})_{i,j}$ is the multiplicity of $C_{i}$ at $P_j$.
\end{definition}
According to Konrad M\"{o}hring(\cite{M}), a general fiber $X(C_t, l_t)$ of a smoothing of $X(C,l)$ is blowing-ups of $\mathscr{C}$ along the support of $\mathscr{L}$. Thus, an incidence matrix encodes the intersection relations of $(-1)$-curves on the Milnor fiber.

From the $\delta$-constancy of $\mathscr{C}$ and the flatness of $\mathscr{L}$, we can formulate the necessary condition of the incidence matrices.
\begin{definition}[de Jong-van Straten {\cite{deJ}*{4,12}}]
  \emph{A combinatorial incidence matrix} of a sandwiched surface singularity $X(C, l)$ is a matrix $M = (m_{ij})_{s \times r}$ satisfying the following equations.
      \begin{equation}
      \label{eq:equationsofcombinatorialincidencematrix}
      \begin{aligned}
       \sum\limits_{j=1}^{r} \frac{m_{ij}(m_{ij}-1)}{2} = \delta(C_{i}) \text{ for all } i \\
       \sum\limits_{j=1}^{r} m_{ij}m_{kj} = C_{i}.C_{k} \text{ for all } i \neq k \\
       \sum\limits_{j=1}^{r} m_{ij} = l(i) \text{ for all } i 
      \end{aligned}
      \end{equation}
\end{definition}
Every incidence matrix satisfies Equation~\ref{eq:equationsofcombinatorialincidencematrix}. 
\begin{definition}[Park-Shin{\cite{park2022deformations}*{Definition~2.19}}]
  Let $\mathscr{C}(X)$ be the set of irreducible components of the reduce versal deformation space $\Def(X)$ and let $\mathscr{I}(X)$ be the set of all incidence matrices of $X$ of a given sandwiched structure. The \emph{incidence map} of $X$ is a map
  $$\phi_I : \mathscr{C}(X) \to \mathscr{I}(X)$$
  where, for each $S \in \mathscr{C}(X)$, $\phi_I(S)$ is defined by the incidence matrix corresponding to a picture deformation that parametrizes $S$.
\end{definition}
\section{Deformations of cyclic quotient surface singularities}
In this section, we briefly review deformation theories of cyclic quotient surface singularities. And then we analyze (combinatorial) incidence matrices of cyclic quotient surface singularities.

A cyclic quotient surface singularity $X$ of type $\frac{1}{n}(1, q)$ is a quotient surface singularity $\mathbb{C}^2/G$ where $G = \left\langle \begin{pmatrix} \zeta & 0 \\ 0 & \zeta^{q} \end{pmatrix} \right\rangle$, $\zeta$ is a primitive $n$-th root of unity and $1 \leq q < n$. 

It is well known that the minimal resolution of a cyclic quotient surface singularity of type $\frac{1}{n}(1, q)$ is a chain of $\mathbb{CP}^1$'s of self-intersection numbers $-a_1, \cdots, -a_r$ where $a_1, \dots, a_r$ are Hirzebruch-Jung continued fraction of $\frac{n}{q} = a_1 - \frac{1}{a_2 -\frac{1}{\ddots - \frac{1}{a_{r}}}}$.

We use a dual resolution graph 
$$\begin{tikzpicture}[scale=0.5]
 \node[circle] (20) at (2,0) [label=above:$A_1$,label=below:{$-a_1$}] {};

\node[empty] (250) at (3,0) [] {};
\node[empty] (30) at (3.5,0) [] {};

 \node[circle] (350) at (4.5,0) [label=above:$A_r$,label=below:{$-a_r$}] {};

\draw [-] (20)--(250);
\draw [dotted] (20)--(350);
\draw [-] (30)--(350);
\end{tikzpicture}$$
or a continued fraction $[a_1, \cdots, a_r]$
to denote a cyclic quotient surface singularity.
. The node means a curve $A_i$ with self-intersection number $-a_i$. The line means intersection relations of curves.

\begin{example}
Consider a cyclic quotient surface singularity of type $\frac{1}{19}(1,11)$. Since $\frac{19}{11} = [2, 4, 3]$, it has the dual resolution graph.
\begin{tikzpicture}
  [inner sep=1mm,
R/.style={circle,draw=black!255,fill=white!20,thick},
T/.style={circle,draw=black!255,fill=black!255,thick}]
    \node[R] (a) [label=below:$-2$, , label=above:$A_1$]{};
    \node[R] (b) [right=of a, label=below:$-4$, , label=above:$A_2$] {};
    \node[R] (c) [right=of b, label=below:$-3$, , label=above:$A_3$] {};

    \draw[-] (a)--(b);
    \draw[-] (b)--(c);
\end{tikzpicture}
\label{eg:19/11}
\end{example}

\subsection{P-resolutions}
We summarize a P-resolution description of $\text{Def}(X)$ introduced by J. Koll\'{a}r and N. I. Shepherd-Barron.

\begin{definition}[K-SB \cite{KSB}*{Definition 3.7}]
  A \emph{singularity of class $T$} is a cyclic quotient surface singularity of type $\dfrac{1}{dn^2}(1, dna-1)$ with $d, a \ge 1$, $n \ge 2 $, $(n,a)=1$. A \emph{Wahl singularity} is a singularity of class $T$ with $d=1$.
\end{definition}

Singularity of class $T$ is completely determined by its dual resolution graph.

\begin{proposition}[K-SB \cite{KSB}*{proposition 3.11}]\;
\label{prop:inductive way of singularities of class T}
\begin{enumerate}[(i)]
\item The singularities
\begin{tikzpicture}
  [inner sep=1mm,
R/.style={circle,draw=black!255,fill=white!20,thick},
T/.style={rectangle,draw=black!255,fill=white!255,thick}]
 \node[T] (10) at (1,0) [label=above:{$-4$}] {};
\end{tikzpicture}
and
\begin{tikzpicture}
  [inner sep=1mm,
R/.style={circle,draw=black!255,fill=white!20,thick},
T/.style={rectangle,draw=black!255,fill=white!255,thick}]
 \node[T] (10) at (1,0) [label=above:{$-3$}] {};
 \node[T] (20) at (2,0) [label=above:{$-2$}] {};

\node[empty] (250) at (2.5,0) [] {};
\node[empty] (30) at (3,0) [] {};

 \node[T] (350) at (3.5,0) [label=above:{$-2$}] {};
 \node[T] (450) at (4.5,0) [label=above:{$-3$}] {};

\draw [-] (10)--(20);
\draw [-] (20)--(250);
\draw [dotted] (20)--(350);
\draw [-] (30)--(350);
\draw [-] (350)--(450);
\end{tikzpicture}
are of class $T$

\item If the singularity
\begin{tikzpicture}
  [inner sep=1mm,
R/.style={circle,draw=black!255,fill=white!20,thick},
T/.style={rectangle,draw=black!255,fill=white!255,thick}]
 \node[T] (20) at (2,0) [label=above:{$-b_1$}] {};

\node[empty] (250) at (2.5,0) [] {};
\node[empty] (30) at (3,0) [] {};

 \node[T] (350) at (3.5,0) [label=above:{$-b_r$}] {};

\draw [-] (20)--(250);
\draw [dotted] (20)--(350);
\draw [-] (30)--(350);
\end{tikzpicture}
is of class $T$, then so are
\begin{equation*}
\begin{tikzpicture}
  [inner sep=1mm,
R/.style={circle,draw=black!255,fill=white!20,thick},
T/.style={rectangle,draw=black!255,fill=white!255,thick}]
 \node[T] (10) at (1,0) [label=above:{$-2$}] {};
 \node[T] (20) at (2,0) [label=above:{$-b_1$}] {};

\node[empty] (250) at (2.5,0) [] {};
\node[empty] (30) at (3,0) [] {};

 \node[T] (350) at (3.5,0) [label=above:{$-b_{r-1}$}] {};
 \node[T] (450) at (4.5,0) [label=above:{$-b_r-1$}] {};

\draw [-] (10)--(20);
\draw [-] (20)--(250);
\draw [dotted] (20)--(350);
\draw [-] (30)--(350);
\draw [-] (350)--(450);
\end{tikzpicture}
\text{~and~}
\begin{tikzpicture}  
[inner sep=1mm,
R/.style={circle,draw=black!255,fill=white!20,thick},
T/.style={rectangle,draw=black!255,fill=white!255,thick}]
 \node[T] (10) at (1,0) [label=above:{$-b_1-1$}] {};
 \node[T] (20) at (2,0) [label=above:{$-b_2$}] {};

\node[empty] (250) at (2.5,0) [] {};
\node[empty] (30) at (3,0) [] {};

 \node[T] (350) at (3.5,0) [label=above:{$-b_r$}] {};
 \node[T] (450) at (4.5,0) [label=above:{$-2$}] {};

\draw [-] (10)--(20);
\draw [-] (20)--(250);
\draw [dotted] (20)--(350);
\draw [-] (30)--(350);
\draw [-] (350)--(450);
\end{tikzpicture}
\end{equation*}

\item Every singularity of class $T$ that is not a rational double point can be obtained by starting with one of the singularities described in (i) and iterating the steps described in (ii).
\end{enumerate}
\end{proposition}
A Wahl singularity is a singularity of clss $T$ that obtained from 
\begin{tikzpicture}
  [inner sep=1mm,
R/.style={circle,draw=black!255,fill=white!20,thick},
T/.style={rectangle,draw=black!255,fill=white!255,thick}]
 \node[T] (10) at (1,0) [label=above:{$-4$}] {};
\end{tikzpicture}.
In a dual resolution graph, we use rectangle nodes to indicate exceptional curves that contracted to a singularity of class T.

\begin{definition}[K-SB \cite{KSB}*{Definition 3.12}]
\label{definition:P-resolution}
A \emph{$P$-resolution} $f \colon Y \rightarrow X$ of a quotient singularity $X$ is a modification such that $Y$ has at most rational double points or singularities of class $T$ as singularities, and $K_Y$ is $f$-relative ample, i.e., $K_{Y} \cdot E_i >0$ for all exceptional divisors $E_i$ of $f$.
\end{definition}
\begin{example}[Continued from \ref{eg:19/11}]
\label{P-resolution of 19/11}
For $X$, we have three $P$-resolutions.\\
\begin{center}
\begin{tikzpicture}
  [inner sep=1mm,
R/.style={circle,draw=black!255,fill=white!20,thick},
T/.style={rectangle,draw=black!255,fill=white!255,thick}]
    \node[R] (a) [label=below:$-2$]{};
    \node[R] (b) [right=of a, label=below:$-4$] {};
    \node[R] (c) [right=of b, label=below:$-3$] {};

    \draw[-] (a)--(b);
    \draw[-] (b)--(c);
\end{tikzpicture}
\;
\begin{tikzpicture}
  [inner sep=1mm,
R/.style={circle,draw=black!255,fill=white!20,thick},
T/.style={rectangle,draw=black!255,fill=white!255,thick}]
    \node[R] (a) [label=below:$-2$]{};
    \node[T] (b) [right=of a, label=below:$-4$] {};
    \node[R] (c) [right=of b, label=below:$-3$] {};

    \draw[-] (a)--(b);
    \draw[-] (b)--(c);
\end{tikzpicture}
\;
\begin{tikzpicture}
  [inner sep=1mm,
R/.style={circle,draw=black!255,fill=white!20,thick},
T/.style={rectangle,draw=black!255,fill=white!255,thick}]
    \node[T] (a) [label=below:$-2$]{};
    \node[T] (b) [right=of a, label=below:$-5$] {};
    \node[R] (c) [right=of b, label=below:$-1$] {};
    \node[T] (d) [right=of c, label=below:$-4$] {};

    \draw[-] (a)--(b);
    \draw[-] (b)--(c);
    \draw[-] (c)--(d);
\end{tikzpicture}
\end{center}
\label{eg:presolutionof19/11}
\end{example}

We can narrow singularities of class T down to Wahl singularities.

\begin{definition}[Behnke-Christophersen \cite{10.2307/2375004}*{p.882}]
An \emph{$M$-resolution} of a quotient surface singularity $X$ is a partial resolution $f : Y_M \rightarrow X$ such that\\
(1) $Y_M$ has only Wahl singularities.\\
(2) $K_{Y_M}$ is nef relative to $f$, i.e., $K_{Y_M}.E \geq 0$ for all $f$-exceptional curves $E$.
\end{definition}

\begin{theorem}[Behnke-Christophersen \cite{10.2307/2375004}]
Let $(X, 0)$ be a quotient surface singularity. Then\\
(1) Each $P$-resolution $Y \rightarrow X$ is dominated by a unique $M$-resolution $Y_M \rightarrow X$, i.e., there is a surjection $g : Y_M \rightarrow Y$, with the property $K_{Y_M} = g^*K_Y$.\\
(2) There is a surjective map $\text{Def}^{\text{QG}}(Y_M) \rightarrow \text{Def}^{\text{QG}}(Y)$ induced by blowing down deformations.\\
(3) There is a one-to-one correspondence between the components of $\text{Def}(X)$ and $M$-resolution of $X$.
\end{theorem}

\subsection{Stevens's description}
We recall a description of irreducible components of $\text{Def}(X)$ : Positive integer sequence $\underline{k}$(\cite{MR1129040}) by J. Stevens.

\begin{definition}[Orlik-Wagreich \cite{Or}]
We call a positive integer sequence $\underline{k} =p (k_1, \cdots, k_s) \in \mathbb{N}^r$ of length $r$ is admissible if the matrix 
$$M(\underline{k}) = \begin{bmatrix} k_1 & 1 & & & \\ 1 & k_2 & 1 & & \\ & 1 & \ddots & \ddots & \\ & & \ddots & k_{r-1} & 1 \\ & & & 1 & k_s \end{bmatrix}$$
 is positive semi-definite with $\text{rank}(M(\underline{k})) \geq s-1$.
\end{definition}
We denote the set of admissible sequence of length $r$ as $\text{Add}_s$

\begin{definition}[Christophersen \cite{Ch}] For $s \geq 1$, we define a following set.
$$K_s = \{(k_1, \cdots, k_s) \in \text{Add}_s ~|~ [k_1, \cdots, k_s] = 0\}$$
That is, the set of all admissible integer sequence of length $s$ that representing zero as the Hirzebruch-Jung continued fraction.
\end{definition}

\begin{proposition}[Stevens \cite{MR1129040}*{Theorem~4.1}]
Let $X$ be a cyclic quotient surface singularity $\frac{1}{n}(1,a)$ with $(n, a) = 1$. Let $n/(n-a) = [b_1, \cdots, b_s]$ be the Hirzebruch-Jung continued fraction. Then $K_s(n/(n-a)) = \{\underline{k} \in K_s ~|~ k_i \leq b_i\}$ parametrizes irreducible components of $\text{Def}(X)$. Therefore $\underline{k}$ corresponds to $P$-resolutions.
\end{proposition}

\begin{example}[Continued from \ref{eg:19/11}] For the cyclic quotient surface singularity of $\frac{1}{19}(1,11)$,
$K_4(19/19-11) = \{(1,2,2,1), (3,1,2,2), (2,1,3,1)\}$.
\label{eg:kof19/11}
\end{example}

Moreover, there is a geometric way to parametrize the set $K_s(n/(n-a))$.

\begin{proposition}[Stevens \cite{MR1129040}*{6.1}]
Let $\mathcal{P}_{s+1}$ be a convex $(s+1)$-gon such that each vertex is named by $b_i$ consecutively in a counterclockwise direction.(there is one unnamed vertex between the vertex $b_1$ and $b_s$) Let $\mathcal{T}(\mathcal{P}_{s+1})$ be the set of triangulations of $\mathcal{P}_{s+1}$. Then there is a bijective map from $\mathcal{T}(\mathcal{P}_{s+1})$ to $K_s$ that assigning $\theta \in \mathcal{T}(\mathcal{P}_{s+1})$ to $(k_1, \cdots, k_s)$ where $k_i$ is the number of the triangles in $\theta$ containing the vertex $b_i$.
\label{Nemethi}
\end{proposition}

\begin{example}[Continued from \ref{eg:kof19/11}]\;\\
\label{eg:triof19/11}
We have the convex $5$-gon whose vertices are named by $3, 2, 3, 2$ counterclockwise as follows.
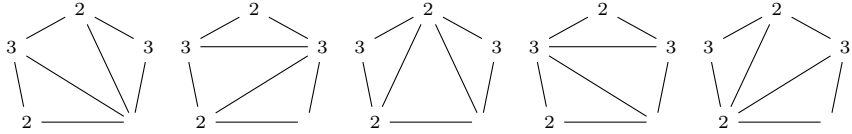
\begin{center}
\begin{tikzpicture}
\node[] (a) at (0,1) [] {2};
\node[] (b) at (-0.9,0.5) [] {3};
\node[] (c) at (-.7,-0.5) [] {2};
\node[] (d) at (.7,-0.5) [] {};
\node[] (e) at (0.9,0.5) [] {3};

\draw [-] (a) -- (b);
\draw [-] (b) -- (c);
\draw [-] (c) -- (d);
\draw [-] (d) -- (e);
\draw [-] (e) -- (a);

\draw [-] (a) -- (d);
\draw [-] (b) -- (d);
\end{tikzpicture}
\begin{tikzpicture}
\node[] (a) at (0,1) [] {2};
\node[] (b) at (-0.9,0.5) [] {3};
\node[] (c) at (-.7,-0.5) [] {2};
\node[] (d) at (.7,-0.5) [] {};
\node[] (e) at (0.9,0.5) [] {3};

\draw [-] (a) -- (b);
\draw [-] (b) -- (c);
\draw [-] (c) -- (d);
\draw [-] (d) -- (e);
\draw [-] (e) -- (a);

\draw [-] (b) -- (e);
\draw [-] (c) -- (e);
\end{tikzpicture}
\begin{tikzpicture}
\node[] (a) at (0,1) [] {2};
\node[] (b) at (-0.9,0.5) [] {3};
\node[] (c) at (-.7,-0.5) [] {2};
\node[] (d) at (.7,-0.5) [] {};
\node[] (e) at (0.9,0.5) [] {3};

\draw [-] (a) -- (b);
\draw [-] (b) -- (c);
\draw [-] (c) -- (d);
\draw [-] (d) -- (e);
\draw [-] (e) -- (a);

\draw [-] (a) -- (c);
\draw [-] (a) -- (d);
\end{tikzpicture}
\begin{tikzpicture}
\node[] (a) at (0,1) [] {2};
\node[] (b) at (-0.9,0.5) [] {3};
\node[] (c) at (-.7,-0.5) [] {2};
\node[] (d) at (.7,-0.5) [] {};
\node[] (e) at (0.9,0.5) [] {3};

\draw [-] (a) -- (b);
\draw [-] (b) -- (c);
\draw [-] (c) -- (d);
\draw [-] (d) -- (e);
\draw [-] (e) -- (a);

\draw [-] (b) -- (e);
\draw [-] (b) -- (d);
\end{tikzpicture}
\begin{tikzpicture}
\node[] (a) at (0,1) [] {2};
\node[] (b) at (-0.9,0.5) [] {3};
\node[] (c) at (-.7,-0.5) [] {2};
\node[] (d) at (.7,-0.5) [] {};
\node[] (e) at (0.9,0.5) [] {3};

\draw [-] (a) -- (b);
\draw [-] (b) -- (c);
\draw [-] (c) -- (d);
\draw [-] (d) -- (e);
\draw [-] (e) -- (a);

\draw [-] (c) -- (a);
\draw [-] (c) -- (e);
\end{tikzpicture}
\end{center}
From these, we obtain five integer sequences $(1, 2, 2, 1)$, $(3,1,2,2)$, $(1, 3, 1, 2)$, $(2, 1, 3, 1)$, $(2, 2, 1, 3)$. Since $k_i < b_i$, sequences that we want are $(1,2,2,1)$, $(3,1,2,2)$, $(2,1,3,1)$ and we check that these are the same with Example~\ref{eg:kof19/11}.
\end{example}

\subsection{Stevens to Incidence matrix}
From a sequence $\underline{k} \in K_r(n/(n-a))$ and its triangulation $\theta$, we can construct an incidence matrix and this incidence matrix corresponds to the P-resolution that parametrized by $\underline{k}$.

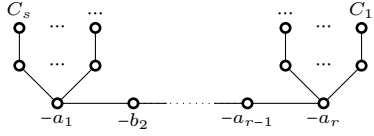
\begin{figure}
\centering
\begin{tikzpicture}
 \node[circle] (10) at (1,0) [label=below:{$-a_1$}] {};
 \node[circle] (20) at (2,0) [label=below:{$-b_2$}] {};

\node[empty] (250) at (2.5,0) [] {};
\node[empty] (30) at (3,0) [] {};

 \node[circle] (350) at (3.5,0) [label=below:{$-a_{r-1}$}] {};
 \node[circle] (450) at (4.5,0) [label=below:{$-a_r$}] {};

 \node[circle] (c1) at (.5,.5)  {};
 \node[diamond] (c2) at (1,.5)  {$\cdots$};
 \node[circle] (c3) at (1.5,.5) {};

 \node[circle] (d1) at (.5,1) [label=above:{$C_s$}] {};
 \node[diamond] (d2) at (1,1)  {$\cdots$};
 \node[circle] (d3) at (1.5,1) [label=above:{$\cdots$}] {}; 

 \node[circle] (c4) at (4,.5)  {};
 \node[diamond] (c5) at (4.5,.5)  {$\cdots$};
 \node[circle] (c6) at (5,.5)  {};

 \node[circle] (d4) at (4,1) [label=above:{$\cdots$}] {};
 \node[diamond] (d5) at (4.5,1)  {$\cdots$};
 \node[circle] (d6) at (5,1) [label=above:{$C_1$}] {}; 

\draw [-] (10)--(20);
\draw [-] (20)--(250);
\draw [dotted] (20)--(350);
\draw [-] (30)--(350);
\draw [-] (350)--(450);

\draw [-] (10)--(c1);
\draw [-] (c1)--(d1);
\draw [-] (10)--(c3);
\draw [-] (c3)--(d3);

\draw [-] (450)--(c4);
\draw [-] (c4)--(d4);
\draw [-] (450)--(c6);
\draw [-] (c6)--(d6);
\end{tikzpicture}
\caption{Usual sandwich structure of a cyclic quotient surface singularity}
\label{fig:Usual Sandwiched structure of a cyclic quotient surface singularity}
\end{figure}

For the minimal resolution of a cyclic quotient surface singularity
\begin{tikzpicture}[scale=0.5]
 \node[circle] (20) at (2,0) [label=above:$A_1$,label=below:{$-a_1$}] {};

\node[empty] (250) at (2.5,0) [] {};
\node[empty] (30) at (3,0) [] {};

 \node[circle] (350) at (3.5,0) [label=above:$A_r$,label=below:{$-a_r$}] {};

\draw [-] (20)--(250);
\draw [dotted] (20)--(350);
\draw [-] (30)--(350);
\end{tikzpicture}
, we attach $(a_1-1)$ $(-1)$-curves to $A_1$, $(a_i-2)$ $(-1)$-curves to $A_i$ for $2 \leq i \leq r$. Then the graph contracts from the left to right and we obtain a sandwiched structure by attaching decorated curves on the $(-1)$-curves(refer Figure~\ref{fig:Usual Sandwiched structure of a cyclic quotient surface singularity}).

Order the triangles in $\theta = \{\triangle_j\}_{j=1}^t$. Then each triangle $\triangle_j$ contains two or three named vertices of the convex $(s+1)$-gon $\mathcal{P}_{s+1}$. We assign $\alpha(b_i, \triangle_j)$ to $0$ if $\triangle_j$ doesn't contain the vertex $b_i$, $1$ if $\triangle_j$ contains the vertex $b_i$ and $b_i$ is the first or third vertex in $\triangle_j$(according to the order of the subscripts of $b_i$ in $\triangle_j$) and $-1$ if the vertex $b_i$ is the second vertex. Then we define a matrix $D(\underline{k})$ by  $D(\underline{k})(i,j) = \alpha(b_i, \triangle_j) \in Mat_{s, t}(\mathbb{Z})$. Furthermore, define $M_{s,b_i-k_i}(i) \in Mat_{r,b_i-k_i}(\mathbb{Z})$ be the matrix that all entries of the $i$th row are $1$ and other entries are all $0$. For a matrix $M$, we define $\int M$ as the matrix that its $i$-th row is the sum from $1$st to $i$-th rows of $M$. Then we have the following theorem.

\begin{theorem}[N\'{e}methi-Popescu-Pampu \cite{NP}*{7.2}]
Define a matrix
$$D(\underline{b};\underline{k}) = (D(\underline{k}) ~|~ M_{s, b_1 - k_1}(1) ~|~ \cdots ~|~ M_{s, b_s - k_s}(s)).$$
Then the matrix $\int D(\underline{b}, \underline{k})$ is the incidence matrix that corresponding to P-resolution parametrized $\underline{k}$.  
\end{theorem}

\begin{example}[Continued from \ref{eg:triof19/11}]\;\\
\begin{minipage}{.2\textwidth}
\begin{tikzpicture}
\node[] (a) at (0,1) [] {2};
\node[] (b) at (-0.9,0.5) [] {3};
\node[] (c) at (-.7,-0.5) [] {2};
\node[] (d) at (.7,-0.5) [] {};
\node[] (e) at (0.9,0.5) [] {3};

\draw [-] (a) -- (b);
\draw [-] (b) -- (c);
\draw [-] (c) -- (d);
\draw [-] (d) -- (e);
\draw [-] (e) -- (a);

\draw [-] (a) -- (d);
\draw [-] (b) -- (d);
\end{tikzpicture}
\end{minipage}
$\rightsquigarrow$
\begin{minipage}{.4\textwidth}
$\begin{bmatrix}[ccc|cc|c|c]
~~1  & ~~0  & ~~0  & 1  & 1  & 0  & 0 \\
-1   & ~~1  & ~~0  & 0  & 0  & 0  & 0 \\
~~0  & -1   & ~~1  & 0  & 0  & 1  & 0 \\
~~0  & ~~0  & -1   & 0  & 0  & 0  & 1
\end{bmatrix}$
\end{minipage}
$\rightsquigarrow$
\begin{minipage}{.4\textwidth}
$\begin{bmatrix}[ccccccc]
1  & 0  & 0  & 1  & 1  & 0  & 0 \\
0  & 1  & 0  & 1  & 1  & 0  & 0 \\
0  & 0  & 1  & 1  & 1  & 1  & 0 \\
0  & 0  & 0  & 1  & 1  & 1  & 1
\end{bmatrix}$
\end{minipage}
\begin{minipage}{.2\textwidth}
\begin{tikzpicture}
\node[] (a) at (0,1) [] {2};
\node[] (b) at (-0.9,0.5) [] {3};
\node[] (c) at (-.7,-0.5) [] {2};
\node[] (d) at (.7,-0.5) [] {};
\node[] (e) at (0.9,0.5) [] {3};

\draw [-] (a) -- (b);
\draw [-] (b) -- (c);
\draw [-] (c) -- (d);
\draw [-] (d) -- (e);
\draw [-] (e) -- (a);

\draw [-] (b) -- (e);
\draw [-] (c) -- (e);
\end{tikzpicture}
\end{minipage}
$\rightsquigarrow$
\begin{minipage}{.4\textwidth}
$\begin{bmatrix}[ccc|c|c]
~~1  & ~~1  & ~~1  & 0  & 0   \\
-1   & ~~0  & ~~0  & 1  & 0   \\
~~1  & -1   & ~~0  & 0  & 1   \\
~~0  & ~~1  & -1   & 0  & 0 
\end{bmatrix}$
\end{minipage}
$\rightsquigarrow$
\begin{minipage}{.4\textwidth}
$\begin{bmatrix}[ccccccc]
1  & 1  & 1  & 0  & 0 \\
0  & 1  & 1  & 1  & 0 \\
1  & 0  & 1  & 1  & 1 \\
1  & 1  & 0  & 1  & 1
\end{bmatrix}$
\end{minipage}
\begin{minipage}{.2\textwidth}
\begin{tikzpicture}
\node[] (a) at (0,1) [] {2};
\node[] (b) at (-0.9,0.5) [] {3};
\node[] (c) at (-.7,-0.5) [] {2};
\node[] (d) at (.7,-0.5) [] {};
\node[] (e) at (0.9,0.5) [] {3};

\draw [-] (a) -- (b);
\draw [-] (b) -- (c);
\draw [-] (c) -- (d);
\draw [-] (d) -- (e);
\draw [-] (e) -- (a);

\draw [-] (b) -- (e);
\draw [-] (b) -- (d);
\end{tikzpicture}
\end{minipage}
$\rightsquigarrow$
\begin{minipage}{.4\textwidth}
$\begin{bmatrix}[ccc|c|c|c]
~~1  & ~~0  & ~~1  & 1  & 0 & 0  \\
-1   & ~~0  & ~~0  & 0  & 1 & 0  \\
~~1  & ~~1  & -1   & 0  & 0 & 0  \\
~~0  & -1   & ~~0  & 0  & 0 & 1
\end{bmatrix}$
\end{minipage}
$\rightsquigarrow$
\begin{minipage}{.4\textwidth}
$\begin{bmatrix}[ccccccc]
1  & 0  & 1  & 1  & 0 & 0 \\
0  & 0  & 1  & 1  & 1 & 0 \\
1  & 1  & 0  & 1  & 1 & 0 \\
1  & 0  & 0  & 1  & 1 & 1
\end{bmatrix}$
\end{minipage}
\end{example}

\subsection{P-resolutions to Incidence matrices}
\label{subsection:P-resolutions to Incidence matrices of sandwiched singularities}
Park-Shin(\cite{park2022deformations}) build an explicit algorithm to obtain the incidence matrix from a P-resolution of a sandwiched surface singularity by using the minimal model program. We summarize Sections~3, 5 and 6 of \cite{park2022deformations}.

Let $C = \bigcup C_i \subset \mathbb{C}^2$ be a decorated curve. There is a natural compactification $D = \bigcup D_i \subset \mathbb{CP}^2$ of the decorated curve $(C, l)$ where $D_i$ is a projective plane curve and $D_i \bigcap \mathbb{C}^2 = C_i$.

Just as we constructed a sandwiched surface singularity $X(C,l)$ from $C = \bigcup C_i$ in Proposition~\ref{pro:X(C,l)}, we can similarly construct a projective singular surface $Y(D,l)$ from $D = \bigcup D_i$. Then we have the following diagram:
\[
\begin{tikzcd}
(V,E) \arrow[r,hook] \arrow[d] & (W,E) \arrow[d] \\
(X,p) \arrow[r,hook]           & (Y,p)
\end{tikzcd}
\]
where $(V,E)$ and $(W,E)$ are minimal resolutions of $(X,p)$ and $(Y,p)$ respectively. Then we have:
\begin{theorem}[Park-Shin {\cite{park2022deformations}*{Theorem 3.2}}]
Any deformation of $X(C,l)$ can be extended to a deformation of $Y(D,l)$.
\end{theorem}
Therefore we work on the compactified decorated curve $(D, l)$ and the singularity $Y(D,l)$. Consider a one parameter smoothing $\mathscr{Y} \to \Delta$ of the sandwiched surface singularity $Y(D,l)$. Assume that there exists an $M$-resolution $Z \to Y$ such that the $\mathbb{Q}$-Gorenstein smoothing $\mathscr{Z} \to \Delta$ blows down to the smoothing $\mathscr{Y} \to \Delta$. To apply the minimal model program, especially flips and divisorial contraction, we consider the morphism $\mathscr{Z} \to \mathscr{Y}$ as an extremal neighborhood.
\begin{definition}[cf.{\cite{MR3606997}*{Proposition~2.1} , \cite{MR3616327}*{Definition~2.5}}]
  Let $(Q \in Y)$ be a two-dimensional germ of a cyclic quotient surface singularity, $f:Z \to Y$ be a partial resolution of $Q \in Y$ such that $f^{-1}(Q)=C$ is a smooth rational curve with one(or two) Wahl singularity(ies) of $Z$ on it. Suppose that $K_{Z}.C < 0$. Let $(Z \subset \mathscr{Z}) \to (0 \in \Delta)$ be a $\mathbb{Q}$-Gorenstein smoothing of $Z$ over small disk $\Delta$. Let $(Y \subset \mathscr{Y}) \to \Delta$ be the corresponding blow-down deformation of $Y$. The induced birational morphism $(C \subset \mathscr{Z}) \to (Q \in \mathscr{Y})$ is called an \emph{extremal neighborhood of type mk1A(or mk2A)}. It is \emph{flipping} if the exceptional set is $C$ and \emph{divisorial} if the exceptional set is of dimension $2$.
\end{definition}
\begin{proposition}[Koll\'{a}r-Mori {\cite{MR1149195}*{\S11 and Theorem~13.5}}]
\label{prop:extremal P-resolution}
  Suppose that $f : (C \subset \mathscr{Z}) \to (Q \in \mathscr{Y})$ is a flipping extremal neighborhood of type $mk1A$ or $mk2A$. Let $f_0:(C \subset Z) \to (Q \in Y)$ be the contraction of $C$ between the central fibers $Z$ and $Y$. Then there exists an extremal $P$-resolution $f^+:(C^+ \subset Z^+) \to (Q \in Y)$ such that the flip $(C^+ \subset \mathscr{Z}^+) \to (Q \in \mathscr{Y})$ is obtained by the blow-down deformation of a $\mathbb{Q}$-Gorenstein smoothing of $Z^+$. That is, we have the commutative diagram
 \begin{center} 
  \begin{tikzcd}
    (C \subset \mathscr{Z}) \arrow[rr,dashrightarrow] \arrow[rd] \arrow[rdd] &                                & (C^+ \subset \mathscr{Z}^+) \arrow[ld] \arrow[ldd] \\
                                                                             & (Q \in \mathscr{Y}), \arrow[d] & \\
                                                                             & (0 \in \Delta)                 & \\
  \end{tikzcd}
\end{center}
which is restricted to the central fibers as follows:
 \begin{center} 
  \begin{tikzcd}
    (C \subset Z) \arrow[rr,dashrightarrow] \arrow[rd]&                                & (C^+ \subset Z^+) \arrow[ld] \\
                                                      & (Q \in Y)                      & \\
   \end{tikzcd}
\end{center}
\end{proposition}
In this paper, we encounter only one type of flips. Consider an extremal neighborhood $\mathscr{Z} \supset C$ where a Wahl singularity $[a_1, \cdots, a_r]$ is on $C$ and $K_Z.C < 0$. In the minimal resolution of $Z$, the curve $C$ becomes a $(-1)$-curve. Assume that the $(-1)$-curve intersects only the exceptional curve $A_r$. Then we have:
\begin{proposition}[Urz\'{u}a {\cite{MR3616327}*{Proposition 2.15}}]
  Assume that $a_i \geq 3$ and $a_j = \cdots = a_r = 2$ for $j > i$ for some $i$. If $a_r \geq 3$, then $r = i$. Then the image of $A_1$ in the extremal $P$-resolution $Z^+$ is the curve $C^+$ and there is a Wahl singularity $[a_2, \dots, a_i -1]$ on $C^+$ if $i \geq 2$.
\end{proposition}
In our situation, a decorated curve $D_i$ intersect the curve $C$. In general, after the flip, the curve $D_i$ degenerates.
\begin{proposition}[Urz\'{u}a {\cite{MR3593546}*{Proposition~4.1}}]
  Let the image of $D_i$ be $D_i^+$ after the flip. Then $D_i^+ = D_i' + A_1$ where $D_i'$ is the strict transform of $D_i$.
\end{proposition}
\begin{example} 
Let $\mathscr{L}$ be an extremal neighborhood such that a Wahl singularity $[a_1, \cdots, a_r]$ is on a curve $C$ and a curve $D$ intersects $C$ at the different point with the singularity. After the flip, the image $C^+$ of $C$ is the curve $A_1$, and the curve $D$ degenerates to $D' + A_1$.
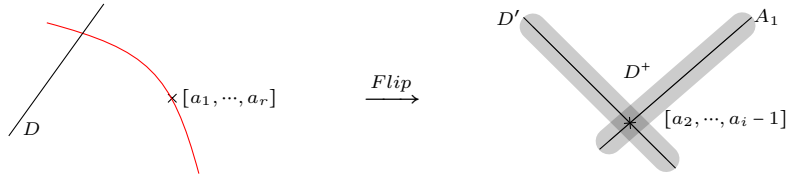
\begin{center}
  \begin{minipage}{0.4\textwidth}
    \begin{center}
    \begin{tikzpicture}
      \node  (0) at (-1, 1) {};
      \node  (1) at (1, -1) {};
      \node  (2) at (-0.25, 1.25) {};
      \node  (3) at (-1.5, -0.5) {};

      \node[cross] [label=right:{$[a_1,\cdots,a_r]$}] at (0.65,0) {};
      \node        [] at (-1.2,-0.4) {$D$};

      \draw[red] [bend left, looseness=1.25] (0.center) to (1.center);
      \draw (2.center) to (3.center);
    \end{tikzpicture}
    \end{center}
    \end{minipage}
    \begin{minipage}{.1\textwidth}
    \begin{center}
    $\xrightarrow{Flip}$
    \end{center}
    \end{minipage}
    \begin{minipage}{0.4\textwidth}
    \begin{center}
    \begin{tikzpicture}
    \node  (0) at (-1, 1) {};
    \node  (1) at (1, -1) {};
    \node  (2) at (2, 1) {};
    \node  (3) at (0, -0.75) {};

    \node[cross, rotate=50] [label={[xshift=2.2cm, yshift=-0.2cm]$[a_2,\cdots,a_i-1]$}] at (0.4,-0.4) {};
    \node        [] at (-1.2,1) {$D'$};
    \node        [] at (2.2,1) {$A_1$};
    \node        [] at (.5,0.3) {$D^+$};

    \draw (0.center) to (1.center);
    \draw (2.center) to (3.center);

    \draw[edge,line width=10pt] (0) -- (1);
    \draw[edge,line width=10pt] (2) -- (3);
    \end{tikzpicture}
    \end{center}
    \end{minipage}
\end{center}
We use the following dual resolution graph notation.
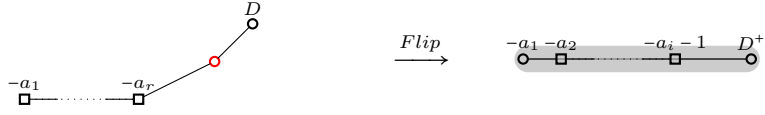
\begin{center}
  \begin{minipage}{.4\textwidth}
    \begin{tikzpicture}
  [inner sep=1mm,
R/.style={circle,draw=red!255,fill=white!20,thick},
T/.style={rectangle,draw=black!255,fill=white!255,thick},
S/.style={circle,draw=black!255,fill=white!255,thick}]

 \node[T] (20) at (2,0) [label=above:{$-a_1$}] {};

\node[empty] (250) at (2.5,0) [] {};
\node[empty] (30) at (3,0) [] {};

 \node[T] (350) at (3.5,0) [label=above:{$-a_r$}] {};
 \node[R] (450) at (4.5,.5) [] {};
 \node[S] (550) at (5,1) [label=above:{$D$}] {};

\draw [-] (20)--(250);
\draw [dotted] (20)--(350);
\draw [-] (30)--(350);
\draw [-] (350)--(450);
\draw [-] (450)--(550);
\end{tikzpicture}
  \end{minipage}
    \begin{minipage}{.1\textwidth}
    $\xrightarrow{Flip}$
    \end{minipage}
    \begin{minipage}{.3\textwidth}
    \begin{tikzpicture}
  [inner sep=1mm,
R/.style={circle,draw=red!255,fill=white!20,thick},
T/.style={rectangle,draw=black!255,fill=white!255,thick},
S/.style={circle,draw=black!255,fill=white!255,thick}]

 \node[S] (20) at (2,0) [label=above:{$-a_1$}] {};
 \node[T] (25) at (2.5,0) [label=above:{$-a_2$}] {};

\node[empty] (250) at (3,0) [] {};
\node[empty] (30) at (3.5,0) [] {};

 \node[T] (350) at (4,0) [label=above:{$-a_i-1$}] {};
 \node[S] (450) at (5,0) [label=above:{$D^+$}] {};

\draw [-] (20)--(25);
\draw [-] (25)--(250);
\draw [dotted] (25)--(350);
\draw [-] (30)--(350);
\draw [-] (350)--(450);

\draw[edge,line width=10pt] (20) -- (450);

\end{tikzpicture}
  \end{minipage}

\end{center}
\end{example}
A divisorial contraction is just a blow-down of a $(-1)$-curve in the special and general fiber of $\mathscr{Z} \to \Delta$.
\begin{proposition}[Urz\'{u}a \cite{MR3616327}]
  If an $mk1A$ or $mk2A$ is divisorial, then $(Q \in Y)$ is a Wahl singularity. In addition, the divisorial contraction $F:\mathscr{Z} \to \mathscr{Y}$ induces the blowing-down of a $(-1)$-curve between the smooth fibers of $\mathscr{Z} \to \Delta$ and $\mathscr{Y} \to \Delta$.
\end{proposition}
An incidence matrix encodes the intersection relations of $(C_{i,t},l_t)$. After blowing up at these points, the resulting object is a complement of the Milnor fiber, i.e., the general fiber of $\mathscr{Z} \to \Delta$. Therefore, if we can locate $(-1)$-curves in the general fiber, then we can induce the corresponding incidence matrix. For this, we use the following results.
\begin{definition}[Urz\'{u}a {\cite{MR3593546}*{Definition~2.1}}]
  A \emph{W-surface} is a normal projective surface $S$ with a proper deformation $\mathscr{S} \to \Delta$ such that
  \begin{enumerate}
    \item $S$ has at most singularities of class $T_0$
    \item $\mathscr{S}$ is a normal complex $3$-fold where the canonical divisor $K_{\mathscr{S}}$ is $\mathbb{Q}$-Cartier
    \item The fiber $S_0$ is reduced and isomorphic to $S$
    \item The fiber $S_t$ is nonsingular for $t \neq 0$
  \end{enumerate}
\end{definition}

\begin{proposition}[Urz\'{u}a {\cite{MR3593546}*{Corollary 3.5}}]
\label{prop:urzua}
If $S_0$ is birational to $S_t$ for $t \neq 0$, then the smoothing $\mathscr{S} \rightarrow \Delta$ can be reduced to a deformation $\mathscr{S}' \rightarrow \Delta$ whose central fiber $S_0'$ is smooth by applying a finite number of the divisorial contractions and the flips.
\end{proposition}
Assume that a compactified decorated curve $(D, l)$, its corresponding singularity $Y(D,l)$ and a $M$-resolution $Z \to Y(D,l)$ are given. The singularity $Z$ is a $W$-surface
with its smoothing $\mathscr{Z} \to \Delta$. Since $Z_0$ and $Z_t$ have a $(+1)$-curve, they are birational to $\mathbb{C}^2$. Therefore we can apply the proposition~\ref{prop:urzua} to the smoothing $\mathscr{Z} \to \Delta$
\begin{proposition}[Park-Shin {\cite{park2022deformations}*{Proposition 6.2}}] 
  By applying the divisorial contractions and flips to $(-1)$-curves on the central fiber $Z_0$ of $\mathscr{Z} \to \Delta$, one can rum MMP to $\mathscr{Z} \to \Delta$ untill one obtains a deformation $\mathscr{Z}' \to \Delta$ whose central fiber $Z_0'$ is smooth.
\end{proposition}
Steps to make the central fiber to be smooth are as follows. For a $(-1)$-curve in $Z_0$, 
\begin{enumerate}
  \item If a Wahl singularity is not on the $(-1)$-curve, then contract it(divisorial contraction).
  \item If a Wahl singularity is on the $(-1)$-curve, then apply the flip.
  \item If the Wahl singularity still remains(in fact, new Wahl singularity), there must be new $(-1)$-curve pass through the singularity. Apply the flip again.
  \item We can apply flips until the Wahl singularity disappear.
\end{enumerate}
Flips do not affect the general fiber but a divisorial contraction is just a blow-down of a $(-1)$-curve on each fibers. Therefore,
\begin{corollary}[Park-Shin {\cite{park2022deformations}*{Corollary~6.3}}] 
  In the previous proposition, a general fiber $Z_t$ of $\mathscr{Z} \to \Delta$ is obtained by blowing up several times a general fiber $Z_t'$ of the smoothing $\mathscr{Z}' \to \Delta$ for $\mathscr{Z}_0'$.
\end{corollary}
Park-Shin prove that a general fiber obtained as the previous corollary is the same with the general fiber of a generic smoothing of $X(C,l)$ comes from a picture deformation $(\mathscr{C}, \mathscr{L})$ by blowing ups.
\begin{theorem}[Park-Shin {\cite{park2022deformations}*{Theorem~6.4}}]
  One can run the semi-stable MMP to $\mathscr{Z} \to \Delta$ until one obtains the corresponding picture deformation $(\mathscr{D}, \mathscr{L})$ of the compactified decorated curve $(D, l)$.
\end{theorem}
\begin{example}[Continued from \ref{P-resolution of 19/11}]
\label{eg:MMP on 19/11}
 We find incidence matrices of the CQSS of $\frac{1}{19}(1,11)$ under the usual sandwiched structure. The left side is a $P$-resolution, and the right side is the general fiber. A divisorial contraction is shorten to be d.c\\
(1) The minimal resolution. We apply only divisorial contractions.
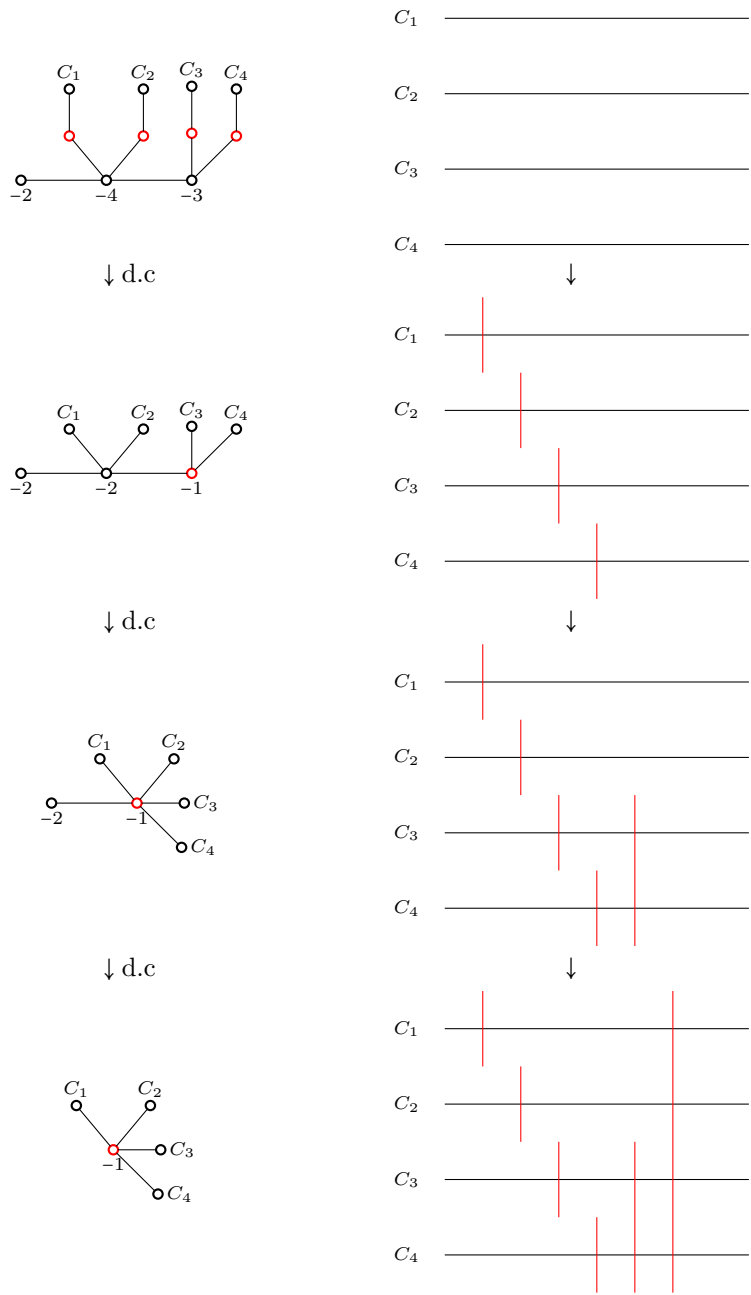
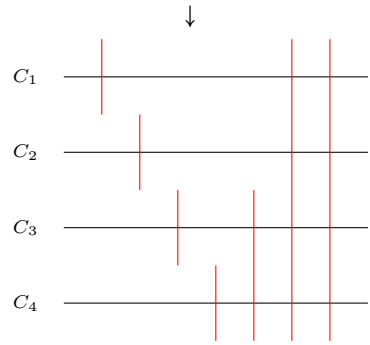
\begin{center}
  \begin{minipage}{.45\textwidth}
  \begin{center}
    \begin{tikzpicture}
        [inner sep=1mm,
      R/.style={circle,draw=black!255,fill=white!20,thick},
      T/.style={circle,draw=red!255,fill=white!255,thick}]
      \node[R] (a) [label=below:$-2$]{};
     \node[R] (b) [right=of a, label=below:$-4$] {};
     \node[R] (c) [right=of b, label=below:$-3$] {};

     \node[T] (d) [above right = .5cm and .5cm of c] {};
     \node[T] (e) [above left = .5cm and .4cm of b] {};
     \node[T] (f) [above right = .5cm and .4cm of b] {};
      \node[T] (g) [above = .5cm of c] {};

      \node[R] (h) [above = .5cm of d, label = above:$C_4$] {};
     \node[R] (i) [above = .5cm of e, label = above:$C_1$] {};
     \node[R] (j) [above = .5cm of f, label = above:$C_2$] {};
     \node[R] (k) [above = .5cm of g, label = above:$C_3$] {};

     \draw[-] (a)--(b);
      \draw[-] (b)--(c);

     \draw[-] (c)--(d);
     \draw[-] (b)--(e);
     \draw[-] (b)--(f);
     \draw[-] (c)--(g);

     \draw[-] (d)--(h);
     \draw[-] (e)--(i);
     \draw[-] (f)--(j);
      \draw[-] (g)--(k);
    \end{tikzpicture}
    \end{center}
  \end{minipage}
  \begin{minipage}{.45\textwidth}
  \begin{center}
    \begin{tikzpicture}
    \node  (0) at (-2, 1.5) {};
    \node  (1) at (2, 1.5) {};
    \node  (2) at (-2, 0.5) {};
    \node  (3) at (2, 0.5) {};
    \node  (4) at (-2, -0.5) {};
    \node  (5) at (2, -0.5) {};
    \node  (6) at (-2, -1.5) {};
    \node  (7) at (2, -1.5) {};
    \node  (8) at (-2.5, 1.5) {};
    \node  (9) at (-2.5, 1.5) {$C_1$};
    \node  (11) at (-2.5, 0.5) {$C_2$};
    \node  (12) at (-2.5, -0.5) {$C_3$};
    \node  (13) at (-2.5, -1.5) {$C_4$};

    \draw (0.center) to (1.center);
    \draw (3.center) to (2.center);
    \draw (4.center) to (5.center);
    \draw (7.center) to (6.center);
    \end{tikzpicture}
    \end{center}
  \end{minipage}
\begin{minipage}{.45\textwidth}
  \begin{center}
  $\downarrow \text{d.c} $
  \end{center}
\end{minipage}
\begin{minipage}{.45\textwidth}
  \begin{center}
  $\downarrow$
  \end{center}
\end{minipage}
\begin{minipage}{.45\textwidth}
\begin{center}
    \begin{tikzpicture}
        [inner sep=1mm,
      R/.style={circle,draw=black!255,fill=white!20,thick},
      T/.style={circle,draw=red!255,fill=white!255,thick}]
      \node[R] (a) [label=below:$-2$]{};
     \node[R] (b) [right=of a, label=below:$-2$] {};
     \node[T] (c) [right=of b, label=below:$-1$] {};

     \node[R] (d) [above right = .5cm and .5cm of c, label=above:$C_4$] {};
     \node[R] (e) [above left = .5cm and .4cm of b, label=above:$C_1$] {};
     \node[R] (f) [above right = .5cm and .4cm of b, label=above:$C_2$] {};
     \node[R] (g) [above = .5cm of c, label=above:$C_3$] {};

     \draw[-] (a)--(b);
      \draw[-] (b)--(c);

     \draw[-] (c)--(d);
     \draw[-] (b)--(e);
     \draw[-] (b)--(f);
     \draw[-] (c)--(g);
    \end{tikzpicture}
\end{center}
  \end{minipage}
  \begin{minipage}{.45\textwidth}
  \begin{center}
    \begin{tikzpicture}
    \node  (0) at (-2, 1.5) {};
    \node  (1) at (2, 1.5) {};
    \node  (2) at (-2, 0.5) {};
    \node  (3) at (2, 0.5) {};
    \node  (4) at (-2, -0.5) {};
    \node  (5) at (2, -0.5) {};
    \node  (6) at (-2, -1.5) {};
    \node  (7) at (2, -1.5) {};
    \node  (11) at (-2.5, 0.5) {$C_2$};
    \node  (12) at (-2.5, -0.5) {$C_3$};
    \node  (13) at (-2.5, -1.5) {$C_4$};
    \node  (14) at (-2.5, 1.5) {$C_1$};
    \node  (15) at (-1.5, 2) {};
    \node  (16) at (-1.5, 1) {};
    \node  (17) at (-1, 1) {};
    \node  (18) at (-1, 0) {};
    \node  (19) at (-0.5, 0) {};
    \node  (20) at (-0.5, -1) {};
    \node  (21) at (0, -1) {};
    \node  (22) at (0, -2) {};

    \draw (0.center) to (1.center);
    \draw (3.center) to (2.center);
    \draw (4.center) to (5.center);
    \draw (7.center) to (6.center);
    \draw [red] (15.center) to (16.center);
    \draw [red] (17.center) to (18.center);
    \draw [red] (19.center) to (20.center);
    \draw [red] (21.center) to (22.center);
    \end{tikzpicture}
    \end{center}
  \end{minipage}
\begin{minipage}{.45\textwidth}
  \begin{center}
    $\downarrow\text{d.c}$
  \end{center}
\end{minipage}
\begin{minipage}{.45\textwidth}
  \begin{center}
  $\downarrow$
  \end{center}
\end{minipage}
\begin{minipage}{.45\textwidth}
    \begin{center}
    \begin{tikzpicture}
        [inner sep=1mm,
      R/.style={circle,draw=black!255,fill=white!20,thick},
      T/.style={circle,draw=red!255,fill=white!255,thick}]
      \node[R] (a) [label=below:$-2$]{};
     \node[T] (b) [right=of a, label=below:$-1$] {};

     \node[R] (d) [below right = .5cm and .5cm of b, label=right:$C_4$] {};
     \node[R] (e) [above left = .5cm and .4cm of b, label=above:$C_1$] {};
     \node[R] (f) [above right = .5cm and .4cm of b, label=above:$C_2$] {};
     \node[R] (g) [right = .5cm of b, label=right:$C_3$] {};

     \draw[-] (a)--(b);
     \draw[-] (b)--(e);
     \draw[-] (b)--(f);
     \draw[-] (b)--(d);
     \draw[-] (b)--(g);
    \end{tikzpicture}
    \end{center}
  \end{minipage}
  \begin{minipage}{.45\textwidth}
  \begin{center}
    \begin{tikzpicture}
    \node  (0) at (-2, 1.5) {};
    \node  (1) at (2, 1.5) {};
    \node  (2) at (-2, 0.5) {};
    \node  (3) at (2, 0.5) {};
    \node  (4) at (-2, -0.5) {};
    \node  (5) at (2, -0.5) {};
    \node  (6) at (-2, -1.5) {};
    \node  (7) at (2, -1.5) {};
    \node  (11) at (-2.5, 0.5) {$C_2$};
    \node  (12) at (-2.5, -0.5) {$C_3$};
    \node  (13) at (-2.5, -1.5) {$C_4$};
    \node  (14) at (-2.5, 1.5) {$C_1$};
    \node  (15) at (-1.5, 2) {};
    \node  (16) at (-1.5, 1) {};
    \node  (17) at (-1, 1) {};
    \node  (18) at (-1, 0) {};
    \node  (19) at (-0.5, 0) {};
    \node  (20) at (-0.5, -1) {};
    \node  (21) at (0, -1) {};
    \node  (22) at (0, -2) {};
    \node  (23) at (0.5, 0) {};
    \node  (24) at (0.5, -2) {};

    \draw (0.center) to (1.center);
    \draw (3.center) to (2.center);
    \draw (4.center) to (5.center);
    \draw (7.center) to (6.center);
    \draw [red] (15.center) to (16.center);
    \draw [red] (17.center) to (18.center);
    \draw [red] (19.center) to (20.center);
    \draw [red] (21.center) to (22.center);
    \draw [red] (23.center) to (24.center);
    \end{tikzpicture}
    \end{center}
  \end{minipage}
\begin{minipage}{.45\textwidth}
  \begin{center}
  $\downarrow\text{d.c}$
  \end{center}
\end{minipage}
\begin{minipage}{.45\textwidth}
  \begin{center}
  $\downarrow$
  \end{center}
\end{minipage}
\begin{minipage}{.45\textwidth}
    \begin{center}
    \begin{tikzpicture}
        [inner sep=1mm,
      R/.style={circle,draw=black!255,fill=white!20,thick},
      T/.style={circle,draw=red!255,fill=white!255,thick}]
      \node[T] (a) [label=below:$-1$]{};

     \node[R] (d) [below right = .5cm and .5cm of a, label=right:$C_4$] {};
     \node[R] (e) [above left = .5cm and .4cm of a, label=above:$C_1$] {};
     \node[R] (f) [above right = .5cm and .4cm of a, label=above:$C_2$] {};
     \node[R] (g) [right = .5cm of a, label=right:$C_3$] {};

     \draw[-] (a)--(d);
     \draw[-] (a)--(e);
     \draw[-] (a)--(f);
     \draw[-] (a)--(g);
    \end{tikzpicture}
    \end{center}
  \end{minipage}
  \begin{minipage}{.45\textwidth}
  \begin{center}
    \begin{tikzpicture}
    \node  (0) at (-2, 1.5) {};
    \node  (1) at (2, 1.5) {};
    \node  (2) at (-2, 0.5) {};
    \node  (3) at (2, 0.5) {};
    \node  (4) at (-2, -0.5) {};
    \node  (5) at (2, -0.5) {};
    \node  (6) at (-2, -1.5) {};
    \node  (7) at (2, -1.5) {};
    \node  (11) at (-2.5, 0.5) {$C_2$};
    \node  (12) at (-2.5, -0.5) {$C_3$};
    \node  (13) at (-2.5, -1.5) {$C_4$};
    \node  (14) at (-2.5, 1.5) {$C_1$};
    \node  (15) at (-1.5, 2) {};
    \node  (16) at (-1.5, 1) {};
    \node  (17) at (-1, 1) {};
    \node  (18) at (-1, 0) {};
    \node  (19) at (-0.5, 0) {};
    \node  (20) at (-0.5, -1) {};
    \node  (21) at (0, -1) {};
    \node  (22) at (0, -2) {};
    \node  (23) at (0.5, 0) {};
    \node  (24) at (0.5, -2) {};
    \node  (25) at (1, 2) {};
    \node  (26) at (1, -2) {};

    \draw (0.center) to (1.center);
    \draw (3.center) to (2.center);
    \draw (4.center) to (5.center);
    \draw (7.center) to (6.center);
    \draw [red] (15.center) to (16.center);
    \draw [red] (17.center) to (18.center);
    \draw [red] (19.center) to (20.center);
    \draw [red] (21.center) to (22.center);
    \draw [red] (23.center) to (24.center);
    \draw [red] (25.center) to (26.center);
    \end{tikzpicture}
    \end{center}
  \end{minipage}
\begin{minipage}{.45\textwidth}
  \begin{center}
  $\downarrow\text{d.c}$
  \end{center}
\end{minipage}
\begin{minipage}{.45\textwidth}
  \begin{center}
  $\downarrow$
  \end{center}
\end{minipage}
\begin{minipage}{.45\textwidth}
    \;
  \end{minipage}
  \begin{minipage}{.45\textwidth}
  \begin{center}
    \begin{tikzpicture}
    \node  (0) at (-2, 1.5) {};
    \node  (1) at (2, 1.5) {};
    \node  (2) at (-2, 0.5) {};
    \node  (3) at (2, 0.5) {};
    \node  (4) at (-2, -0.5) {};
    \node  (5) at (2, -0.5) {};
    \node  (6) at (-2, -1.5) {};
    \node  (7) at (2, -1.5) {};
    \node  (11) at (-2.5, 0.5) {$C_2$};
    \node  (12) at (-2.5, -0.5) {$C_3$};
    \node  (13) at (-2.5, -1.5) {$C_4$};
    \node  (14) at (-2.5, 1.5) {$C_1$};
    \node  (15) at (-1.5, 2) {};
    \node  (16) at (-1.5, 1) {};
    \node  (17) at (-1, 1) {};
    \node  (18) at (-1, 0) {};
    \node  (19) at (-0.5, 0) {};
    \node  (20) at (-0.5, -1) {};
    \node  (21) at (0, -1) {};
    \node  (22) at (0, -2) {};
    \node  (23) at (0.5, 0) {};
    \node  (24) at (0.5, -2) {};
    \node  (25) at (1, 2) {};
    \node  (26) at (1, -2) {};
    \node  (27) at (1.5, 2) {};
    \node  (28) at (1.5, -2) {};

    \draw (0.center) to (1.center);
    \draw (3.center) to (2.center);
    \draw (4.center) to (5.center);
    \draw (7.center) to (6.center);
    \draw [red] (15.center) to (16.center);
    \draw [red] (17.center) to (18.center);
    \draw [red] (19.center) to (20.center);
    \draw [red] (21.center) to (22.center);
    \draw [red] (23.center) to (24.center);
    \draw [red] (25.center) to (26.center);
    \draw [red] (27.center) to (28.center);
    \end{tikzpicture}
    \end{center}
  \end{minipage}
\end{center}
It corresponds to the incidence matrix 
$$\begin{bmatrix}[c|cccccccc]
         C_{1} & 1 & 1 & 1 &   &   &   &     \\
         C_{2} & 1 & 1 &   & 1 &   &   &     \\
         C_{3} & 1 & 1 &   &   & 1 & 1 &     \\
         C_{4} & 1 & 1 &   &   & 1 &   & 1         
        \end{bmatrix}$$
(2) The $P$-resolution with the Wahl singularity $[4]$.
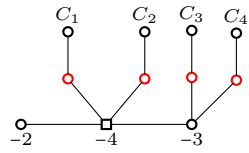
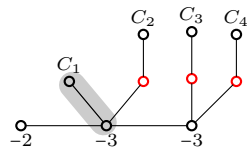
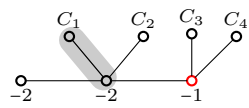
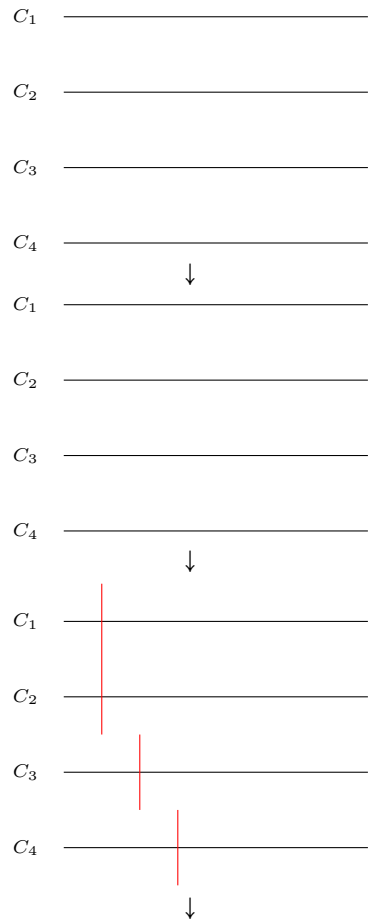
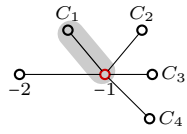
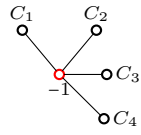
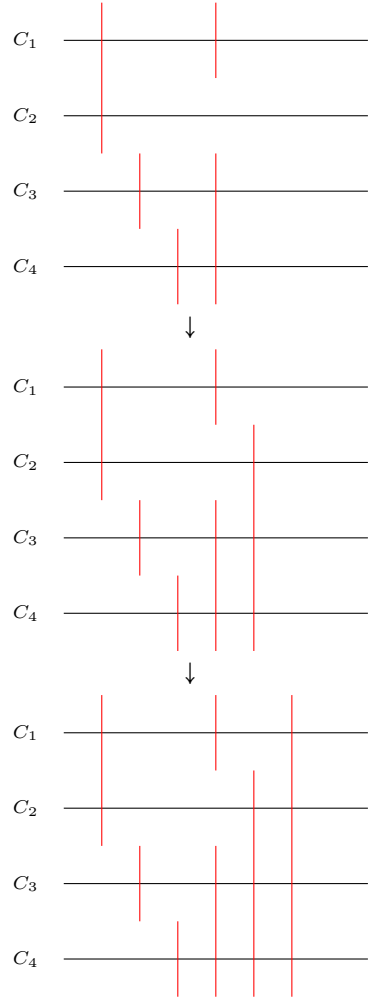
\begin{center}
  \begin{minipage}{.45\textwidth}
  \begin{center}
    \begin{tikzpicture}
        [inner sep=1mm,
      R/.style={circle,draw=black!255,fill=white!20,thick},
      T/.style={circle,draw=red!255,fill=white!255,thick},
      S/.style={rectangle,draw=black!255,fill=white!255,thick}]
     \node[R] (a) [label=below:$-2$]{};
     \node[S] (b) [right=of a, label=below:$-4$] {};
     \node[R] (c) [right=of b, label=below:$-3$] {};

     \node[T] (d) [above right = .5cm and .5cm of c] {};
     \node[T] (e) [above left = .5cm and .4cm of b] {};
     \node[T] (f) [above right = .5cm and .4cm of b] {};
      \node[T] (g) [above = .5cm of c] {};

      \node[R] (h) [above = .5cm of d, label = above:$C_4$] {};
     \node[R] (i) [above = .5cm of e, label = above:$C_1$] {};
     \node[R] (j) [above = .5cm of f, label = above:$C_2$] {};
     \node[R] (k) [above = .5cm of g, label = above:$C_3$] {};

     \draw[-] (a)--(b);
      \draw[-] (b)--(c);

     \draw[-] (c)--(d);
     \draw[-] (b)--(e);
     \draw[-] (b)--(f);
     \draw[-] (c)--(g);

     \draw[-] (d)--(h);
     \draw[-] (e)--(i);
     \draw[-] (f)--(j);
      \draw[-] (g)--(k);
    \end{tikzpicture}
    \end{center}
  \end{minipage}
  \begin{minipage}{.45\textwidth}
  \begin{center}
    \begin{tikzpicture}
    \node  (0) at (-2, 1.5) {};
    \node  (1) at (2, 1.5) {};
    \node  (2) at (-2, 0.5) {};
    \node  (3) at (2, 0.5) {};
    \node  (4) at (-2, -0.5) {};
    \node  (5) at (2, -0.5) {};
    \node  (6) at (-2, -1.5) {};
    \node  (7) at (2, -1.5) {};
    \node  (8) at (-2.5, 1.5) {};
    \node  (9) at (-2.5, 1.5) {$C_1$};
    \node  (11) at (-2.5, 0.5) {$C_2$};
    \node  (12) at (-2.5, -0.5) {$C_3$};
    \node  (13) at (-2.5, -1.5) {$C_4$};

    \draw (0.center) to (1.center);
    \draw (3.center) to (2.center);
    \draw (4.center) to (5.center);
    \draw (7.center) to (6.center);
    \end{tikzpicture}
    \end{center}
  \end{minipage}
\begin{minipage}{.45\textwidth}
  \begin{center}
  $\downarrow\text{flip}$
  \end{center}
\end{minipage}
\begin{minipage}{.45\textwidth}
  \begin{center}
  $\downarrow$
  \end{center}
\end{minipage}
\begin{minipage}{.45\textwidth}
  \begin{center}
    \begin{tikzpicture}
        [inner sep=1mm,
      R/.style={circle,draw=black!255,fill=white!20,thick},
      T/.style={circle,draw=red!255,fill=white!255,thick},
      S/.style={rectangle,draw=black!255,fill=white!255,thick}]
     \node[R] (a) [label=below:$-2$]{};
     \node[R] (b) [right=of a, label=below:$-3$] {};
     \node[R] (c) [right=of b, label=below:$-3$] {};

     \node[T] (d) [above right = .5cm and .5cm of c] {};
     \node[R] (e) [above left = .5cm and .4cm of b, label = above:$C_1$] {};
     \node[T] (f) [above right = .5cm and .4cm of b] {};
      \node[T] (g) [above = .5cm of c] {};

      \node[R] (h) [above = .5cm of d, label = above:$C_4$] {};
     \node[R] (j) [above = .5cm of f, label = above:$C_2$] {};
     \node[R] (k) [above = .5cm of g, label = above:$C_3$] {};

     \draw[-] (a)--(b);
     \draw[-] (b)--(c);
     \draw[-] (c)--(d);
     \draw[-] (b)--(e);
     \draw[-] (b)--(f);
     \draw[-] (c)--(g);
     \draw[-] (d)--(h);
     \draw[-] (f)--(j);
     \draw[-] (g)--(k);

     \draw[edge,line width=10pt] (b) -- (e);
    \end{tikzpicture}
    \end{center}
  \end{minipage}
  \begin{minipage}{.45\textwidth}
  \begin{center}
    \begin{tikzpicture}
    \node  (0) at (-2, 1.5) {};
    \node  (1) at (2, 1.5) {};
    \node  (2) at (-2, 0.5) {};
    \node  (3) at (2, 0.5) {};
    \node  (4) at (-2, -0.5) {};
    \node  (5) at (2, -0.5) {};
    \node  (6) at (-2, -1.5) {};
    \node  (7) at (2, -1.5) {};
    \node  (8) at (-2.5, 1.5) {};
    \node  (9) at (-2.5, 1.5) {$C_1$};
    \node  (11) at (-2.5, 0.5) {$C_2$};
    \node  (12) at (-2.5, -0.5) {$C_3$};
    \node  (13) at (-2.5, -1.5) {$C_4$};

    \draw (0.center) to (1.center);
    \draw (3.center) to (2.center);
    \draw (4.center) to (5.center);
    \draw (7.center) to (6.center);
    \end{tikzpicture}
    \end{center}
  \end{minipage}
\begin{minipage}{.45\textwidth}
  \begin{center}
  $\downarrow\text{d.c}$
  \end{center}
\end{minipage}
\begin{minipage}{.45\textwidth}
  \begin{center}
  $\downarrow$
  \end{center}
\end{minipage}
\begin{minipage}{.45\textwidth}
\begin{center}
    \begin{tikzpicture}
        [inner sep=1mm,
      R/.style={circle,draw=black!255,fill=white!20,thick},
      T/.style={circle,draw=red!255,fill=white!255,thick},
      S/.style={rectangle,draw=black!255,fill=white!255,thick}]
     \node[R] (a) [label=below:$-2$]{};
     \node[R] (b) [right=of a, label=below:$-2$] {};
     \node[T] (c) [right=of b, label=below:$-1$] {};

     \node[R] (d) [above right = .5cm and .5cm of c, label = above:$C_4$] {};
     \node[R] (e) [above left = .5cm and .4cm of b, label = above:$C_1$] {};
     \node[R] (f) [above right = .5cm and .4cm of b, label = above:$C_2$] {};
     \node[R] (g) [above = .5cm of c, label = above:$C_3$] {};
  
     \draw[-] (a)--(b);
     \draw[-] (b)--(c);
     \draw[-] (c)--(d);
     \draw[-] (b)--(e);
     \draw[-] (b)--(f);
     \draw[-] (c)--(g);

     \draw[edge,line width=10pt] (b) -- (e);
    \end{tikzpicture}
\end{center}
  \end{minipage}
  \begin{minipage}{.45\textwidth}
  \begin{center}
    \begin{tikzpicture}
    \node  (0) at (-2, 1.5) {};
    \node  (1) at (2, 1.5) {};
    \node  (2) at (-2, 0.5) {};
    \node  (3) at (2, 0.5) {};
    \node  (4) at (-2, -0.5) {};
    \node  (5) at (2, -0.5) {};
    \node  (6) at (-2, -1.5) {};
    \node  (7) at (2, -1.5) {};
    \node  (11) at (-2.5, 0.5) {$C_2$};
    \node  (12) at (-2.5, -0.5) {$C_3$};
    \node  (13) at (-2.5, -1.5) {$C_4$};
    \node  (14) at (-2.5, 1.5) {$C_1$};
    \node  (15) at (-1.5, 2) {};
    \node  (16) at (-1.5, 0) {};
    \node  (17) at (-1, 0) {};
    \node  (18) at (-1, -1) {};
    \node  (19) at (-0.5, -1) {};
    \node  (20) at (-0.5, -2) {};

    \draw (0.center) to (1.center);
    \draw (3.center) to (2.center);
    \draw (4.center) to (5.center);
    \draw (7.center) to (6.center);
    \draw [red] (15.center) to (16.center);
    \draw [red] (17.center) to (18.center);
    \draw [red] (19.center) to (20.center);
    \end{tikzpicture}
    \end{center}
  \end{minipage}
\begin{minipage}{.45\textwidth}
  \begin{center}
  $\downarrow\text{d.c}$
  \end{center}
\end{minipage}
\begin{minipage}{.45\textwidth}
  \begin{center}
  $\downarrow$
  \end{center}
\end{minipage}
\begin{minipage}{.45\textwidth}
    \begin{center}
    \begin{tikzpicture}
        [inner sep=1mm,
      R/.style={circle,draw=black!255,fill=white!20,thick},
      T/.style={circle,draw=red!255,fill=white!255,thick}]
      \node[R] (a) [label=below:$-2$]{};
     \node[T] (b) [right=of a, label=below:$-1$] {};

     \node[R] (d) [below right = .5cm and .5cm of b, label=right:$C_4$] {};
     \node[R] (e) [above left = .5cm and .4cm of b, label=above:$C_1$] {};
     \node[R] (f) [above right = .5cm and .4cm of b, label=above:$C_2$] {};
     \node[R] (g) [right = .5cm of b, label=right:$C_3$] {};

     \draw[-] (a)--(b);
     \draw[-] (b)--(e);
     \draw[-] (b)--(f);
     \draw[-] (b)--(d);
     \draw[-] (b)--(g);

     \draw[edge,line width=10pt] (b) -- (e);
    \end{tikzpicture}
    \end{center}
  \end{minipage}
  \begin{minipage}{.45\textwidth}
  \begin{center}
    \begin{tikzpicture}
    \node  (0) at (-2, 1.5) {};
    \node  (1) at (2, 1.5) {};
    \node  (2) at (-2, 0.5) {};
    \node  (3) at (2, 0.5) {};
    \node  (4) at (-2, -0.5) {};
    \node  (5) at (2, -0.5) {};
    \node  (6) at (-2, -1.5) {};
    \node  (7) at (2, -1.5) {};
    \node  (11) at (-2.5, 0.5) {$C_2$};
    \node  (12) at (-2.5, -0.5) {$C_3$};
    \node  (13) at (-2.5, -1.5) {$C_4$};
    \node  (14) at (-2.5, 1.5) {$C_1$};
    \node  (15) at (-1.5, 2) {};
    \node  (16) at (-1.5, 0) {};
    \node  (17) at (-1, 0) {};
    \node  (18) at (-1, -1) {};
    \node  (19) at (-0.5, -1) {};
    \node  (20) at (-0.5, -2) {};
    \node  (21) at (0, 0) {};
    \node  (22) at (0, -2) {};
    \node  (27) at (0, 2) {};
    \node  (28) at (0, 1) {};

    \draw (0.center) to (1.center);
    \draw (3.center) to (2.center);
    \draw (4.center) to (5.center);
    \draw (7.center) to (6.center);
    \draw [red] (15.center) to (16.center);
    \draw [red] (17.center) to (18.center);
    \draw [red] (19.center) to (20.center);
    \draw [red] (21.center) to (22.center);
    \draw [red] (27.center) to (28.center);
    \end{tikzpicture}
    \end{center}
  \end{minipage}
\begin{minipage}{.45\textwidth}
  \begin{center}
  $\downarrow\text{d.c}$
  \end{center}
\end{minipage}
\begin{minipage}{.45\textwidth}
  \begin{center}
  $\downarrow$
  \end{center}
\end{minipage}
\begin{minipage}{.45\textwidth}
    \begin{center}
    \begin{tikzpicture}
        [inner sep=1mm,
      R/.style={circle,draw=black!255,fill=white!20,thick},
      T/.style={circle,draw=red!255,fill=white!255,thick}]
      \node[T] (a) [label=below:$-1$]{};

     \node[R] (d) [below right = .5cm and .5cm of a, label=right:$C_4$] {};
     \node[R] (e) [above left = .5cm and .4cm of a, label=above:$C_1$] {};
     \node[R] (f) [above right = .5cm and .4cm of a, label=above:$C_2$] {};
     \node[R] (g) [right = .5cm of a, label=right:$C_3$] {};

     \draw[-] (a)--(d);
     \draw[-] (a)--(e);
     \draw[-] (a)--(f);
     \draw[-] (a)--(g);
    \end{tikzpicture}
    \end{center}
  \end{minipage}
  \begin{minipage}{.45\textwidth}
  \begin{center}
    \begin{tikzpicture}
    \node  (0) at (-2, 1.5) {};
    \node  (1) at (2, 1.5) {};
    \node  (2) at (-2, 0.5) {};
    \node  (3) at (2, 0.5) {};
    \node  (4) at (-2, -0.5) {};
    \node  (5) at (2, -0.5) {};
    \node  (6) at (-2, -1.5) {};
    \node  (7) at (2, -1.5) {};
    \node  (11) at (-2.5, 0.5) {$C_2$};
    \node  (12) at (-2.5, -0.5) {$C_3$};
    \node  (13) at (-2.5, -1.5) {$C_4$};
    \node  (14) at (-2.5, 1.5) {$C_1$};
    \node  (15) at (-1.5, 2) {};
    \node  (16) at (-1.5, 0) {};
    \node  (17) at (-1, 0) {};
    \node  (18) at (-1, -1) {};
    \node  (19) at (-0.5, -1) {};
    \node  (20) at (-0.5, -2) {};
    \node  (21) at (0, 0) {};
    \node  (22) at (0, -2) {};
    \node  (23) at (0.5, 1) {};
    \node  (24) at (0.5, -2) {};
    \node  (27) at (0, 2) {};
    \node  (28) at (0, 1) {};

    \draw (0.center) to (1.center);
    \draw (3.center) to (2.center);
    \draw (4.center) to (5.center);
    \draw (7.center) to (6.center);
    \draw [red] (15.center) to (16.center);
    \draw [red] (17.center) to (18.center);
    \draw [red] (19.center) to (20.center);
    \draw [red] (21.center) to (22.center);
    \draw [red] (23.center) to (24.center);
    \draw [red] (27.center) to (28.center);
    \end{tikzpicture}
    \end{center}
  \end{minipage}
\begin{minipage}{.45\textwidth}
  \begin{center}
  $\downarrow$
  \end{center}
\end{minipage}
\begin{minipage}{.45\textwidth}
  \begin{center}
  $\downarrow$
  \end{center}
\end{minipage}
\begin{minipage}{.45\textwidth}
    \;
  \end{minipage}
  \begin{minipage}{.45\textwidth}
  \begin{center}
    \begin{tikzpicture}
    \node  (0) at (-2, 1.5) {};
    \node  (1) at (2, 1.5) {};
    \node  (2) at (-2, 0.5) {};
    \node  (3) at (2, 0.5) {};
    \node  (4) at (-2, -0.5) {};
    \node  (5) at (2, -0.5) {};
    \node  (6) at (-2, -1.5) {};
    \node  (7) at (2, -1.5) {};
    \node  (11) at (-2.5, 0.5) {$C_2$};
    \node  (12) at (-2.5, -0.5) {$C_3$};
    \node  (13) at (-2.5, -1.5) {$C_4$};
    \node  (14) at (-2.5, 1.5) {$C_1$};
    \node  (15) at (-1.5, 2) {};
    \node  (16) at (-1.5, 0) {};
    \node  (17) at (-1, 0) {};
    \node  (18) at (-1, -1) {};
    \node  (19) at (-0.5, -1) {};
    \node  (20) at (-0.5, -2) {};
    \node  (21) at (0, 0) {};
    \node  (22) at (0, -2) {};
    \node  (23) at (0.5, 1) {};
    \node  (24) at (0.5, -2) {};
    \node  (25) at (1, 2) {};
    \node  (26) at (1, -2) {};
    \node  (27) at (0, 2) {};
    \node  (28) at (0, 1) {};

    \draw (0.center) to (1.center);
    \draw (3.center) to (2.center);
    \draw (4.center) to (5.center);
    \draw (7.center) to (6.center);
    \draw [red] (15.center) to (16.center);
    \draw [red] (17.center) to (18.center);
    \draw [red] (19.center) to (20.center);
    \draw [red] (21.center) to (22.center);
    \draw [red] (23.center) to (24.center);
    \draw [red] (25.center) to (26.center);
    \draw [red] (27.center) to (28.center);
    \end{tikzpicture}
    \end{center}
  \end{minipage}
\end{center}
It corresponds to the incidence matrix 
$$\begin{bmatrix}[c|cccccccc]
         C_{1} & 1 & 1 & 1 &   &   &  \\
         C_{2} & 1 & 1 &   & 1 &   &  \\
         C_{3} & 1 &   & 1 & 1 & 1 &  \\
         C_{4} & 1 &   & 1 & 1 &   & 1         
        \end{bmatrix}$$
(3) The $P$-resolution with two Wahl singularities $[2,5]$ and $[4]$
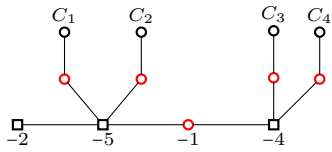
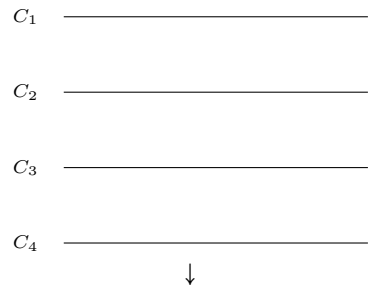
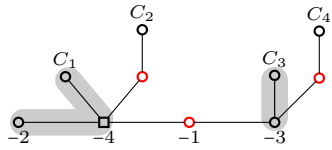
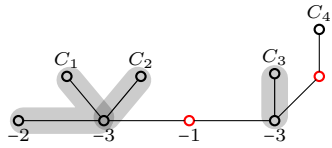
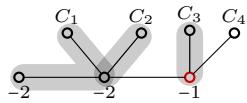
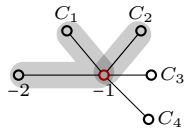
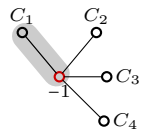
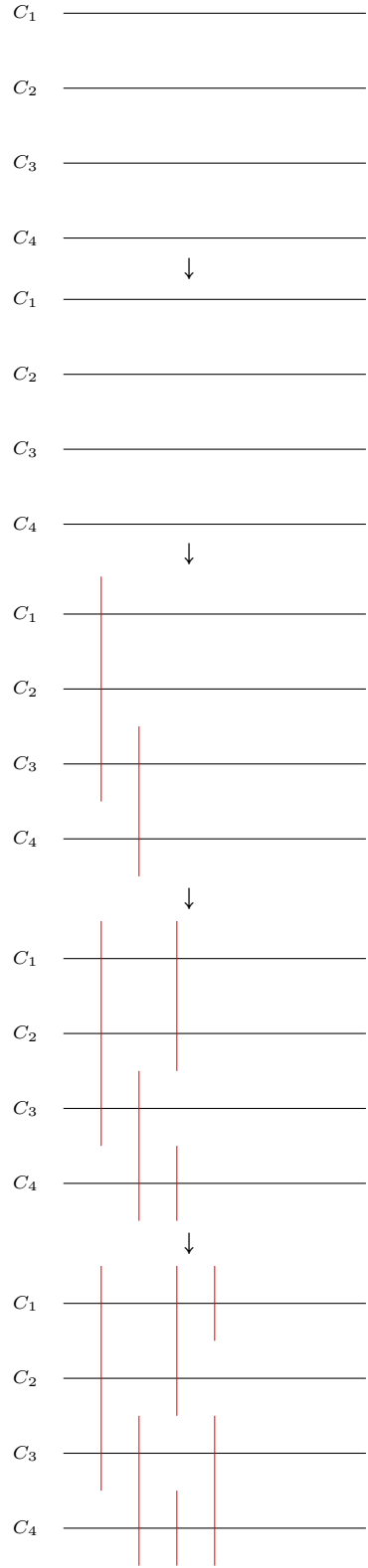
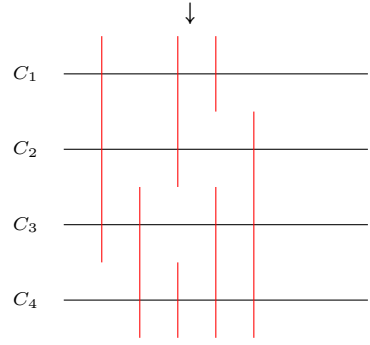
\begin{center}
  \begin{minipage}{.45\textwidth}
  \begin{center}
    \begin{tikzpicture}
        [inner sep=1mm,
      R/.style={circle,draw=black!255,fill=white!20,thick},
      T/.style={circle,draw=red!255,fill=white!255,thick},
      S/.style={rectangle,draw=black!255,fill=white!255,thick}]
     \node[S] (a) [label=below:$-2$]{};
     \node[S] (b) [right=of a, label=below:$-5$] {};
     \node[T] (bc) [right=of b, label=below:$-1$] {};
     \node[S] (c) [right=of bc, label=below:$-4$] {};

     \node[T] (d) [above right = .5cm and .5cm of c] {};
     \node[T] (e) [above left = .5cm and .4cm of b] {};
     \node[T] (f) [above right = .5cm and .4cm of b] {};
      \node[T] (g) [above = .5cm of c] {};

      \node[R] (h) [above = .5cm of d, label = above:$C_4$] {};
     \node[R] (i) [above = .5cm of e, label = above:$C_1$] {};
     \node[R] (j) [above = .5cm of f, label = above:$C_2$] {};
     \node[R] (k) [above = .5cm of g, label = above:$C_3$] {};

     \draw[-] (a)--(b);
     \draw[-] (b)--(bc);
      \draw[-] (bc)--(c);

     \draw[-] (c)--(d);
     \draw[-] (b)--(e);
     \draw[-] (b)--(f);
     \draw[-] (c)--(g);

     \draw[-] (d)--(h);
     \draw[-] (e)--(i);
     \draw[-] (f)--(j);
      \draw[-] (g)--(k);
    \end{tikzpicture}
    \end{center}
  \end{minipage}
  \begin{minipage}{.45\textwidth}
  \begin{center}
    \begin{tikzpicture}
    \node  (0) at (-2, 1.5) {};
    \node  (1) at (2, 1.5) {};
    \node  (2) at (-2, 0.5) {};
    \node  (3) at (2, 0.5) {};
    \node  (4) at (-2, -0.5) {};
    \node  (5) at (2, -0.5) {};
    \node  (6) at (-2, -1.5) {};
    \node  (7) at (2, -1.5) {};
    \node  (8) at (-2.5, 1.5) {};
    \node  (9) at (-2.5, 1.5) {$C_1$};
    \node  (11) at (-2.5, 0.5) {$C_2$};
    \node  (12) at (-2.5, -0.5) {$C_3$};
    \node  (13) at (-2.5, -1.5) {$C_4$};

    \draw (0.center) to (1.center);
    \draw (3.center) to (2.center);
    \draw (4.center) to (5.center);
    \draw (7.center) to (6.center);
    \end{tikzpicture}
    \end{center}
  \end{minipage}
\begin{minipage}{.45\textwidth}
  \begin{center}
  $\downarrow\text{flip}$
  \end{center}
\end{minipage}
\begin{minipage}{.45\textwidth}
  \begin{center}
  $\downarrow$
  \end{center}
\end{minipage}
\begin{minipage}{.45\textwidth}
  \begin{center}
    \begin{tikzpicture}
        [inner sep=1mm,
      R/.style={circle,draw=black!255,fill=white!20,thick},
      T/.style={circle,draw=red!255,fill=white!255,thick},
      S/.style={rectangle,draw=black!255,fill=white!255,thick}]
     \node[R] (a) [label=below:$-2$]{};
     \node[S] (b) [right=of a, label=below:$-4$] {};
     \node[T] (bc) [right=of b, label=below:$-1$] {};
     \node[R] (c) [right=of bc, label=below:$-3$] {};

     \node[T] (d) [above right = .5cm and .5cm of c] {};
     \node[R] (e) [above left = .5cm and .4cm of b, label = above:$C_1$] {};
     \node[T] (f) [above right = .5cm and .4cm of b] {};
      \node[R] (g) [above = .5cm of c, label = above:$C_3$] {};

      \node[R] (h) [above = .5cm of d, label = above:$C_4$] {};

     \node[R] (j) [above = .5cm of f, label = above:$C_2$] {};

     \draw[-] (a)--(b);
     \draw[-] (b)--(bc);
      \draw[-] (bc)--(c);

     \draw[-] (c)--(d);
     \draw[-] (b)--(e);
     \draw[-] (b)--(f);
     \draw[-] (c)--(g);

     \draw[-] (d)--(h);

     \draw[-] (f)--(j);

     \draw[edge,line width=10pt] (e) -- (b) -- (a);
     \draw[edge,line width=10pt] (g) -- (c);
    \end{tikzpicture}
    \end{center}
  \end{minipage}
  \begin{minipage}{.45\textwidth}
  \begin{center}
    \begin{tikzpicture}
    \node  (0) at (-2, 1.5) {};
    \node  (1) at (2, 1.5) {};
    \node  (2) at (-2, 0.5) {};
    \node  (3) at (2, 0.5) {};
    \node  (4) at (-2, -0.5) {};
    \node  (5) at (2, -0.5) {};
    \node  (6) at (-2, -1.5) {};
    \node  (7) at (2, -1.5) {};
    \node  (8) at (-2.5, 1.5) {};
    \node  (9) at (-2.5, 1.5) {$C_1$};
    \node  (11) at (-2.5, 0.5) {$C_2$};
    \node  (12) at (-2.5, -0.5) {$C_3$};
    \node  (13) at (-2.5, -1.5) {$C_4$};

    \draw (0.center) to (1.center);
    \draw (3.center) to (2.center);
    \draw (4.center) to (5.center);
    \draw (7.center) to (6.center);
    \end{tikzpicture}
    \end{center}
  \end{minipage}
\begin{minipage}{.45\textwidth}
  \begin{center}
  $\downarrow\text{flip}$
  \end{center}
\end{minipage}
\begin{minipage}{.45\textwidth}
  \begin{center}
  $\downarrow$
  \end{center}
\end{minipage}
\begin{minipage}{.45\textwidth}
  \begin{center}
    \begin{tikzpicture}
        [inner sep=1mm,
      R/.style={circle,draw=black!255,fill=white!20,thick},
      T/.style={circle,draw=red!255,fill=white!255,thick},
      S/.style={rectangle,draw=black!255,fill=white!255,thick}]
     \node[R] (a) [label=below:$-2$]{};
     \node[R] (b) [right=of a, label=below:$-3$] {};
     \node[T] (bc) [right=of b, label=below:$-1$] {};
     \node[R] (c) [right=of bc, label=below:$-3$] {};

     \node[T] (d) [above right = .5cm and .5cm of c] {};
     \node[R] (e) [above left = .5cm and .4cm of b, label = above:$C_1$] {};
     \node[R] (f) [above right = .5cm and .4cm of b, label = above:$C_2$] {};
     \node[R] (g) [above = .5cm of c, label = above:$C_3$] {};
     \node[R] (h) [above = .5cm of d, label = above:$C_4$] {};

     \draw[-] (a)--(b);
     \draw[-] (b)--(bc);
     \draw[-] (bc)--(c);
     \draw[-] (c)--(d);
     \draw[-] (b)--(e);
     \draw[-] (b)--(f);
     \draw[-] (c)--(g);
     \draw[-] (d)--(h);

     \draw[edge,line width=10pt] (e) -- (b) -- (a);
     \draw[edge,line width=10pt] (g) -- (c);
     \draw[edge,line width=10pt] (b) -- (f);
    \end{tikzpicture}
    \end{center}
  \end{minipage}
  \begin{minipage}{.45\textwidth}
  \begin{center}
    \begin{tikzpicture}
    \node  (0) at (-2, 1.5) {};
    \node  (1) at (2, 1.5) {};
    \node  (2) at (-2, 0.5) {};
    \node  (3) at (2, 0.5) {};
    \node  (4) at (-2, -0.5) {};
    \node  (5) at (2, -0.5) {};
    \node  (6) at (-2, -1.5) {};
    \node  (7) at (2, -1.5) {};
    \node  (8) at (-2.5, 1.5) {};
    \node  (9) at (-2.5, 1.5) {$C_1$};
    \node  (11) at (-2.5, 0.5) {$C_2$};
    \node  (12) at (-2.5, -0.5) {$C_3$};
    \node  (13) at (-2.5, -1.5) {$C_4$};

    \draw (0.center) to (1.center);
    \draw (3.center) to (2.center);
    \draw (4.center) to (5.center);
    \draw (7.center) to (6.center);
    \end{tikzpicture}
    \end{center}
  \end{minipage}
\begin{minipage}{.45\textwidth}
  \begin{center}
  $\downarrow\text{d.c}$
  \end{center}
\end{minipage}
\begin{minipage}{.45\textwidth}
  \begin{center}
  $\downarrow$
  \end{center}
\end{minipage}
\begin{minipage}{.45\textwidth}
\begin{center}
    \begin{tikzpicture}
        [inner sep=1mm,
      R/.style={circle,draw=black!255,fill=white!20,thick},
      T/.style={circle,draw=red!255,fill=white!255,thick},
      S/.style={rectangle,draw=black!255,fill=white!255,thick}]
     \node[R] (a) [label=below:$-2$]{};
     \node[R] (b) [right=of a, label=below:$-2$] {};
     \node[T] (c) [right=of b, label=below:$-1$] {};

     \node[R] (d) [above right = .5cm and .5cm of c, label = above:$C_4$] {};
     \node[R] (e) [above left = .5cm and .4cm of b, label = above:$C_1$] {};
     \node[R] (f) [above right = .5cm and .4cm of b, label = above:$C_2$] {};
     \node[R] (g) [above = .5cm of c, label = above:$C_3$] {};

     \draw[-] (a)--(b);
     \draw[-] (b)--(c);
     \draw[-] (c)--(d);
     \draw[-] (b)--(e);
     \draw[-] (b)--(f);
     \draw[-] (c)--(g);

     \draw[edge,line width=10pt] (e) -- (b) -- (a);
     \draw[edge,line width=10pt] (g) -- (c);
     \draw[edge,line width=10pt] (b) -- (f);
    \end{tikzpicture}
\end{center}
  \end{minipage}
  \begin{minipage}{.45\textwidth}
  \begin{center}
    \begin{tikzpicture}
    \node  (0) at (-2, 0.5) {};
    \node  (1) at (-2, -0.5) {};
    \node  (2) at (2, 0.5) {};
    \node  (3) at (2, -0.5) {};
    \node  (4) at (-2, 1.5) {};
    \node  (5) at (2, 1.5) {};
    \node  (6) at (-2, -1.5) {};
    \node  (7) at (2, -1.5) {};
    \node  (8) at (-2.5, 1.5) {$C_1$};
    \node  (9) at (-2.5, 0.5) {$C_2$};
    \node  (10) at (-2.5, -0.5) {$C_3$};
    \node  (11) at (-2.5, -1.5) {$C_4$};
    \node  (12) at (-1.5, 2) {};
    \node  (13) at (-1.5, -1) {};
    \node  (14) at (-1, 0) {};
    \node  (15) at (-1, -2) {};

    \draw (4.center) to (5.center);
    \draw (0.center) to (2.center);
    \draw (1.center) to (3.center);
    \draw (6.center) to (7.center);
    \draw [red] (12.center) to (13.center);
    \draw [red] (14.center) to (15.center);
    \end{tikzpicture}
    \end{center}
  \end{minipage}
\begin{minipage}{.45\textwidth}
  \begin{center}
  $\downarrow\text{d.c}$
  \end{center}
\end{minipage}
\begin{minipage}{.45\textwidth}
  \begin{center}
  $\downarrow$
  \end{center}
\end{minipage}
\begin{minipage}{.45\textwidth}
    \begin{center}
    \begin{tikzpicture}
        [inner sep=1mm,
      R/.style={circle,draw=black!255,fill=white!20,thick},
      T/.style={circle,draw=red!255,fill=white!255,thick}]
      \node[R] (a) [label=below:$-2$]{};
     \node[T] (b) [right=of a, label=below:$-1$] {};

     \node[R] (d) [below right = .5cm and .5cm of b, label=right:$C_4$] {};
     \node[R] (e) [above left = .5cm and .4cm of b, label=above:$C_1$] {};
     \node[R] (f) [above right = .5cm and .4cm of b, label=above:$C_2$] {};
     \node[R] (g) [right = .5cm of b, label=right:$C_3$] {};

     \draw[-] (a)--(b);
     \draw[-] (b)--(e);
     \draw[-] (b)--(f);
     \draw[-] (b)--(d);
     \draw[-] (b)--(g);

     \draw[edge,line width=10pt] (e) -- (b) -- (a);
     \draw[edge,line width=10pt] (b) -- (f);
    \end{tikzpicture}
    \end{center}
  \end{minipage}
  \begin{minipage}{.45\textwidth}
  \begin{center}
    \begin{tikzpicture}
    \node  (0) at (-2, 0.5) {};
    \node  (1) at (-2, -0.5) {};
    \node  (2) at (2, 0.5) {};
    \node  (3) at (2, -0.5) {};
    \node  (4) at (-2, 1.5) {};
    \node  (5) at (2, 1.5) {};
    \node  (6) at (-2, -1.5) {};
    \node  (7) at (2, -1.5) {};
    \node  (8) at (-2.5, 1.5) {$C_1$};
    \node  (9) at (-2.5, 0.5) {$C_2$};
    \node  (10) at (-2.5, -0.5) {$C_3$};
    \node  (11) at (-2.5, -1.5) {$C_4$};
    \node  (12) at (-1.5, 2) {};
    \node  (13) at (-1.5, -1) {};
    \node  (14) at (-1, 0) {};
    \node  (15) at (-1, -2) {};
    \node  (16) at (-0.5, 2) {};
    \node  (17) at (-0.5, 0) {};
    \node  (18) at (-0.5, -1) {};
    \node  (19) at (-0.5, -2) {};

    \draw (4.center) to (5.center);
    \draw (0.center) to (2.center);
    \draw (1.center) to (3.center);
    \draw (6.center) to (7.center);
    \draw [red] (12.center) to (13.center);
    \draw [red] (14.center) to (15.center);
    \draw [red] (16.center) to (17.center);
    \draw [red] (18.center) to (19.center);
    \end{tikzpicture}
    \end{center}
  \end{minipage}
\begin{minipage}{.45\textwidth}
  \begin{center}
  $\downarrow\text{d.c}$
  \end{center}
\end{minipage}
\begin{minipage}{.45\textwidth}
  \begin{center}
  $\downarrow$
  \end{center}
\end{minipage}
\begin{minipage}{.45\textwidth}
    \begin{center}
    \begin{tikzpicture}
        [inner sep=1mm,
      R/.style={circle,draw=black!255,fill=white!20,thick},
      T/.style={circle,draw=red!255,fill=white!255,thick}]
      \node[T] (a) [label=below:$-1$]{};

     \node[R] (d) [below right = .5cm and .5cm of a, label=right:$C_4$] {};
     \node[R] (e) [above left = .5cm and .4cm of a, label=above:$C_1$] {};
     \node[R] (f) [above right = .5cm and .4cm of a, label=above:$C_2$] {};
     \node[R] (g) [right = .5cm of a, label=right:$C_3$] {};

     \draw[-] (a)--(d);
     \draw[-] (a)--(e);
     \draw[-] (a)--(f);
     \draw[-] (a)--(g);

     \draw[edge,line width=10pt] (a) -- (e);
    \end{tikzpicture}
    \end{center}
  \end{minipage}
  \begin{minipage}{.45\textwidth}
  \begin{center}
    \begin{tikzpicture}
    \node  (0) at (-2, 0.5) {};
    \node  (1) at (-2, -0.5) {};
    \node  (2) at (2, 0.5) {};
    \node  (3) at (2, -0.5) {};
    \node  (4) at (-2, 1.5) {};
    \node  (5) at (2, 1.5) {};
    \node  (6) at (-2, -1.5) {};
    \node  (7) at (2, -1.5) {};
    \node  (8) at (-2.5, 1.5) {$C_1$};
    \node  (9) at (-2.5, 0.5) {$C_2$};
    \node  (10) at (-2.5, -0.5) {$C_3$};
    \node  (11) at (-2.5, -1.5) {$C_4$};
    \node  (12) at (-1.5, 2) {};
    \node  (13) at (-1.5, -1) {};
    \node  (14) at (-1, 0) {};
    \node  (15) at (-1, -2) {};
    \node  (16) at (-0.5, 2) {};
    \node  (17) at (-0.5, 0) {};
    \node  (18) at (-0.5, -1) {};
    \node  (19) at (-0.5, -2) {};
    \node  (20) at (0, 2) {};
    \node  (21) at (0, 1) {};
    \node  (22) at (0, 0) {};
    \node  (23) at (0, -2) {};

    \draw (4.center) to (5.center);
    \draw (0.center) to (2.center);
    \draw (1.center) to (3.center);
    \draw (6.center) to (7.center);
    \draw [red] (12.center) to (13.center);
    \draw [red] (14.center) to (15.center);
    \draw [red] (16.center) to (17.center);
    \draw [red] (18.center) to (19.center);
    \draw [red] (20.center) to (21.center);
    \draw [red] (22.center) to (23.center);
    \end{tikzpicture}
    \end{center}
  \end{minipage}
\begin{minipage}{.45\textwidth}
  \begin{center}
  $\downarrow\text{d.c}$
  \end{center}
\end{minipage}
\begin{minipage}{.45\textwidth}
  \begin{center}
  $\downarrow$
  \end{center}
\end{minipage}
\begin{minipage}{.45\textwidth}
    \;
  \end{minipage}
  \begin{minipage}{.45\textwidth}
  \begin{center}
    \begin{tikzpicture}
    \node  (0) at (-2, 0.5) {};
    \node  (1) at (-2, -0.5) {};
    \node  (2) at (2, 0.5) {};
    \node  (3) at (2, -0.5) {};
    \node  (4) at (-2, 1.5) {};
    \node  (5) at (2, 1.5) {};
    \node  (6) at (-2, -1.5) {};
    \node  (7) at (2, -1.5) {};
    \node  (8) at (-2.5, 1.5) {$C_1$};
    \node  (9) at (-2.5, 0.5) {$C_2$};
    \node  (10) at (-2.5, -0.5) {$C_3$};
    \node  (11) at (-2.5, -1.5) {$C_4$};
    \node  (12) at (-1.5, 2) {};
    \node  (13) at (-1.5, -1) {};
    \node  (14) at (-1, 0) {};
    \node  (15) at (-1, -2) {};
    \node  (16) at (-0.5, 2) {};
    \node  (17) at (-0.5, 0) {};
    \node  (18) at (-0.5, -1) {};
    \node  (19) at (-0.5, -2) {};
    \node  (20) at (0, 2) {};
    \node  (21) at (0, 1) {};
    \node  (22) at (0, 0) {};
    \node  (23) at (0, -2) {};
    \node  (24) at (0.5, 1) {};
    \node  (25) at (0.5, -2) {};

    \draw (4.center) to (5.center);
    \draw (0.center) to (2.center);
    \draw (1.center) to (3.center);
    \draw (6.center) to (7.center);
    \draw [red] (12.center) to (13.center);
    \draw [red] (14.center) to (15.center);
    \draw [red] (16.center) to (17.center);
    \draw [red] (18.center) to (19.center);
    \draw [red] (20.center) to (21.center);
    \draw [red] (22.center) to (23.center);
    \draw [red] (24.center) to (25.center);
    \end{tikzpicture}
    \end{center}
  \end{minipage}
\end{center}
It corresponds to the incidence matrix 
$$\begin{bmatrix}[c|cccccccc]
         C_{1} & 1 & 1 & 1 &   &    \\
         C_{2} & 1 & 1 &   & 1 &    \\
         C_{3} & 1 &   & 1 & 1 & 1  \\
         C_{4} &   & 1 & 1 & 1 & 1         
        \end{bmatrix}$$
\end{example}
Let $[a_1, \cdots, a_r]$ be a cyclic quotient surface singularity with the usual sandwiched structure(cf Figure~\ref{fig:Usual Sandwiched structure of a cyclic quotient surface singularity}). For future reference, we present two lemmas that deal with specific situations.
\begin{lemma}
\label{lem:No free point ensure that the last curve is an exceptional of Wahl singularity}
  Assume that $a_i \geq 3$. Then there exists at least one decorated curve $C_i$ connected to the exceptional curve $A_i$ through a $(-1)$-curve $E$. If $C_i$ has not free points in the picture deformation, then the curve $A_r$ is an exceptional curve of a Wahl singularity in the corresponding $P$-resolution.
\end{lemma}
\begin{proof}
  If the decorated curve $C_i$ does not have a free point, then it means that if a $(-1)$-curve in the process of the flips and divisorial contractions is connected to $C_i$, then it must be connected to other decorated curve. Particularly, the first $(-1)$-curve connects the curve $C_i$ and other curve. In the case, other decorated curve must degenerate to the exceptional curve $A_i$. Therefore $A_i$ is an exceptional curve of a Wahl singularity.
\end{proof}
\begin{lemma}
\label{lem:The 1-column ensures that the last curve is not an exceptional of Wahl singularity}
  If there exists a column that all entries are $1$ in an incidence matrix, then the curve $A_r$ is not an exceptional curve of a Wahl singularity in the corresponding $P$-resolution.
\end{lemma}
\begin{proof}
  Suppose that the exceptional curve $A_r$ is a exceptional curve of a Wahl singularity. After the flips until the Wahl singularity disappears, we arrive at the following step.
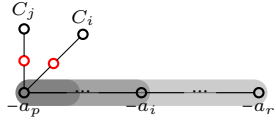
\begin{center}
\begin{tikzpicture}
  [inner sep=1mm,
R/.style={circle,draw=black!255,fill=white!20,thick},
T/.style={circle,draw=red!255,fill=white!255,thick}]
    \node[R] (a) [label=below:$-a_p$]{};
    \node    (b) [right=.5cm of a] {$\cdots$};
    \node[R] (c) [right=.5cm of b, label=below:$-a_{i}$]{};
    \node    (d) [right=.5cm of c] {$\cdots$};
    \node[R] (e) [right=.5cm of d, label=below:$-a_r$]{};
    \node[T] (f) [above right=.3cm and .3cm of a]{};
    \node[R] (g) [above right=.3cm and .3cm of f, label=above:$C_{i}$]{};
    \node[T] (h) [above=.3cm of a]{};
    \node[R] (i) [above=.3cm of h, label=above:$C_{j}$]{};

    \draw[-] (a)--(b);
    \draw[-] (b)--(c);
    \draw[-] (c)--(d);
    \draw[-] (d)--(e);
    \draw[-] (a)--(f);
    \draw[-] (f)--(g);
    \draw[-] (a)--(h);
    \draw[-] (h)--(i);

    \draw[edge,line width=10pt] (a) -- (e);
    \draw[edge,line width=10pt] (a) -- (c);
    \draw[edge,line width=10pt] (a) -- (b);
\end{tikzpicture}
\end{center}
We assume that $A_p$ is the initial curve of the Wahl singularity and that the decorated curves $C_{i}$ and $C_j$ are connected to $A_p$ through a $(-1)$-curve respectively. Since $A_p$ is the initial curve, $a_p$ is greater than or equal to $4$. Therefore, we can assume the two decorated curves. We follow the decorated curves $C_i$ and $C_j$. In this step, $C_i$ and $C_j$ are not connected. After divisorial contractions, we obtain:
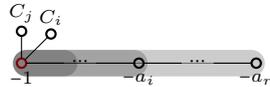
\begin{center}
\begin{tikzpicture}
  [inner sep=1mm,
R/.style={circle,draw=black!255,fill=white!20,thick},
T/.style={circle,draw=red!255,fill=white!255,thick}]
    \node[T] (a) [label=below:$-1$]{};
    \node    (b) [right=.5cm of a] {$\cdots$};
    \node[R] (c) [right=.5cm of b, label=below:$-a_{i}$]{};
    \node    (d) [right=.5cm of c] {$\cdots$};
    \node[R] (e) [right=.5cm of d, label=below:$-a_r$]{};
    \node[R] (f) [above right=.3cm and .3cm of a, label=above:$C_{i}$]{};

    \node[R] (h) [above=.3cm of a, label=above:$C_{j}$]{};

    \draw[-] (a)--(b);
    \draw[-] (b)--(c);
    \draw[-] (c)--(d);
    \draw[-] (d)--(e);
    \draw[-] (a)--(f);

    \draw[-] (a)--(h);

    \draw[edge,line width=10pt] (a) -- (e);
    \draw[edge,line width=10pt] (a) -- (c);
    \draw[edge,line width=10pt] (a) -- (b);
\end{tikzpicture}
\end{center}
The $(-1)$-curve does not connect $C_i$ and the decorated curve that degenerates to the curve $A_p$. Note that during divisorial contractions, every $(-1)$-curve does not connect $C_i$ and some decorated curves that degenerate. Therefore, there is no $(-1)$-curve that connects all decorated curves. In the aspect of an incidence matrix, this means that there dose not exist a column whose entries are all $1$.
\end{proof}

\section{Incidence matrices under the different sandwiched structure}
\label{subsection:Incidence matrices under different sandwiched structure}

In this section, we figure out incidence matrices of cyclic quotient surface singularities under a different sandwiched structure.

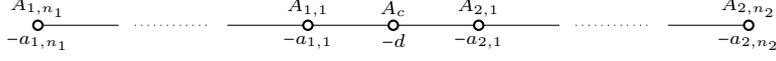
\begin{figure} 
\centering
\begin{tikzpicture}
  [inner sep=1mm,
R/.style={circle,draw=black!255,fill=white!20,thick},
T/.style={circle,draw=black!255,fill=black!255,thick}]
    \node[R] (a) [label=below:$-a_{1,n_1}$, , label=above:$A_{1,n_1}$]{};
    \node    (b) [right=of a] {};

    \node    (h) [right=of b] {};

    \node[R] (c) [right=of h, label=below:$-a_{1,1}$, label=above:$A_{1,1}$] {};

    \node[R] (d) [right=of c, label=below:$-d$, label=above:$A_c$] {};

    \node[R] (e) [right=of d, label=below:$-a_{2,1}$, label=above:$A_{2,1}$] {};
    \node    (f) [right=of e] {};

    \node    (i) [right=of f] {};

    \node[R] (g) [right=of i, label=below:$-a_{2,n_2}$, , label=above:$A_{2,n_2}$] {};

    \draw[-] (a)--(b);
    \draw[dotted] (b)--(h);
    \draw[-] (h)--(c);
    \draw[-] (c)--(d);
    \draw[-] (d)--(e);
    \draw[-] (e)--(f);
    \draw[dotted] (f)--(i);
    \draw[-] (i)--(g);
\end{tikzpicture}
\caption{dual resolution graph of $X$}
\label{fig:resolution}
\end{figure}

Let $(X, 0)$ be a cyclic quotient surface singularity $\frac{1}{n}(1,q)$ where $n/q = [a_{1,n_1}, \dots, a_{1,1}, d, a_{2,1}, \dots, a_{2,n_2}]$ and assume that $d \geq 4$. Then the dual resolution graph is as shown in Figure~\ref{fig:resolution}.
We denote exceptional $(-a_{i,j})$-curves as capital letters $A_{i,j}$.  We say that the curves $A_{1,j}$ are in the first branch and $A_{2,j}$ are in the second branch. We call the curve of degree $-d$ as the central curve $A_c$. We use these notations to broaden our discussion to weighted homogeneous surface singularities.

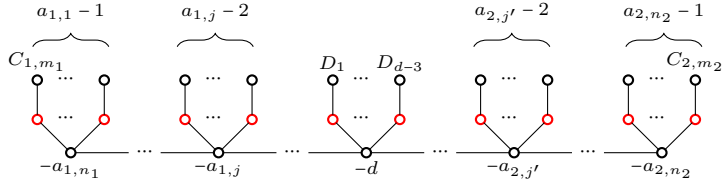
\begin{figure} 
\centering
\begin{tikzpicture}
  [inner sep=1mm,
R/.style={circle,draw=black!255,fill=white!20,thick},
T/.style={circle,draw=red!255,fill=white!255,thick}]
    \node[R] (a) [label=below:$-a_{1,n_1}$]{};
    \node    (b) [right=.7cm of a]{$\cdots$};
    \node[R] (c) [right=.7cm of b, label=below:$-a_{1,j}$]{};
    \node    (d) [right=.7cm of c]{$\cdots$};
    \node[R] (e) [right=.7cm of d, label=below:$-d$]{};
    \node    (f) [right=.7cm of e]{$\cdots$};
    \node[R] (g) [right=.7cm of f, label=below:$-a_{2,j'}$]{};
    \node    (h) [right=.7cm of g]{$\cdots$};
    \node[R] (i) [right=.7cm of h, label=below:$-a_{2,n_2}$]{};

    \node[T] (j) [above left=.5cm of a]{};
    \node    (k) [right=.1cm of j]{$\cdots$};
    \node[T] (l) [above right=.5cm of a]{};
    \node[R] (m) [above=.4cm of j, label=above:$C_{1,m_1}$]{};
    \node    (k) [right=.1cm of m]{$\cdots$};
    \node[R] (n) [above=.4cm of l]{};

    \node[T] (o) [above left=.5cm of c]{};
    \node    (p) [right=.1cm of o]{$\cdots$};
    \node[T] (q) [above right=.5cm of c]{};
    \node[R] (r) [above=.4cm of o]{};
    \node    (k) [right=.1cm of r]{$\cdots$};
    \node[R] (s) [above=.4cm of q]{};

    \node[T] (t) [above left=.5cm of g]{};
    \node    (u) [right=.1cm of t]{$\cdots$};
    \node[T] (v) [above right=.5cm of g]{};
    \node[R] (w) [above=.4cm of t]{};
    \node    (k) [right=.1cm of w]{$\cdots$};
    \node[R] (x) [above=.4cm of v]{};

    \node[T] (y) [above left=.5cm of i]{};
    \node    (z) [right=.1cm of y]{$\cdots$};
    \node[T] (aa) [above right=.5cm of i]{};
    \node[R] (ab) [above=.4cm of y]{};
    \node    (k) [right=.1cm of ab]{$\cdots$};
    \node[R] (ac) [above=.4cm of aa, label=above:$C_{2,m_2}$]{};

    \node[T] (ad) [above left=.5cm of e]{};
    \node    (ae) [right=.1cm of ad]{$\cdots$};
    \node[T] (af) [above right=.5cm of e]{};
    \node[R] (ag) [above=.4cm of ad, label=above:$D_{1}$]{};
    \node    (k) [right=.1cm of ag]{$\cdots$};
    \node[R] (ah) [above=.4cm of af, label=above:$D_{d-3}$]{};

    \draw[-] (a)--(b);
    \draw[-] (b)--(c);
    \draw[-] (c)--(d);
    \draw[-] (d)--(e);
    \draw[-] (e)--(f);
    \draw[-] (f)--(g);
    \draw[-] (g)--(h);
    \draw[-] (h)--(i);
    
    \draw[-] (a)--(j);
    \draw[-] (a)--(l);
    \draw[-] (j)--(m);
    \draw[-] (l)--(n);

    \draw[-] (c)--(o);
    \draw[-] (c)--(q);
    \draw[-] (r)--(o);
    \draw[-] (s)--(q);

    \draw[-] (t)--(g);
    \draw[-] (v)--(g);
    \draw[-] (w)--(t);
    \draw[-] (x)--(v);

    \draw[-] (y)--(i);
    \draw[-] (aa)--(i);
    \draw[-] (ab)--(y);
    \draw[-] (ac)--(aa);

    \draw[-] (ad)--(e);
    \draw[-] (af)--(e);
    \draw[-] (ag)--(ad);
    \draw[-] (ah)--(af);

\draw [decorate,decoration={brace,amplitude=5pt,raise=4ex}]
  (-0.5,.8) -- (0.5,.8) node[midway,yshift=3em]{$a_{1,1}-1$};

\draw [decorate,decoration={brace,amplitude=5pt,raise=4ex}]
  (1.4,.8) -- (2.4,.8) node[midway,yshift=3em]{$a_{1,j}-2$};

\draw [decorate,decoration={brace,amplitude=5pt,raise=4ex}]
  (5.3,.8) -- (6.3,.8) node[midway,yshift=3em]{$a_{2,j'}-2$};

\draw [decorate,decoration={brace,amplitude=5pt,raise=4ex}]
  (7.3,.8) -- (8.3,.8) node[midway,yshift=3em]{$a_{2,n_2}-1$};
\end{tikzpicture}

\caption{Sandwiched structure of $X$}
\label{fig:sandwichedsatructureofcqss}
\end{figure}

If we attach $(a_{i,n_i}-1)$ $(-1)$-curves on $A_{i,n_i}$ for $i = 1, 2$, $(a_{1,j}-2)$ $(-1)$-curves on $A_{1,j}$ for $j < n_1$, $(a_{2,j}-2)$ $(-1)$-curves on $A_{2,j}$ for $j < n_2$ and $(d-3)$ $(-1)$-curves on the central curve, the graph(Figure~\ref{fig:resolution}) is contracted to the central curve and finally a smooth point. The graph is therefore sandwiched. By attaching a decorated curve on each $(-1)$-curve, we obtain a sandwiched structure of $X$. We denote the decorated curves on the first branch as $C_{1,j}$, second branch as $C_{2,j}$ and the central curve as $D_k$. The second subscript is ordered inside out. The number of decorated curves on the first branch is the length of the dual Hirzebruch-Jung continued fraction of $[a_{1,1}, \dots, a_{1,n_1}]$, denoted as $m_1$. Similarly the number of decorated curves on the second branch is the length of the dual Hirzebruch-Jung continued fraction of $[a_{2,1}, \dots, a_{2,n_2}]$, denoted as $m_2$. The number of decorated curves on the central curve is $d-3$. See Figure~\ref{fig:sandwichedsatructureofcqss}. 

We want to classify the configurations of the incidence matrices of cyclic quotient surface singularities by analyzing the MMP-algorithm. We start from two types of P-resolutions(more precisely, M-resolutions). The first type is that the central curve is not an exceptional curve of a Wahl singularity. The second type is that the central curve is an exceptional curve of a Wahl singularity. Let 
$$\begin{tikzpicture}
  [inner sep=1mm,
R/.style={circle,draw=black!255,fill=white!20,thick},
T/.style={rectangle,draw=black!255,fill=white!255,thick}]

\node[T] (a) [label=below:$-a_{1,n_1'}$] {};
\node[]          (b) [right=of a] {$\cdots$};
\node[T] (c) [right=of b, label=below:$-a_{1,1}$] {};
\node[T] (d) [right=of c, label=below:$-d$] {};
\node[T] (e) [right=of d, label=below:$-a_{2,1}$] {};
\node[]          (f) [right=of e] {$\cdots$};
\node[T] (g) [right=of f, label=below:$-a_{2,n_2'}$] {};

\draw[-] (a) -- (b);
\draw[-] (b) -- (c);
\draw[-] (c) -- (d);
\draw[-] (d) -- (e);
\draw[-] (e) -- (f);
\draw[-] (f) -- (g);
\end{tikzpicture}$$
be the Wahl singularity with $n_1' < n_1$ and $n_2' < n_2$. If
$$\begin{tikzpicture}
  [inner sep=1mm,
R/.style={circle,draw=black!255,fill=white!20,thick},
T/.style={rectangle,draw=black!255,fill=white!255,thick}]

\node[T] (a) [label=below:$-a_{1,n_1'}$] {};
\node[]          (b) [right=of a] {$\cdots$};
\node[T] (c) [right=of b, label=below:$-a_{1,1}$] {};
\node[T] (d) [right=of c, label=below:$-d$] {};

\draw[-] (a) -- (b);
\draw[-] (b) -- (c);
\draw[-] (c) -- (d);
\end{tikzpicture}$$
 $[a_{1,n_1'}, \dots, a_{1,1}, d, a_{2,1}, \cdots, a_c]$ is a Wahl singularity, then we say $n_2' = 0$. We may assume that the initial curve is in the first branch or the central curve.

\begin{definition}[Incidence matrix of type 1 and 2]
Let $X$ be a cyclic quotient surface singularity with the sandwiched structure as in Figure~\ref{fig:sandwichedsatructureofcqss}. We call an incidence matrix of $X$ is \emph{type 1} if it is induced from a P-resolution that the central curve of the minimal resolution is not an exceptional curve of a Wahl singularity. Otherwise, we call it \emph{type 2-1} if $n_2' > 0$ and \emph{type 2-2} if $n_2' = 0$.
\label{def:Incidence matrix of type 1 and 2}
\end{definition}

\begin{example}[Continued from \ref{eg:19/11}]
  Consider a sandwiched structure on the CQSS $\frac{1}{19}(1,11)$ as Figure~\ref{fig:sandwiched structure of 19/11}. Incidence matrices under the sandwiched structure are : 
  \begin{center}
        $\begin{bmatrix}[c|cccccccc]
         C_{1} & 1 &   &   &   & 1 &   & 1   \\
         C_{2} &   & 1 &   &   &   &   & 1   \\
         C_{3} &   &   & 1 &   &   & 1 & 1   \\
         C_{4} &   &   &   & 1 &   & 1 & 1         
        \end{bmatrix}$
        $\begin{bmatrix}[c|cccccccc]
         C_{1} & 1 &   &   & 1 &   & 1 \\
         C_{2} &   &   &   & 1 & 1 &   \\
         C_{3} &   & 1 &   &   & 1 & 1 \\
         C_{4} &   &   & 1 &   & 1 & 1         
        \end{bmatrix}$
        $\begin{bmatrix}[c|cccccccc]
         C_{1} & 1 & 1 &   & 1 &    \\
         C_{2} & 1 &   &   &   & 1  \\
         C_{3} &   & 1 & 1 &   & 1  \\
         C_{4} &   &   & 1 & 1 & 1         
        \end{bmatrix}$
    \end{center}
  The first one is of type 1. The second and third one are of type 2-2.

\end{example}

\begin{figure} 
\centering
\begin{tikzpicture}
  [inner sep=1mm,
R/.style={circle,draw=black!255,fill=white!20,thick},
T/.style={circle,draw=red!255,fill=white!255,thick}]
    \node[R] (a) [label=below:$-2$, label=left:$A_{1,1}$]{};

    \node[R] (b) [right=of a, label=below:$-4$] {};

    \node[R] (c) [right=of b, label=below:$-3$, label=right:$A_{2,1}$]{};

    \node[R] (d) [above left=of a, label=above:$C_{1,1}$]{};
    \node[T] (d1) [above left= .5cm and .5cm of a]{};

    \node[R] (e) [above=of b, label=above:$D_1$]{};
    \node[T] (e1) [above=.5cm and .5cm of b]{};

    \node[R] (f) [above=of c, label=above:$C_{2,1}$]{};
    \node[R] (g) [above right=of c, label=above:$C_{2,2}$]{};    
    \node[T] (f1) [above=.5cm and .5cm of c]{};
    \node[T] (g1) [above right=.5cm and .5cm of c]{};

    \draw[-] (a)--(b);
    \draw[-] (b)--(c);
    \draw[-] (a)--(d1);
    \draw[-] (d1)--(d);
    \draw[-] (b)--(e1);
    \draw[-] (e1)--(e);
    \draw[-] (c)--(f1);
    \draw[-] (f1)--(f);
    \draw[-] (c)--(g1);
    \draw[-] (g1)--(g);
\end{tikzpicture}
\caption{sandwiched structure of $\frac{1}{19}(1,11)$}
\label{fig:sandwiched structure of 19/11}
\end{figure}

We investigate incidence matrices of each type. First, we assume that the central curve is not an exceptional curve of a Wahl singularity. Then in the procedure of the MMP algorithm, any decorated curve does not degenerate to the central curve. Then the $(-1)$-curve that attached on the decorated curve $D_k$ can not be connected to other decorated curves. Therefore the $(-1)$-curve corresponds to a free point $p_k$ on $D_k$. After flips and divisorial contractions, the central curve becomes a $(-1)$-curve and all decorated curves are attached on the $(-1)$-curve. This $(-1)$-curve corresponds to a point $p_0$ that all decorated curves pass through. See the Figure~\ref{fig:final divisorial contraction}.

\begin{figure} 
\centering
\begin{minipage}{.3\textwidth}
\begin{tikzpicture}
[inner sep=1mm,
R/.style={circle,draw=black!255,fill=white!20,thick},
T/.style={circle,draw=red!255,fill=white!255,thick}]
    \node[R] (a) [label={[yshift=.1cm]:$-d+2$}]{};
    \node[R] (b) [above left=.1cm and .7cm of a, label=left:$C_{1,1}$]{};
    \node    (c) [left=.6cm of a]{$\vdots$};
    \node[R] (d) [below left=.5cm and .7cm of a, label=left:$C_{1,m_1}$]{};
    \node[T] (e) [below left=.5cm and .3cm of a]{};
    \node    (f) [right=.0000000005cm of e]{$\cdots$};
    \node[T] (g) [below right=.5cm and .3cm of a]{};
    \node[R] (h) [below =.5cm of e, label=below:$D_1$]{};
    \node[R] (i) [below =.5cm of g, label=below:$D_{d-3}$]{};
    \node[R] (j) [above right=.1cm and .7cm of a, label=right:$C_{2,1}$]{};
    \node    (k) [right=.6cm of a]{$\vdots$};
    \node[R] (l) [below right=.5cm and .7cm of a, label=right:$C_{2,m_2}$]{};

    \draw[-] (a) -- (b);
    \draw[-] (a) -- (d);
    \draw[-] (a) -- (e);
    \draw[-] (a) -- (g);
    \draw[-] (e) -- (h);
    \draw[-] (i) -- (g);
    \draw[-] (a) -- (j);
    \draw[-] (a) -- (l);
\end{tikzpicture}
\end{minipage}
$\rightarrow$
\begin{minipage}{.3\textwidth}
\begin{tikzpicture}
[inner sep=1mm,
R/.style={circle,draw=black!255,fill=white!20,thick},
T/.style={circle,draw=red!255,fill=white!255,thick}]
    \node[R] (a) [label={[yshift=.1cm]:$-1$}]{};
    \node[R] (b) [above left=.1cm and .7cm of a, label=left:$C_{1,1}$]{};
    \node    (c) [left=.6cm of a]{$\vdots$};
    \node[R] (d) [below left=.5cm and .7cm of a, label=left:$C_{1,m_1}$]{};
    \node[R] (e) [below left=.5cm and .3cm of a, label=below:$D_1$]{};
    \node    (f) [right=.0000000005cm of e]{$\cdots$};
    \node[R] (g) [below right=.5cm and .3cm of a, label=below:$D_{d-3}$]{};
    \node[R] (j) [above right=.1cm and .7cm of a, label=right:$C_{2,1}$]{};
    \node    (k) [right=.6cm of a]{$\vdots$};
    \node[R] (l) [below right=.5cm and .7cm of a, label=right:$C_{2,m_2}$]{};

    \draw[-] (a) -- (b);
    \draw[-] (a) -- (d);
    \draw[-] (a) -- (e);
    \draw[-] (a) -- (g);
    \draw[-] (a) -- (j);
    \draw[-] (a) -- (l);
\end{tikzpicture}
\end{minipage}
\caption{Divisorial contractions in the type 1}
\label{fig:final divisorial contraction}
\end{figure}
Moreover, the $(-1)$-curves that appear in each branch do not connect two decorated curves in the different branches. This means that any two decorated curves in the different branches do not intersect except $p_0$. From the discussions so far, we obtain an incidence matrix of Figure~\ref{fig:Type 1}.
\begin{figure} 
\centering
$\begin{bmatrix}[c|ccccccccccccccccccc]
               \:      & p_0    & p_1    & \cdots & p_{d-3} \\
               \hline
               C_{1,1}     & 1      & \;     & \;     & \; & *      & \cdots & *      & \;      & \;      & \;   \\
               \vdots  & \vdots & \;     & \;     & \; & \vdots & \ddots & \vdots & \;      & \;      & \;       \\
               C_{1,m_1} & 1      & \;     & \;     & \; & *      & \cdots & *      & \;      & \;      & \;     \\
                \hline
               C_{2,1}     & 1      & \;     & \;     & \; & \;     & \;     & \;     & *       & \cdots  & *    \\
               \vdots  & \vdots & \;     & \;     & \; & \;     & \;     & \;     & \vdots  & \ddots  & \vdots   \\
               C_{2,m_2} & 1      & \;     & \;     & \; & \;     & \;     & \;     & *       & \cdots  & *      \\  
                \hline
               D_1     & 1      & 1      & \;     & \;  & \;     & \;     & \;     & \;      & \;      & \;       \\
               \vdots  & \vdots & \;     & \ddots & \;  & \;     & \;     & \;     & \;      & \;      & \;       \\
               D_{d-3} & 1      & \;     & \;     & 1  & \;     & \;     & \;     & \;      & \;      & \;       \\
\end{bmatrix}$
\caption{Type 1}
\label{fig:Type 1}
\end{figure}
In Figure~\ref{fig:Type 1}, blank entries mean $0$ entries. And $*$ entries mean that the entries are $0$ or $1$ but for each column consisting of  $*$ entries, at least one of the $*$ entries in the column is $1$.

\begin{lemma}[Incidence matrix of type 1]
An incidence matrix of type 1 is of the form as shown in Figure~\ref{fig:Type 1}.
\end{lemma}

Second, we assume that the central curve is an exceptional curve of a Wahl singularity. Then some decorated curves degenerate to the central curve. By the following lemma, we can assume that only decorated curves in the first branch degenerate to the central curve.

\begin{lemma}[Park-Shin \cite{park2022deformations}*{Lemma 5.16}]
Let 
\begin{tikzpicture}[scale=0.5]
 \node[rectangle] (10) at (1,0) [label=above:{$A_{1,n_1'}$}] {};

\node[empty] (20) at (2,0) [] {};
\node[empty] (250) at (2.5,0) [] {};

\node[rectangle] (350) at (3.5,0) [label=above:{$A_{1,p}$}] {};

\node[empty] (450) at (4.5,0) [] {};
\node[empty] (50) at (5,0) [] {};

\node[rectangle] (60) at (6,0) [label=above:{$A_{c}$}] {};

\node[empty] (70) at (7,0) [] {};
\node[empty] (750) at (7.5,0) [] {};

\node[rectangle] (850) at (8.5,0) [label=above:{$A_{2,n_2'}$}] {};

\draw [-] (10)--(20);
\draw [dotted] (20)--(250);
\draw [-] (250)--(350);
\draw [-] (350)--(450);
\draw [dotted] (450)--(50);
\draw [-] (50)--(60);
\draw [-] (60)--(70);
\draw [dotted] (70)--(750);
\draw [-] (750)--(850);
\end{tikzpicture}
be the dual resolution graph of a Wahl singularity with $A_{1, p}$ be its initial curve and $A_c$ be the central curve. Let $[a_{1,n_1'}, \dots, a_{2,n_2'}]$ be its Hirzebruch-Jung continued fraction. We consider the sandwiched structure that $(a_{i,n_i'} - 1)$ $(-1)$-curves attached to $A_{i,n_i'}$ for $i = 1, 2$ ; $(a_{i,j} - 2)$ $(-1)$-curves attached to $A_{i,j}$ for $i = 1, 2$ and $1 \leq j \leq n_i'$ ; $(a_{1,p} - 3)$ $(-1)$-curves attached to $A_{1,p}$. Let $\mathfrak{L} = [a_{1,n_1'}, \dots, a_{2,n_2'}]$ be the extremal neighborhood with the $(-1)$-curves. Then we can apply the usual flips to $\mathfrak{L}$ successively starting from the $(-1)$-curves intersecting $A_{1,n_1'}$ to $(-1)$-curves intersecting $A_{1,p}$ until we obtain
\begin{tikzpicture}[scale=0.5]
 \node[circle] (10) at (1,0) [label=above:{$A_{1,p}$}, label=below:{$-a_{1,p}-1$}] {};

\node[empty] (20) at (2,0) [] {};
\node[empty] (250) at (2.5,0) [] {};

\node[circle] (350) at (3.5,0) [label=above:{$A_{c}$}, label=below:{$-a_{c}$}] {};

\node[empty] (450) at (4.5,0) [] {};
\node[empty] (50) at (5,0) [] {};

\node[circle] (60) at (6,0) [label=above:{$A_{2,n_2'}$}, label=below:{$-a_{2,n_2'}$}] {};


\draw [-] (10)--(20);
\draw [dotted] (20)--(250);
\draw [-] (250)--(350);
\draw [-] (350)--(450);
\draw [dotted] (450)--(50);
\draw [-] (50)--(60);
\end{tikzpicture}
without no singularity.
\label{flips on Tsing contaiing the central}
\end{lemma}
\begin{proof}
The proof is similar to Lemma 5.16 of Park-Shin(\cite{park2022deformations}).
\end{proof}

We follow the MMP algorithm on 
$$\begin{tikzpicture}
  [inner sep=1mm,
R/.style={circle,draw=black!255,fill=white!20,thick},
T/.style={circle,draw=red!255,fill=white!255,thick}]

\node[rectangle] (a) [label=below:$-a_{1,n_1'}$, label=above:$A_{1,n_1'}$] {};
\node[]          (b) [right=of a] {$\cdots$};
\node[rectangle] (c) [right=of b, label=below:$-a_{1,p}$, label=above:$A_{1,p}$] {};
\node[]          (d) [right=of c] {$\cdots$};
\node[rectangle] (e) [right=of d, label=below:$-d$, , label=above:$A_c$] {};
\node[rectangle] (f) [right=of e, label=below:$-a_{2,1}$, label=above:$A_{2,1}$] {};
\node[]          (g) [right=of f] {$\cdots$};
\node[rectangle] (h) [right=of g, label=below:$-a_{2,n_2'}$, label=above:$A_{2,n_2'}$] {};

\draw[-] (a) -- (b);
\draw[-] (b) -- (c);
\draw[-] (c) -- (d);
\draw[-] (d) -- (e);
\draw[-] (e) -- (f);
\draw[-] (f) -- (g);
\draw[-] (g) -- (h);
\end{tikzpicture}$$
precisely.

A $(-1)$-curve passes through the singularity and a decorated curve $C_{1,m_1}$ intersects the $(-1)$-curve. We apply the flip to the $(-1)$-curve. If $a_{1,n_1'} > 2$, then we obtain
$$\begin{tikzpicture}
  [inner sep=1mm,
R/.style={circle,draw=black!255,fill=white!20,thick},
T/.style={circle,draw=red!255,fill=white!255,thick}]

\node[rectangle] (a) [label=below:$-a_{1,n_1'}-1$, label=above:$A_{1,n_1'}$] {};
\node[]          (b) [right=of a] {$\cdots$};
\node[rectangle] (c) [right=of b, label=below:$-a_{1,p}$, label=above:$A_{1,p}$] {};
\node[]          (d) [right=of c] {$\cdots$};
\node[rectangle] (e) [right=of d, label=below:$-d$, , label=above:$A_c$] {};
\node[rectangle] (f) [right=of e, label=below:$-a_{2,1}$, label=above:$A_{2,1}$] {};
\node[]          (g) [right=of f] {$\cdots$};
\node[rectangle] (h) [right=of g, label=below:$-a_{2,n_2'-1}$, label=above:$A_{2,n_2'-1}$] {};
\node[circle] (i) [right=of h, label=below:$-a_{2,n_2'}$, label=above:$A_{2,n_2'}$] {};

\draw[-] (a) -- (b);
\draw[-] (b) -- (c);
\draw[-] (c) -- (d);
\draw[-] (d) -- (e);
\draw[-] (e) -- (f);
\draw[-] (f) -- (g);
\draw[-] (g) -- (h);
\draw[-] (h) -- (i);
\end{tikzpicture}$$
with a degeneration $C_{1,m_1}^+ = C_{1,m_1} + A_{2,n_2'}$. If $a_{1,n_1} = 2$, then we obtain
$$\begin{tikzpicture}
  [inner sep=1mm,
R/.style={circle,draw=black!255,fill=white!20,thick},
T/.style={circle,draw=red!255,fill=white!255,thick}]

\node[rectangle] (a) [label=below:$-a_{1,n_1''}-1$, label=above:$A_{1,n_1''}$] {};
\node[]          (b) [right=of a] {$\cdots$};
\node[rectangle] (c) [right=of b, label=below:$-a_{1,p}$, label=above:$A_{1,p}$] {};
\node[]          (d) [right=of c] {$\cdots$};
\node[rectangle] (e) [right=of d, label=below:$-d$, , label=above:$A_c$] {};
\node[rectangle] (f) [right=of e, label=below:$-a_{2,1}$, label=above:$A_{2,1}$] {};
\node[]          (g) [right=of f] {$\cdots$};
\node[rectangle] (h) [right=of g, label=below:$-a_{2,n_2'-1}$, label=above:$A_{2,n_2'-1}$] {};
\node[circle] (i) [right=of h, label=below:$-a_{2,n_2'}$, label=above:$A_{2,n_2'}$] {};

\draw[-] (a) -- (b);
\draw[-] (b) -- (c);
\draw[-] (c) -- (d);
\draw[-] (d) -- (e);
\draw[-] (e) -- (f);
\draw[-] (f) -- (g);
\draw[-] (g) -- (h);
\draw[-] (h) -- (i);
\end{tikzpicture}$$
with the same degeneration. In any cases, a $(-1)$-curve passes through the new singularity and a decorated curve $C_{1,m_1-1}$ intersects the new $(-1)$-curve. We apply the flip to the $(-1)$-curve. Similarly, we obtain
$$\begin{tikzpicture}
  [inner sep=1mm,
R/.style={circle,draw=black!255,fill=white!20,thick},
T/.style={circle,draw=red!255,fill=white!255,thick}]

\node[rectangle] (a) [label=below:$-a_{1,n_1'''}-1$, label=above:$A_{1,n_1'''}$] {};
\node[]          (b) [right=of a] {$\cdots$};
\node[rectangle] (c) [right=of b, label=below:$-a_{1,p}$, label=above:$A_{1,p}$] {};
\node[]          (d) [right=of c] {$\cdots$};
\node[rectangle] (e) [right=of d, label=below:$-d$, , label=above:$A_c$] {};
\node[rectangle] (f) [right=of e, label=below:$-a_{2,1}$, label=above:$A_{2,1}$] {};
\node[]          (g) [right=of f] {$\cdots$};
\node[rectangle] (h) [right=of g, label=below:$-a_{2,n_2'-2}$, label=above:$A_{2,n_2'-2}$] {};
\node[circle] (i) [right=of h, label=below:$-a_{2,n_2'-1}$, label=above:$A_{2,n_2'-1}$] {};
\node[circle] (j) [right=of i, label=below:$-a_{2,n_2'}$, label=above:$A_{2,n_2'}$] {};

\draw[-] (a) -- (b);
\draw[-] (b) -- (c);
\draw[-] (c) -- (d);
\draw[-] (d) -- (e);
\draw[-] (e) -- (f);
\draw[-] (f) -- (g);
\draw[-] (g) -- (h);
\draw[-] (h) -- (i);
\draw[-] (i) -- (j);
\end{tikzpicture}$$
with a degeneration $C_{1,m_1-1}^+ = C_{1,m_1-1} + A_{2,n_2'-1}$. We continue until we obtain
$$\begin{tikzpicture}
  [inner sep=1mm,
R/.style={circle,draw=black!255,fill=white!20,thick},
T/.style={circle,draw=red!255,fill=white!255,thick}]

\node[circle] (c) [label=below:$-a_{1,p}-1$, label=above:$A_{1,p}$] {};
\node[]       (d) [right=of c] {$\cdots$};
\node[circle] (e) [right=of d, label=below:$-d$, , label=above:$A_c$] {};
\node[circle] (f) [right=of e, label=below:$-a_{2,1}$, label=above:$A_{2,1}$] {};
\node[]       (g) [right=of f] {$\cdots$};
\node[circle] (h) [right=of g, label=below:$-a_{2,n_2'}$, label=above:$A_{2,n_2'}$] {};

\draw[-] (c) -- (d);
\draw[-] (d) -- (e);
\draw[-] (e) -- (f);
\draw[-] (f) -- (g);
\draw[-] (g) -- (h);
\end{tikzpicture}$$
with degenerations $C_{1,m_1}^+ = C_{1,m_1} + A_{2,n_2'}$, $\dots$, $C_{1,m_1-\square}^+ = C_{1,m_1-\square} + A_c$. There are more degenerations but we only consider them. Note that the degeneration occurs from the decorated curve that the second subscript is largest.

\begin{figure} 
\centering
\begin{minipage}{.3\textwidth}
\begin{tikzpicture}
[inner sep=1mm,
R/.style={circle,draw=black!255,fill=white!20,thick},
T/.style={circle,draw=red!255,fill=white!255,thick}]
    \node[R] (a) [label={[yshift=.3cm]:$-d+2$}]{};
    \node[R] (b) [above left=.1cm and 1cm of a, label=left:$C_{1,e}$]{};
    \node    (c) [left=.9cm of a]{$\vdots$};
    \node[R] (d) [below left=.5cm and 1cm of a, label=left:$C_{1,m_1}$]{};
    \node[T] (e) [below left=.7cm and .5cm of a]{};
    \node    (f) [below=.6cm of a]{$\cdots$};
    \node[T] (g) [below right=.7cm and .5cm of a]{};
    \node[R] (h) [below =.5cm of e, label=below:$D_1$]{};
    \node[R] (i) [below =.5cm of g, label=below:$D_{d-3}$]{};
    \node[R] (j) [above right=.1cm and 1cm of a, label=right:$C_{2,1}$]{};
    \node    (k) [right=.9cm of a]{$\vdots$};
    \node[R] (l) [below right=.5cm and 1cm of a, label=right:$C_{2,m_2}$]{};

    \node[R] (m) [above left=.7cm and .7cm of a, label=above:$C_{1,e-1}$]{};
    \node[R] (n) [above right=.7cm and .7cm of a, label=above:$C_{1,1}$]{};
    \node    (o) [above=.6cm of a]{$\cdots$};

    \draw[-] (a) -- (b);
    \draw[-] (a) -- (d);
    \draw[-] (a) -- (e);
    \draw[-] (a) -- (g);
    \draw[-] (e) -- (h);
    \draw[-] (i) -- (g);
    \draw[-] (a) -- (j);
    \draw[-] (a) -- (l);
    \draw[-] (a) -- (m);
    \draw[-] (a) -- (n);

    \draw[edge,line width=10pt] (a) -- (b);
    \draw[edge,line width=10pt] (a) -- (d);
\end{tikzpicture}
\end{minipage}
$\rightarrow$
\begin{minipage}{.3\textwidth}
\begin{tikzpicture}
[inner sep=1mm,
R/.style={circle,draw=black!255,fill=white!20,thick},
T/.style={circle,draw=red!255,fill=white!255,thick}]
    \node[R] (a) [label={[yshift=.3cm]:$-1$}]{};
    \node[R] (b) [above left=.1cm and 1cm of a, label=left:$C_{1,1}$]{};
    \node    (c) [left=.9cm of a]{$\vdots$};
    \node[R] (d) [below left=.5cm and 1cm of a, label=left:$C_{1,m_1}$]{};
    \node[R] (e) [below left=.7cm and .5cm of a, label=below:$D_1$]{};
    \node    (f) [below=.6cm of a]{$\cdots$};
    \node[R] (g) [below right=.7cm and .5cm of a, label=below:$D_{d-3}$]{};
    \node[R] (j) [above right=.1cm and 1cm of a, label=right:$C_{2,1}$]{};
    \node    (k) [right=.9cm of a]{$\vdots$};
    \node[R] (l) [below right=.5cm and 1cm of a, label=right:$C_{2,m_2}$]{};

    \node[R] (m) [above left=.7cm and .7cm of a, label=above:$C_{1,e-1}$]{};
    \node[R] (n) [above right=.7cm and .7cm of a, label=above:$C_{1,1}$]{};
    \node    (o) [above=.6cm of a]{$\cdots$};

    \draw[-] (a) -- (b);
    \draw[-] (a) -- (d);
    \draw[-] (a) -- (e);
    \draw[-] (a) -- (g);
    \draw[-] (a) -- (j);
    \draw[-] (a) -- (l);
    \draw[-] (a) -- (m);
    \draw[-] (a) -- (n);

    \draw[edge,line width=10pt] (a) -- (b);
    \draw[edge,line width=10pt] (a) -- (d);
\end{tikzpicture}
\end{minipage}

\caption{Final divisorial contraction}
\label{fig:Final divisorial contraction2}
\end{figure}
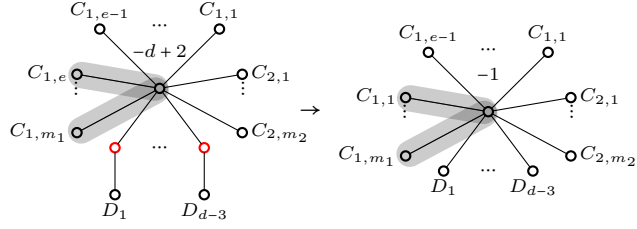

In our case, there are more curves on the singularity. That is, we have
$$\begin{tikzpicture}
  [inner sep=1mm,
R/.style={circle,draw=black!255,fill=white!20,thick},
T/.style={circle,draw=red!255,fill=white!255,thick}]
\node[circle] (a) [label=below:$-a_{1,n_1}$, label=above:$A_{1,n_1}$] {};
\node[]          (b) [right=of a] {$\cdots$};
\node[rectangle] (c) [right=of b, label=below:$-a_{1,n_1'}$, label=above:$A_{1,n_1'}$] {};
\node[]          (d) [right=of c] {$\cdots$};
\node[rectangle] (e) [right=of d, label=below:$-d$, , label=above:$A_c$] {};
\node[]          (f) [right=of e] {$\cdots$};
\node[rectangle] (g) [right=of f, label=below:$-a_{2,n_2'}$, label=above:$A_{2,n_2'}$] {};
\node[]          (h) [right=of g] {$\cdots$};
\node[circle]    (i) [right=of h, label=below:$-a_{2,n_2}$, label=above:$A_{2,n_2}$] {};

\draw[-] (a) -- (b);
\draw[-] (b) -- (c);
\draw[-] (c) -- (d);
\draw[-] (d) -- (e);
\draw[-] (e) -- (f);
\draw[-] (f) -- (g);
\draw[-] (g) -- (h);
\draw[-] (h) -- (i);
\end{tikzpicture}$$
.
But since we apply the MMP algorithm step by step from the left, we may assume that the degeneration occurs from the decorated curves on the first branch that the second subscript is the largest. Therefore We may assume that $C_{1,m_1}, \cdots, C_{1,e}$ be the decorated curves that their degeneration contains the central curve $A_c$ for $1 \leq e \leq m_1$. Then after flips and divisorial contractions, we arrive at a divisorial contraction in Figure~\ref{fig:Final divisorial contraction2}. The gray means a degeneration.

Each $(-1)$-curve in the left one in Figure~\ref{fig:Final divisorial contraction2} connects all $C_{1,m_1}, \dots C_{1,e}$ and $D_k$ for $k = 1, \dots, d-3$. Therefore there are points $p_k$ that $C_{1,e}, \dots, C_{1,m_1}$ and $D_k$ pass through $p_k$ for $k = 1, \dots, d-3$ respectively. After divisorial contractions, we obtain the right one in Figure~\ref{fig:Final divisorial contraction2} and the $(-1)$-curve connects all decorated curves except $C_{1,m_1}, \dots, C_{1,e}$. We denote the corresponding point as $p_0$. As an incidence matrix, we obtain Figure~\ref{fig:Type 2}.

\begin{figure} 
\centering
$\begin{bmatrix}[c|ccccccccccccccccccc]
                           & p_0    & p_1    & \cdots & p_{d-3}&         &        &        & q_1    & \cdots  & q_g   \\
               \cmidrule(lr){1-17}
               C_{1,1}   & 1      & \;     & \;     & \;     & *       & \cdots & *      & \;     & \;      & \;     &        & \;     & \;     & *      & \cdots & *      \\
               \vdots      & \vdots & \;     & \;     & \;     & \vdots  & \ddots & \vdots & \;     & \;      & \;     &        & \;     & \;     & \vdots & \ddots & \vdots \\
               C_{1,e-1}   & 1      & \;     & \;     & \;     & *       & \cdots & *      & \;     & \;      & \;     &        & \;     & \;     & *      & \cdots & *      \\
               \cmidrule(lr){1-17}
               C_{1,e}     & \;     & 1      & \cdots & 1      & *       & \cdots & *      & *      & \cdots  &  *     & *      & \cdots & *      & \;     & \;     & \;     \\    
               \vdots      & \;     & \vdots & \ddots & \vdots & \vdots  & \ddots & \vdots & \vdots & \ddots  & \vdots & \vdots & \ddots & \vdots & \;     & \;     & \;     \\
               C_{1,m_1}     & \;     & 1      & \cdots & 1      & *       & \cdots & *      & *      & \cdots  &  *     & *      & \cdots & *      & \;     & \;     & \;     \\   
               \cmidrule(lr){1-17}

               C_{2,1}     & 1      & \;     & \;     & \;     & \;      & \;     & \;     & *      & \cdots  & *      &        & \;     & \;     & \;     & \;     & \;     \\     
               \vdots      & \vdots & \;     & \;     & \;     & \;      & \;     & \;     & \vdots & \ddots  & \vdots &        & \;     & \;     & \;     & \;     & \;     \\
               C_{2,m_2}   & 1      & \;     & \;     & \;     & \;      & \;     & \;     & *      & \cdots  & *      &        & \;     & \;     & \;     & \;     & \;     \\
               \cmidrule(lr){1-17}
               D_1         & 1      & 1      & 0      & 0      & \;      & \;     & \;     & \;     & \;      & \;     &        & \;     & \;     & \;     & \;     &        \\    
               \vdots      & \vdots & 0      & \ddots & 0      & \;      & \;     & \;     & \;     & \;      & \;     &        & \;     & \;     & \;     & \;     &        \\  
               D_{d-3}     & 1      & 0      & 0      & 1      & \;      & \;     & \;     & \;     & \;      & \;     &        & \;     & \;     & \;     & \;     &
\end{bmatrix}$
\caption{Type 2}
\label{fig:Type 2}
\end{figure}
In the aspect of a combinatorial incidence matrix, we know that there are points for the intersection relations $C_{1,j'}.C_{2,j''} = 1$. We denote these points as $q_1, \cdots, q_g$. The index $g$ depends on the P-resolution.

\begin{figure} 
\centering
\begin{minipage}{.4\textwidth}
\begin{tikzpicture}
  [inner sep=1mm,
R/.style={circle,draw=black!255,fill=white!20,thick},
T/.style={circle,draw=red!255,fill=white!255,thick}]
    \node[R] (a) [label=below:$-d$, , label=above:$A_c$]{};
    \node    (b) [right=.5cm of a] {$\cdots$};
    \node[R] (c) [right=.5cm of b, label=below:$-a_{2, n_2'-1}$, label=above:$A_{2,n_2'-1}$]{};
    \node[R] (d) [right=1cm of c, label=below:$-a_{2, n_2'}$, label=above:$A_{2,n_2'}$]{};
    \node[T] (e) [right=.5cm of d]{};
    \node[R] (f) [above right=.3cm and .3cm of e, label=above:$D_{2,j}$]{};
    \node    (g) [right=.5cm of e]{$\cdots$};

    \draw[-] (a)--(b);
    \draw[-] (b)--(c);
    \draw[-] (c)--(d);
    \draw[-] (d)--(e);
    \draw[-] (e)--(f);
    \draw[-] (e)--(g);

    \draw[edge,line width=10pt] (a) -- (d);
    \draw[edge,line width=10pt] (a) -- (c);
\end{tikzpicture}
\end{minipage}
$\rightarrow$
\begin{minipage}{.4\textwidth}
\begin{tikzpicture}
  [inner sep=1mm,
R/.style={circle,draw=black!255,fill=white!20,thick},
T/.style={circle,draw=red!255,fill=white!255,thick}]
    \node[R] (a) [label=below:$-d$, , label=above:$A_c$]{};
    \node    (b) [right=.5cm of a] {$\cdots$};
    \node[R] (c) [right=.5cm of b, label=below:$-a_{2, n_2'-1}$, label=above:$A_{2,n_2'-1}$]{};
    \node[R] (d) [right=1cm of c, label=below:$-a_{2, n_2'}$, label=above:$A_{2,n_2'}$]{};
    \node[T] (e) [right=.5cm of d]{};
    \node[R] (f) [above right=.3cm and .3cm of e, label=above:$D_{2,j}$]{};
    \node    (g) [right=.5cm of e]{$\cdots$};

    \draw[-] (a)--(b);
    \draw[-] (b)--(c);
    \draw[-] (c)--(d);
    \draw[-] (d)--(e);
    \draw[-] (e)--(f);
    \draw[-] (e)--(g);

    \draw[edge,line width=10pt] (a) -- (e);
    \draw[edge,line width=10pt] (a) -- (d);
    \draw[edge,line width=10pt] (a) -- (c);
\end{tikzpicture}
\end{minipage}
\caption{$n_2' \geq 1$}
\label{fig:n2<=1}
\end{figure}
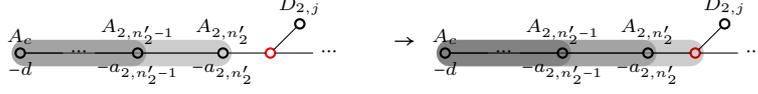

In principle, a decorated curve $C_{1,j}$ on the first branch and $C_{2,j'}$ on the second branch cannot connected through a $(-1)$-curve except the $(-1)$-curve that comes from the central curve. Another possible case is that $C_{1,j}$ degenerates to a curve of the second branch. In our case, the decorated curves $C_{1,m_1}, \dots, C_{1,e}$ are the case.
Let 
$$
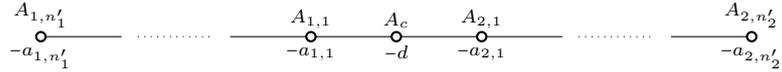
\begin{tikzpicture}
  [inner sep=1mm,
R/.style={circle,draw=black!255,fill=white!20,thick},
T/.style={circle,draw=black!255,fill=white!255,thick}]
    \node[T] (a) [label=below:$-a_{1,n_1'}$, , label=above:$A_{1,n_1'}$]{};
    \node    (b) [right=of a] {};

    \node    (h) [right=of b] {};

    \node[T] (c) [right=of h, label=below:$-a_{1,1}$, label=above:$A_{1,1}$] {};

    \node[T] (d) [right=of c, label=below:$-d$, label=above:$A_c$] {};

    \node[T] (e) [right=of d, label=below:$-a_{2,1}$, label=above:$A_{2,1}$] {};
    \node    (f) [right=of e] {};

    \node    (i) [right=of f] {};

    \node[T] (g) [right=of i, label=below:$-a_{2,n_2'}$, , label=above:$A_{2,n_2'}$] {};

    \draw[-] (a)--(b);
    \draw[dotted] (b)--(h);
    \draw[-] (h)--(c);
    \draw[-] (c)--(d);
    \draw[-] (d)--(e);
    \draw[-] (e)--(f);
    \draw[dotted] (f)--(i);
    \draw[-] (i)--(g);
\end{tikzpicture}$$
 be the dual resolution graph of a Wahl singularity where $n_1' < n_1$ and $n_2' < n_2$ as in lemma~\ref{flips on Tsing contaiing the central}. If $n_2' \geq 1$, then we obtain Figure~\ref{fig:n2<=1} during the MMP algorithm.
If decorated curves $C_{1,m_1}, \cdots, C_{1,e_1}$($e \leq e_1 \leq m_1$) degenerate to $A_{2,n_2'}$, then the $(-1)$-curve in the left in Figure~\ref{fig:n2<=1} connects all $C_{1,m_1}, \cdots, C_{1,e_1}$ and a decorated curve in the second branch. After the divisorial contraction, if the decorated curves $C_{1,e_1-1}, \cdots, C_{1, e_2}$($e \leq e_2 \leq e_1$) degenerate to $C_{2,n_2''}$, then the $(-1)$-curve in the right connects the decorated curves $C_{1,e_1-1}, \cdots, C_{1, e_2}$ and the decorated curve in the second branch. This $(-1)$-curve is not connected to any of the decorated curves of $C_{1,m_1}, \cdots, C_{1,e_1}$. This process continues until the central curve becomes a $(-1)$-curve. 

Therefore if we let $\mathcal{C}_j$ be the set of decorated curves that degenerate to $A_{2,j}$ but not $A_{2,j+1}$ for $j = 1, \cdots n_2'$, then$ \{\mathcal{C}_j\}_{j=1}^{g'}$ is a partition of $\{C_{1,m_1}, \cdots, C_{1, e}\}$. By the above observation, there are $g'$ number of points $q_1, \cdots, q_{g'}$ that all decorated curves in $\mathcal{C}_j$ pass through only $q_j$. In the aspect of an incidence matrix, there are stair-shaped sub-matrix of the Figure~\ref{fig:Stair-shaped sub-matrix} in the columns $q_1, \cdots, q_g$.
\begin{figure} 
\centering
$\begin{bmatrix}
q_1             & q_2    & q_3    & \cdots & q_{g'} \\
\cmidrule(lr){1-5}
1      & 0      & 0      & \cdots & 0      \\
\vdots & \vdots & \vdots & \cdots & 0      \\
1      & 0      & 0      & \cdots & 0      \\
0      & 1      & 0      & \cdots & 0      \\
0      & \vdots & \vdots & \cdots & 0      \\ 
0      & 1      & 0      & \cdots & 0      \\
\vdots & \ddots & \ddots & \ddots & \vdots \\
\vdots & \ddots & \ddots & \ddots & 1      \\
\vdots & \ddots & \ddots & \ddots & \vdots \\
0      & \cdots & \cdots & \cdots & 1
\end{bmatrix}$
\caption{Stair-shaped sub-matrix}
\label{fig:Stair-shaped sub-matrix}
\end{figure}

If $n_2' = 0$, after the MMP algorithm on the second branch, we obtain Figure~\ref{fig:$n_2' = 0$}. From the $(-1)$-curve, we know that the decorated curves $C_{1,m_1}, \dots, C_{1,e}$ and $C_{2,1}, \dots, C_{2,m_2}$ pass through a point $q_1$. Therefore the column $q_1$ consists of $1$ for rows $C_{1,m_1}, \dots, C_{1,e}$ and $C_{2,1}, \dots, C_{2,m_2}$.
\begin{figure} 
\centering
\begin{tikzpicture}
  [inner sep=1mm,
R/.style={circle,draw=black!255,fill=white!20,thick},
T/.style={circle,draw=red!255,fill=white!255,thick}]
    \node    (b) [] {$\cdots$};
    \node[R] (c) [right=.5cm of b, label=below:$-a_{1,1}$, label=above:$A_{1,1}$]{};
    \node[R] (d) [right=1cm of c, label=below:$-d$, label=above:$A_c$]{};
    \node[T] (e) [right=.5cm of d]{};
    \node[R] (f) [above right=.3cm and .3cm of e, label=above:$D_{2,j}$]{};
    \node    (g) [right=.5cm of e]{$\cdots$};

    \draw[-] (b)--(c);
    \draw[-] (c)--(d);
    \draw[-] (d)--(e);
    \draw[-] (e)--(f);
    \draw[-] (e)--(g);

    \draw[edge,line width=7pt] (b) -- (d);
\end{tikzpicture}
\caption{$n_2' = 0$}
\label{fig:$n_2' = 0$}
\end{figure}
That is, every $(-1)$-curve that connected to the central curve connects all decorated curves $C_{1,m_1}, \cdots, C_{1,e}$. Therefore the entries 
of the columns $q_1, \cdots, q_g$ are all $1$.

\begin{lemma}[Incidence matrix of type 2-1 and 2-2]
Incidence matrices of type 2-1 and 2-2 are of the form shown in Figure~\ref{fig:Type 2}. In addition, type 2-1 contains a stair-shaped sub-matrix in the columns $q_1, \dots, q_{g'}$. For type 2-2, the column $q_1$ consists of $1$
\label{lem:configurationoftype2-1andtype2-2}
\end{lemma}

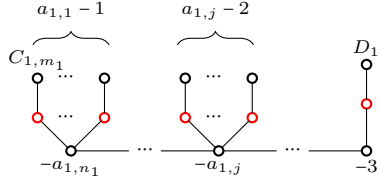
\begin{figure} 
\centering
\begin{tikzpicture}
  [inner sep=1mm,
R/.style={circle,draw=black!255,fill=white!20,thick},
T/.style={circle,draw=red!255,fill=white!255,thick}]
    \node[R] (a) [label=below:$-a_{1,n_1}$]{};
    \node    (b) [right=.7cm of a]{$\cdots$};
    \node[R] (c) [right=.7cm of b, label=below:$-a_{1,j}$]{};
    \node    (d) [right=.7cm of c]{$\cdots$};
    \node[R] (e) [right=.7cm of d, label=below:$-3$]{};

    \node[T] (j) [above left=.5cm of a]{};
    \node    (k) [right=.1cm of j]{$\cdots$};
    \node[T] (l) [above right=.5cm of a]{};
    \node[R] (m) [above=.4cm of j, label=above:$C_{1,m_1}$]{};
    \node    (k) [right=.1cm of m]{$\cdots$};
    \node[R] (n) [above=.4cm of l]{};

    \node[T] (o) [above left=.5cm of c]{};
    \node    (p) [right=.1cm of o]{$\cdots$};
    \node[T] (q) [above right=.5cm of c]{};
    \node[R] (r) [above=.4cm of o]{};
    \node    (k) [right=.1cm of r]{$\cdots$};
    \node[R] (s) [above=.4cm of q]{};

    \node[T] (ad) [above=.5cm of e]{};
    \node[R] (ag) [above=.4cm of ad, label=above:$D_{1}$]{};

    \draw[-] (a)--(b);
    \draw[-] (b)--(c);
    \draw[-] (c)--(d);
    \draw[-] (d)--(e);

    \draw[-] (o)--(c);
    \draw[-] (q)--(c);
    \draw[-] (o)--(r);
    \draw[-] (s)--(q);

    \draw[-] (a)--(j);
    \draw[-] (a)--(l);
    \draw[-] (j)--(m);
    \draw[-] (l)--(n);

    \draw[-] (ad)--(e);
    \draw[-] (ag)--(ad);

\draw [decorate,decoration={brace,amplitude=5pt,raise=4ex}]
  (-0.5,.8) -- (0.5,.8) node[midway,yshift=3em]{$a_{1,1}-1$};

\draw [decorate,decoration={brace,amplitude=5pt,raise=4ex}]
  (1.4,.8) -- (2.4,.8) node[midway,yshift=3em]{$a_{1,j}-2$};

\end{tikzpicture}
\caption{Resolution graph of $[D_1,C_1]$}
\label{fig:Resolution graph of D_1;C1}
\end{figure}

For the type 2-1, we can show that the stair-shaped sub-matrix is unique in the given incidence matrix. We consider the sub-matrix $[D_1, C_1]$(Figure~\ref{fig:D1;C1}).
\begin{figure} 
\centering
$\begin{bmatrix}[c|ccccccccccccccccccc]
                           & p_0    & p_1    & \cdots & p_{d-3}&         &        &        & q_1    & \cdots  & q_g   \\
              \cmidrule(lr){1-17}
               D_1         & 1      & 1      &        &        & \;      & \;     & \;     & \;     & \;      & \;     &        & \;     & \;     & \;     & \;     &        \\    
              \cmidrule(lr){1-17}
               C_{1,1}   & 1      & \;     & \;     & \;     & *       & \cdots & *      & \;     & \;      & \;     &        & \;     & \;     & *      & \cdots & *      \\
               \vdots      & \vdots & \;     & \;     & \;     & \vdots  & \ddots & \vdots & \;     & \;      & \;     &        & \;     & \;     & \vdots & \ddots & \vdots \\
               C_{1,e-1}   & 1      & \;     & \;     & \;     & *       & \cdots & *      & \;     & \;      & \;     &        & \;     & \;     & *      & \cdots & *      \\
              \cmidrule(lr){1-17}
               C_{1,e}     & \;     & 1      & \cdots & 1      & *       & \cdots & *      & *      & \cdots  &  *     & *      & \cdots & *      & \;     & \;     & \;     \\  
               \vdots      & \;     & \vdots & \ddots & \vdots & \vdots  & \ddots & \vdots & \vdots & \ddots  & \vdots & \vdots & \ddots & \vdots & \;     & \;     & \;     \\
               C_{1,m_1}     & \;     & 1      & \cdots & 1      & *       & \cdots & *      & *      & \cdots  &  *     & *      & \cdots & *      & \;     & \;     & \;     \\   
\end{bmatrix}$
\caption{$[D_1,C_1]$}
\label{fig:D1;C1}
\end{figure}
In Figure~\ref{fig:sandwichedsatructureofcqss}, if we only contract the $(-1)$-curves on the first branch and the central curve, then we know that the first branch is contracted and the central curve becomes a $(-d+2)$-curve. Therefore if we ignore the second branch and if $d = 3$, then only by the $(-1)$-curves that we mentioned now, the graph will be contracted to a smooth point. See the Figure~\ref{fig:Resolution graph of D_1;C1}. We consider $[D_1, C_1]$ as an incidence matrix of Figure~\ref{fig:Resolution graph of D_1;C1}.

Let $[2, b_{1,1}, \cdots, b_{1,m_1}]$ be the dual H-J continued fraction of $[3, a_{1,1}, \cdots, a_{1,n_1}]$. By Proposition~\ref{Nemethi}, there is an integer sequence $\underline{k}$ and triangulation $\theta$ that generate the incidence matrix $[D_1,C_1]$. Here, we assign the vertices of the convex $(m_1 + 2)$-gon to $d, C_{1,1}, \dots, C_{1,m_1}, N$ counterclockwise. 

\begin{figure} 
\centering
\begin{tikzpicture}
  [inner sep=1mm,
R/.style={circle,draw=black!255,fill=white!20,thick},
T/.style={circle,draw=red!255,fill=white!255,thick}]

    \node[R] (a) [label=left:$C_{1,e+1}$]{};
    \node[R] (b) [above right=.7cm and .3cm of a, label=left:$C_{1,e}$]{};
    \node[R] (c) [above right=.5cm and .7cm of b, label=left:$C_{1,e-1}$]{};
    \node    (d) [above right=.5cm of c]{$\cdots$};
    \node[R] (e) [below right=.5cm of d, label=right:$C_{1,2}$]{};
    \node[R] (f) [below right=.7cm of e, label=right:$C_{1,1}$]{};
    \node[R] (g) [below right=.7cm and .3cm of f, label=right:$d$]{};
    \node[R] (h) [below left=.7cm and .3cm of g, label=right:$N$]{};
    \node[R] (i) [below left=.7cm and .7cm of h, label=below:$C_{1,m_1}$]{};
    \node    (j) [below right=.7cm and .3cm of a]{$\ddots$};

    \draw[-] (a)--(b);
    \draw[-] (b)--(c);
    \draw[-] (c)--(d);
    \draw[-] (d)--(e);        
    \draw[-] (e)--(f);
    \draw[-] (f)--(g);
    \draw[-] (g)--(h);
    \draw[-] (h)--(i);
    \draw[-] (b)--(f);
    \draw[-] (b)--(g);
    \draw[-] (b)--(h);
    \draw[-] (a)--(j);
\begin{scope}[on background layer] 
    \fill[pattern=north east lines, pattern color=red] (b.center) -- (f.center) -- (g.center) -- cycle;
    \fill[pattern=north west lines, pattern color=blue] (b.center) -- (g.center) -- (h.center) -- cycle;  
\end{scope}
\end{tikzpicture}
\\
$\begin{bmatrix}[c|cc]
                           & p_0    & p_1   \\
              \cmidrule(lr){1-3}
               D_1         & 1     & 1     \\    
              \cmidrule(lr){1-3}
               C_{1,1}     & \;     & -1    \\
               \vdots      & \;     & \;    \\
               C_{1,e-1}     & \;     & \;    \\
              \cmidrule(lr){1-3}
               C_{1,e}   & -1     & 1     \\    
               \vdots      & \;     &       \\
               C_{1,m_1}   & \;     &       \\   
\end{bmatrix}$
$\xrightarrow{\int}$
$\begin{bmatrix}[c|ccccccccccccccccccc]
                           & p_0    & p_1    \\
              \cmidrule(lr){1-3}
               D_1         & 1      & 1      \\    
              \cmidrule(lr){1-3}
               C_{1,1}     & 1      & \;     \\
               \vdots      & \vdots & \;     \\
               C_{1,e-1}     & 1      & \;     \\
              \cmidrule(lr){1-3}
               C_{1,e}   & \;     & 1      \\    
               \vdots      & \;     & \vdots \\
               C_{1,m_1}   & \;     & 1      \\   
\end{bmatrix}$
\caption{triangles that induce columns $p_0$ and $p_1$}
\label{triangles that induce columns p0 and p1}
\end{figure}
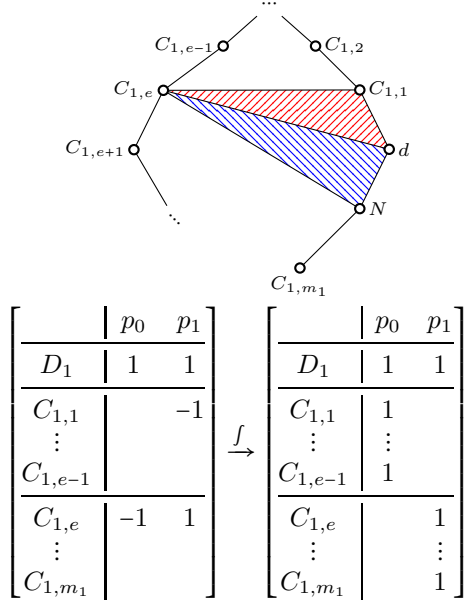

\begin{figure}
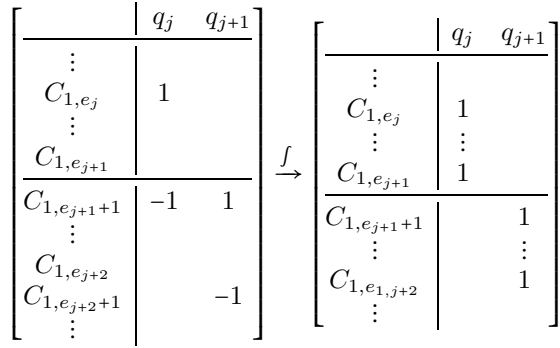
 
\centering
$\begin{bmatrix}[c|cc]
                                 & q_j    & q_{j+1}    \\
               \cmidrule(lr){1-3}
               \vdots          &        &        \\
               C_{1,e_j}       & 1      & \;     \\
               \vdots          &        & \;     \\
               C_{1,e_{j+1}}   &        & \;     \\
               \cmidrule(lr){1-3}
               C_{1,e_{j+1}+1} & -1     & 1      \\    
               \vdots          & \;     &        \\
               C_{1,e_{j+2}}   & \;     &       \\
               C_{1,e_{j+2}+1} & \;     & -1      \\
               \vdots          & \;     &      \\   
\end{bmatrix}$
$\xrightarrow{\int}$
$\begin{bmatrix}[c|cc]
                              & q_j    & q_{j+1}    \\
               \cmidrule(lr){1-3}
               \vdots          &        &        \\
               C_{1,e_j}       & 1      & \;     \\
               \vdots          & \vdots & \;     \\
               C_{1,e_{j+1}}   & 1      & \;     \\
               \cmidrule(lr){1-3}
               C_{1,e_{j+1}+1} & \;     & 1      \\    
               \vdots          & \;     & \vdots \\
               C_{1,e_{1,j+2}} & \;     & 1      \\
               \vdots          & \;     &       \\    
\end{bmatrix}$
\caption{Two triangles induce a 'stair'}
\label{fig:Two triangles induce a 'stair'}
\end{figure}

Since the matrix in Figure~\ref{fig:D1;C1} contains two columns, $p_0$ and $p_1$, we can deduce the existence of two triangles, namely $\triangle(d, C_{1,1}, C_{1,e})$ and $\triangle(d,C_{1,e},N)$, as shown in Figure~\ref{triangles that induce columns p0 and p1}.
Due to the diagonal $\overline{C_{1,e},N}$, there exist triangles $\triangle(N,C_{1,e},C_{1,e_1}), \triangle(N,C_{1,e_1},C_{1,e_2}), \ldots, \triangle(N,C_{1,e_o},C_{1,m_1})$ for some $e < e_1 < \cdots <e_o \leq m_1$. Two consecutive triangles $\triangle(N,C_{1,e_j},N,C_{1,e_{j+1}})$ and $\triangle(N,C_{1,e_{j+1}},N,C_{1,e_{j+2}})$ form a 'stair' pattern, as depicted in Figure~\ref{fig:Two triangles induce a 'stair'}.
Therefore we find the stair-shaped sub-matrix in the triangulation method that we observed in Lemma~\ref{lem:configurationoftype2-1andtype2-2}. Furthermore, we see that the set of triangles that make the stair-shaped sub-matrix is the only one we found.

\section{Deformations of weighted homogeneous surface singularities with big central node}
  In this section, we introduce combinatorial incidence matrices, which are denoted as cases $A$ and $B$, of a weighted homogeneous surface singularity. Then, we prove that every combinatorial incidence matrix of a weighted homogeneous surface singularity with $d \geq t + 3$ is only one of the cases. And we construct $P$-resolutions only from the combinatorial information of the cases. Finally, we show that the constructed $P$-resolutions actually induce the given combinatorial incidence matrices.
\subsection{Weighted homogeneous surface singularities}
In this section, $(X, 0)$ is a weighted homogeneous surface singularity. The singularity $(X,0)$ is a two dimensional singularity with a good $\mathbb{C}^*$-action(Orlik-Wagreich \cite{Orl}).

The dual resolution graph of the singularity $(X, 0)$ is star-shaped. That is, there exist a central node of degree $-d$ and $t$-branches. Each branch is the dual resolution graph of a cyclic quotient surface singularity. Therefore we assign the singularity $(X, 0)$ to $(d, (n_1, q_1), \dots, (n_t, q_t))$ with $n_i / q_i = [a_{i,1}, \dots, a_{i, n_i}]$. We assume that the $a_{i,1}$-curve is connected to the central curve.

Assume that $d \geq t+1$. If we attach $(a_{i,n_i} - 1)$ $(-1)$-curve to $(a_{i,n_i})$-curve, $(a_{i,j}-2)$ $(-1)$-curve to $a_{i,j}$-curve for $j < n_i$ and $(d-t-1)$ $(-1)$-curve to the central curve, then the graph contracts to a smooth point(Refer Figure~\ref{fig:Sandwiched structure of X}). Therefore the graph is sandwiched and we obtain a sandwiched structure of the singularity by attaching decorated curves to the $(-1)$-curves.

Decorated curves connected to a curve of $i$-th branch through the $(-1)$-curve is labeled by $C_{i,j}$. The second sub script is labeled as in cyclic quotient surface singularities.

We frequently examine sub-matrices of a (combinatorial) incidence matrix $M$ that are composed of certain rows representing decorated curves. We indicate the sub-matrix that consists of decorated curves on the $i$-th branch as $M_i$. Furthermore, we use the notation $[M_i, M_j]$ for the sub-matrix that comprises decorated curves on both the $i$-th and $j$-th branches, despite it resembling the parallel sum of two matrices.

Let the combinatorial equations of the cyclic quotient surface singularity of $\frac{1}{n_i}(1, q_i)$ be 
$l(C_{i,j}) = a$ and $C_{i,j}.C_{i,j'} = b$. Then the combinatorial equations of the sandwiched structure of a weighted homogeneous surface singularity is 
\begin{equation}
\begin{split}
l(C_{i,j}) = a + 1 \\ l(D_k) = 2 \\ C_{i,j}.C_{i,j'}= b + 1 \\ C_{i,j}.C_{i',j'} = 1 \\ C_{i,j}.D_k = 1  
\end{split}
\label{eq:combinatorial equations of wighted homogeneous surface singularities}
\end{equation}
for $i, i' = 1, \dots, m_i$, $j, j' = 1, \dots, t$, $k = 1, \dots, d-t-1$, $i \neq i'$, $j \neq j'$.
The difference comes from the central curve. From this observation, we expect that a combinatorial incidence matrix of $(X, 0)$ contains an incidence matrix of a cyclic quotient surface singularity.

\subsection{Incidence matrices of weighted homogeneous surface singularities}
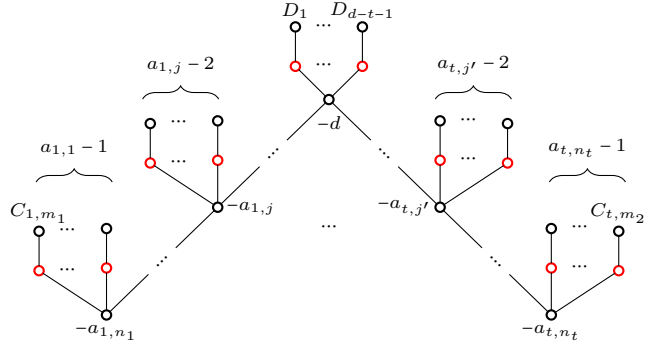
\begin{figure} 
\centering
\begin{tikzpicture}
  [inner sep=1mm,
R/.style={circle,draw=black!255,fill=white!20,thick},
T/.style={circle,draw=red!255,fill=white!255,thick}]
    \node[R] (e) [label={[yshift=-.5cm]:$-d$}]{};
    \node    (d) [below left=.7cm of e]{$\udots$};
    \node[R] (c) [below left=.7cm of d, label=right:$-a_{1,j}$]{};
    \node    (b) [below left=.7cm of c]{$\udots$};
    \node[R] (a) [below left=.7cm of b, label=below:$-a_{1,n_1}$]{};
    \node    (f) [below right=.7cm of e]{$\ddots$};
    \node[R] (g) [below right=.7cm of f, label=left:$-a_{t,j'}$]{};
    \node    (h) [below right=.7cm of g]{$\ddots$};
    \node[R] (i) [below right=.7cm of h, label=below:$-a_{t,n_t}$]{};

    \node[T] (j) [above left=.5cm and 0.8cm of a]{};
    \node    (k) [right=.1cm of j]{$\cdots$};
    \node[T] (l) [above= 0.5cm of a]{};
    \node[R] (m) [above=.4cm of j, label=above:$C_{1,m_1}$]{};
    \node    (k) [right=.1cm of m]{$\cdots$};
    \node[R] (n) [above=.4cm of l]{};

    \node[T] (o) [above left=.5cm and 0.8cm of c]{};
    \node    (p) [right=.1cm of o]{$\cdots$};
    \node[T] (q) [above=.5cm of c]{};
    \node[R] (r) [above=.4cm of o]{};
    \node    (k) [right=.1cm of r]{$\cdots$};
    \node[R] (s) [above=.4cm of q]{};

    \node[T] (t) [above=.5cm of g]{};
    \node    (u) [right=.1cm of t]{$\cdots$};
    \node[T] (v) [above right=.5cm and 0.8cm of g]{};
    \node[R] (w) [above=.4cm of t]{};
    \node    (k) [right=.1cm of w]{$\cdots$};
    \node[R] (x) [above=.4cm of v]{};

    \node[T] (y) [above=.5cm of i]{};
    \node    (z) [right=.1cm of y]{$\cdots$};
    \node[T] (aa) [above right=.5cm and 0.8cm of i]{};
    \node[R] (ab) [above=.4cm of y]{};
    \node    (k) [right=.1cm of ab]{$\cdots$};
    \node[R] (ac) [above=.4cm of aa, label=above:$C_{t,m_2}$]{};

    \node[T] (ad) [above left=.5cm of e]{};
    \node    (ae) [right=.1cm of ad]{$\cdots$};
    \node[T] (af) [above right=.5cm of e]{};
    \node[R] (ag) [above=.4cm of ad, label=above:$D_{1}$]{};
    \node    (k) [right=.1cm of ag]{$\cdots$};
    \node[R] (ah) [above=.4cm of af, label=above:$D_{d-t-1}$]{};

    \node    (CE) [below=1.5cm of e] {$\cdots$};

    \draw[-] (a)--(b);
    \draw[-] (b)--(c);
    \draw[-] (c)--(d);
    \draw[-] (d)--(e);
    \draw[-] (e)--(f);
    \draw[-] (f)--(g);
    \draw[-] (g)--(h);
    \draw[-] (h)--(i);
    
    \draw[-] (a)--(j);
    \draw[-] (a)--(l);
    \draw[-] (j)--(m);
    \draw[-] (l)--(n);

    \draw[-] (c)--(o);
    \draw[-] (c)--(q);
    \draw[-] (r)--(o);
    \draw[-] (s)--(q);

    \draw[-] (t)--(g);
    \draw[-] (v)--(g);
    \draw[-] (w)--(t);
    \draw[-] (x)--(v);

    \draw[-] (y)--(i);
    \draw[-] (aa)--(i);
    \draw[-] (ab)--(y);
    \draw[-] (ac)--(aa);

    \draw[-] (ad)--(e);
    \draw[-] (af)--(e);
    \draw[-] (ag)--(ad);
    \draw[-] (ah)--(af);

\draw [decorate,decoration={brace,amplitude=5pt,raise=4ex}]
  (-3.85,-1.7) -- (-2.85,-1.7) node[midway,yshift=3em]{$a_{1,1}-1$};

\draw [decorate,decoration={brace,amplitude=5pt,raise=4ex}]
  (-2.45,-0.6) -- (-1.45,-0.6) node[midway,yshift=3em]{$a_{1,j}-2$};

\draw [decorate,decoration={brace,amplitude=5pt,raise=4ex}]
  (1.4,.-0.6) -- (2.4,.-0.6) node[midway,yshift=3em]{$a_{t,j'}-2$};

\draw [decorate,decoration={brace,amplitude=5pt,raise=4ex}]
  (2.9,-1.7) -- (3.9,-1.7) node[midway,yshift=3em]{$a_{t,n_t}-1$};
\end{tikzpicture}
\caption{Sandwiched structure of $X$}
\label{fig:Sandwiched structure of X}
\end{figure}

In this subsection, we classify the combinatorial incidence matrices of $X$ based on its sandwiched structure, as illustrated in Figure~\ref{fig:Sandwiched structure of X}.

Let $M$ be a combinatorial incidence matrix of the singularity $X$. We define the sub-matrices of $M$ as $M_i$, which consists of the rows $C_{i,1}, \ldots, C_{i,m_i}$ for each $i = 1, \ldots, t$, and $D$, which consists of $D_1, \ldots, D_{d-t-1}$. We give a lemma about the sub-matrices.
\begin{lemma}
Let $M$ be a combinatorial incidence matrix of a singularity $(X,0)$ with the sandwiched structure as Figure~\ref{fig:Sandwiched structure of X}. Let $[M_i, M_j, D]$ be the sub-matrix of $M$ consisting of $M_i$, $M_j$ and $D$ for $i \neq j$. Then the sub-matrix $[M_i, M_j, D]$ is an incidence matrix of the cyclic quotient surface singularity $[a_{i,n_j}, \dots, a_{i,1}, d-t+2, a_{j,1}, \dots, a_{j,m_j}]$ with the sandwiched structure as in Figure~\ref{fig:Sandwiched structure on $[M_i, M_j, D]$}.
\begin{figure} 
\centering
\begin{tikzpicture}
  [inner sep=1mm,
R/.style={circle,draw=black!255,fill=white!20,thick},
T/.style={circle,draw=red!255,fill=white!255,thick}]
    \node[R] (a) [label=below:$-a_{i,n_i}$]{};
    \node    (b) [right=.7cm of a]{$\cdots$};
    \node[R] (c) [right=.7cm of b, label=below:$-a_{i,j}$]{};
    \node    (d) [right=.7cm of c]{$\cdots$};
    \node[R] (e) [right=.7cm of d, label=below:$-d+t-2$]{};
    \node    (f) [right=.7cm of e]{$\cdots$};
    \node[R] (g) [right=.7cm of f, label=below:$-a_{j,j'}$]{};
    \node    (h) [right=.7cm of g]{$\cdots$};
    \node[R] (i) [right=.7cm of h, label=below:$-a_{j,n_j}$]{};

    \node[T] (j) [above left=.5cm of a]{};
    \node    (k) [right=.1cm of j]{$\cdots$};
    \node[T] (l) [above right=.5cm of a]{};
    \node[R] (m) [above=.4cm of j, label=above:$C_{i,m_i}$]{};
    \node    (k) [right=.1cm of m]{$\cdots$};
    \node[R] (n) [above=.4cm of l]{};

    \node[T] (o) [above left=.5cm of c]{};
    \node    (p) [right=.1cm of o]{$\cdots$};
    \node[T] (q) [above right=.5cm of c]{};
    \node[R] (r) [above=.4cm of o]{};
    \node    (k) [right=.1cm of r]{$\cdots$};
    \node[R] (s) [above=.4cm of q]{};

    \node[T] (t) [above left=.5cm of g]{};
    \node    (u) [right=.1cm of t]{$\cdots$};
    \node[T] (v) [above right=.5cm of g]{};
    \node[R] (w) [above=.4cm of t]{};
    \node    (k) [right=.1cm of w]{$\cdots$};
    \node[R] (x) [above=.4cm of v]{};

    \node[T] (y) [above left=.5cm of i]{};
    \node    (z) [right=.1cm of y]{$\cdots$};
    \node[T] (aa) [above right=.5cm of i]{};
    \node[R] (ab) [above=.4cm of y]{};
    \node    (k) [right=.1cm of ab]{$\cdots$};
    \node[R] (ac) [above=.4cm of aa, label=above:$C_{j,m_j}$]{};

    \node[T] (ad) [above left=.5cm of e]{};
    \node    (ae) [right=.1cm of ad]{$\cdots$};
    \node[T] (af) [above right=.5cm of e]{};
    \node[R] (ag) [above=.4cm of ad, label=above:$D_{1}$]{};
    \node    (k) [right=.1cm of ag]{$\cdots$};
    \node[R] (ah) [above=.4cm of af, label=above:$D_{d-3}$]{};

    \draw[-] (a)--(b);
    \draw[-] (b)--(c);
    \draw[-] (c)--(d);
    \draw[-] (d)--(e);
    \draw[-] (e)--(f);
    \draw[-] (f)--(g);
    \draw[-] (g)--(h);
    \draw[-] (h)--(i);
    
    \draw[-] (a)--(j);
    \draw[-] (a)--(l);
    \draw[-] (j)--(m);
    \draw[-] (l)--(n);

    \draw[-] (c)--(o);
    \draw[-] (c)--(q);
    \draw[-] (r)--(o);
    \draw[-] (s)--(q);

    \draw[-] (t)--(g);
    \draw[-] (v)--(g);
    \draw[-] (w)--(t);
    \draw[-] (x)--(v);

    \draw[-] (y)--(i);
    \draw[-] (aa)--(i);
    \draw[-] (ab)--(y);
    \draw[-] (ac)--(aa);

    \draw[-] (ad)--(e);
    \draw[-] (af)--(e);
    \draw[-] (ag)--(ad);
    \draw[-] (ah)--(af);

\draw [decorate,decoration={brace,amplitude=5pt,raise=4ex}]
  (-0.5,.8) -- (0.5,.8) node[midway,yshift=3em]{$a_{i,1}-1$};

\draw [decorate,decoration={brace,amplitude=5pt,raise=4ex}]
  (1.4,.8) -- (2.4,.8) node[midway,yshift=3em]{$a_{i,j}-2$};

\draw [decorate,decoration={brace,amplitude=5pt,raise=4ex}]
  (5.3,.8) -- (6.3,.8) node[midway,yshift=3em]{$a_{j,j'}-2$};

\draw [decorate,decoration={brace,amplitude=5pt,raise=4ex}]
  (7.3,.8) -- (8.3,.8) node[midway,yshift=3em]{$a_{j,n_2}-1$};
\end{tikzpicture}
\caption{Sandwiched structure on $[M_i, M_j, D]$}
\label{fig:Sandwiched structure on $[M_i, M_j, D]$}
\end{figure}
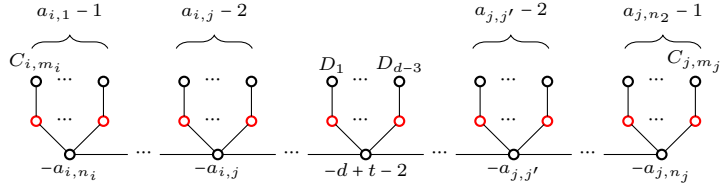
\label{lem:twobranches}
\end{lemma}
\begin{proof}
  Consider the combinatorial equations of $C_{i,1}, \cdots, C_{i, m_i}, C_{j, 1}, \cdots, C_{j, m_j}, D_1, \cdots, D_{d-t-1}$ that obtained from the sandwiched structure of Figure~\ref{fig:Sandwiched structure on $[M_i, M_j, D]$}. It is actually the same with the equations that obtained from Figure~\ref{fig:Sandwiched structure of X}. Therefore the sub-matrix $[M_i, M_j, D]$ satisfies the equations.
\end{proof}

\begin{lemma}
The sub matrix $D$ is a $(d-t-1) \times (d-t)$ matrix
$$\begin{bmatrix}[cccc]
1      & 1 &        &   \\
\vdots &   & \ddots &   \\
1      &   &        & 1
\end{bmatrix}$$ where entries of the first column are all $1$ and the rest is the $(d-t-1) \times (d-t-1)$ identity matrix.
\end{lemma}
\begin{proof}
Decorated curves $D_1, \dots, D_{d-t-1}$ satisfy the following equations.
\begin{equation}
\begin{split}
l(D_k) = 2 \\
D_k.D_l = 1
\end{split}
\label{eq:intersectionrelationofD}
\end{equation}
for all $k, l = 1, \dots, d-t-1$ and $k \neq l$.
Since $d \geq t +3$, we have at least two decorated curves denoted by $D_k$. 

If $d-t-1 = 2$, then we have two decorated curves $D_1$ and $D_2$. There exists only one matrix that satisfies Equation~\ref{eq:intersectionrelationofD}.
\begin{equation}
\begin{bmatrix}[c|cccc]
    & p_0 & p_1 & p_2 \\
\cmidrule(lr){1-4}
D_1 & 1   & 1   &     \\
D_2 & 1   &     & 1
\end{bmatrix}
\label{mat:d-t=2}
\end{equation}

If $d-t-1 = 3$, we have three decorated curves $D_1$, $D_2$ and $D_3$. We can find matrices satisfying Equations~\ref{eq:intersectionrelationofD} by adding the decorated curve $D_3$ to the matrix~\ref{mat:d-t=2} to satisfy Equation~\ref{eq:intersectionrelationofD}. We have two such matrices.
\begin{equation}
\begin{bmatrix}[c|cccc]
    & p_0 & p_1 & p_3 & p_4 \\
\cmidrule(lr){1-5}
D_1 & 1   & 1   &     &   \\
D_2 & 1   &     & 1   &   \\
D_3 & 1   &     &     & 1
\end{bmatrix}
\;\;\;
\begin{bmatrix}[c|cccc]
    & p_0 & p_1 & p_2 \\
\cmidrule(lr){1-4}
D_1 & 1   & 1   &    \\
D_2 & 1   &     & 1  \\
D_3 &     & 1   & 1 
\end{bmatrix}
\label{mat:d-t=3}
\end{equation}
We will now show that the right sub-matrix cannot be a valid sub-matrix of $M$. Consider the decorated curve $C_{1,1}$ in the first branch. From Equation~\ref{eq:combinatorial equations of wighted homogeneous surface singularities}, we have the intersection relation $C_{1,1}.D_i = 1$ for $i = 1, 2, 3$. Suppose $C_{1,1}$ intersects $D_1$ at $p_0$. Then $C_{1,1}$ must also intersect $D_2$ at the same point. To satisfy the intersection relation, $C_{1,1}$ must intersect $D_3$ at $p_1$ or $p_2$. However, this causes $C_{1,1}$ to intersect $D_2$ one more time than it should, violating the intersection relation. Therefore, this sub-matrix cannot appear in any combinatorial incidence matrix.

If $d-t-1 \geq 4$, we have only one choice from the left one of matrices~\ref{mat:d-t=3}
\begin{equation}
\begin{bmatrix}[c|ccccc]
          & p_0    & p_1 & \cdots & p_{d-t-1}  \\
\cmidrule(lr){1-5}
D_1       & 1      & 1   &        &    \\
\vdots    & \vdots &     & \ddots &    \\
D_{d-t-1} & 1      &     &        & 1
\end{bmatrix}
\label{mat:generalD}
\end{equation}
that the lemma claimed.
\end{proof}
Denote the intersection point of all $D_i$(the first column of \ref{mat:generalD} ) as $p_0$ and the others as $p_1, \cdots, p_{d-t-1}$ in the matrix $M$.
\begin{theorem}
\label{thm:caseofcombinatorialincidencematrix}
Every combinatorial incidence matrix of $X$ can be classified into two cases.

$\mathbf{Case~A}$. All entries of the column $p_0$ are $1$. The rest consists of block sub-matrices as follows.
\begin{equation}
\begin{minipage}{.8\textwidth}
\resizebox{\textwidth}{!}
{$\begin{bmatrix}[c|ccccccccccccccccccc]
                       & p_0    & p_1    & \cdots & p_{d-t-1} \\
               \cmidrule(lr){1-15}
               C_{1,1}     & 1      & \;     & \;     & \; & *      & \cdots & *      & \;      & \;      & \;      & \;      & \;      & \; \\
               \vdots  & \vdots & \;     & \;     & \; & \vdots & \ddots & \vdots & -M_1'      & \;      & \;      & \;      & \;      & \; \\
               C_{1,m_1} & 1      & \;     & \;     & \; & *      & \cdots & *      & \;      & \;      & \;      & \;      & \;      & \; \\
               \cmidrule(lr){1-15}
               \vdots  & \vdots & \;     & \;     & \; & \;     & \;     & \;     & \ddots  & \ddots  & \ddots  & \;      & \;      & \; \\
               \vdots  & \vdots & \;     & \;     & \; & \;     & \;     & \;     & \ddots  & \ddots  & \ddots  & \;      & \;      & \; \\
               \vdots  & \vdots & \;     & \;     & \; & \;     & \;     & \;     & \ddots  & \ddots  & \ddots  & \;      & \;      & \; \\  
               \cmidrule(lr){1-15}
               C_{t,1}     & 1      & \;     & \;     & \; & \;     & \;     & \;     & \;      & \;      & \;      & *       & \cdots  & *  \\
               \vdots  & 1      & \;     & \;     & \; & \;     & \;     & \;     & \;      & \;      & \;      & \vdots       & \ddots  & \vdots & - M_t' \\
               C_{t,m_t} & 1      & \;     & \;     & \; & \;     & \;     & \;     & \;      & \;      & \;      & *       & \cdots  & *  \\
               \cmidrule(lr){1-15}
               D_1     & 1      & 1      & \;      & \;  & \;     & \;     & \;     & \;      & \;      & \;      & \;      & \;      & \; \\
               \vdots  & \vdots & \;      & \ddots & \;  & -D'     & \;     & \;     & \;      & \;      & \;      & \;      & \;      & \; \\
               D_{d-t-1} & 1      & \;      & \;      & 1  & \;     & \;     & \;     & \;      & \;      & \;      & \;      & \;      & \; \\
\end{bmatrix}$}
\end{minipage}
\label{mat:caseA}
\end{equation}
$M_i'$ means the corresponding block sub-matrices

$\mathbf{Case~B}$. Some entries are $0$ in the column $p_0$. Rows containing $0$-entries in the column $p_0$ come from only one branch. We may assume that the branch is the first branch and the rows are $C_{1,e}, \dots, C_{1,m_1}$ for $1 \leq e \leq m_1$. Moreover, there is at most one sub-matrix $M_i$ such that $[M_1, M_i, D]$ is type 2-1 in definition~\ref{def:Incidence matrix of type 1 and 2}. We may assume that $i = 2$. Then each sub-matrix $[C_1, C_j, D]$ is type 2-2 for $j = 3, \dots, t$.

If every sub-matrix $[M_1, M_i, D]$ is type 2-2 for all $i = 2, \cdots, t$, then the combinatorial incidence matrix is of the following form.
\begin{equation}
\resizebox{.91\textwidth}{!}
{
$\begin{bmatrix}[c|cccccccccccccccccccc]
                         & p_0    & p_1    &\cdots&p_{d-t-1}& q_{2,1}& \cdots & q_{2,g_2}& q_{3,1}& \cdots & q_{3,g_3} &q_{4,1} & \cdots & q_{t,g_t} & & & & \cdots\\
              \cmidrule(lr){1-18}
               C_{1,1}   & 1      & \;     & \;     & \;     & \;     & \;     & \;      & \;     & \;     & \;     \\
               \vdots    & \vdots & \;     & \;     & \;     & \;     & \;     & \;      & \;     & \;     & \;     \\
               C_{1,e-1} & 1      & \;     & \;     & \;     & \;     & \;     & \;      & \;     & \;     & \;     \\  
               \cmidrule(lr){1-18}
               C_{1,e}   & \;     & 1      & \cdots & 1      & 1      & \cdots & 1       & 1      & \cdots & 1      & 1      & \cdots & 1      & *      & \cdots & * \\
               \vdots    & \;     & \vdots & \ddots & \vdots & \vdots & \ddots & \vdots  & \vdots & \ddots & \vdots & \vdots & \ddots & \vdots & \vdots & \ddots & \vdots  \\
               C_{1,m_1}   & \;     & 1      & \cdots & 1      & 1      & \cdots & 1     & 1      & \cdots & 1      & 1      & \cdots & 1      & *      & \cdots & * \\
               \cmidrule(lr){1-18}

               C_{2,1}   & 1      & \;     & \;     & \;     & *      & \cdots & *       & \;     & \;     & \;     \\
               \vdots    & \vdots & \;     & \;     & \;     & \vdots & \ddots & \vdots  & \;     & \;     & \;     \\
               C_{2,m_2} & 1      & \;     & \;     & \;     & *      & \cdots & *       & \;     & \;     & \;     \\
               \cmidrule(lr){1-18}
               C_{3,1}   & 1      & \;     & \;     & \;     & \;     & \;     & \;      & *      & \cdots & *      \\
               \vdots    & \vdots & \;     & \;     & \;     & \;     & \;     & \;      & \vdots & \ddots & \vdots  \\
               C_{3,m_3} & 1      & \;     & \;     & \;     & \;     & \;     & \;      & *      & \cdots & *       \\
               \cmidrule(lr){1-18}
               \vdots    & \vdots & \;     & \;     & \;     & \;     & \;     & \;      &        &        &        & *      & \cdots & *      \\
               \vdots    & \vdots & \;     & \;     & \;     & \;     & \;     & \;      &        &        &        & \vdots & \ddots & \vdots \\
               \vdots    & \vdots & \;     & \;     & \;     & \;     & \;     & \;      &        &        &        & *      & \cdots & *      \\
               \cmidrule(lr){1-18}
               D_1       & 1      & 1      & 0      & 0      & \;     & \;     & \;      & \;     & \;     & \;     \\
               \vdots    & \vdots & 0      & \ddots & 0      & \;     & \;     & \;      & \;     & \;     & \;     \\
               D_{d-t-1} & 1      & 0      & 0      & 1      & \;     & \;     & \;      & \;     & \;     & \;     \\
\end{bmatrix}$}
\label{mat:caseB-1}
\end{equation}
If the sub-matrix $[M_1, M_2, D]$ is type 2-1, then the combinatorial incidence matrix is of the following form.
\begin{equation}
\resizebox{.91\textwidth}{!}
{
$\begin{bmatrix}[c|cccccccccccccccccccc]
                         & p_0    & p_1    &\cdots&p_{d-t-1}& q_{2,1}& \cdots & q_{2,g'}& q_{3,1}& \cdots & q_{3,g_3} & & & &q_{4,1} & \cdots & q_{t,g_t} & \cdots\\
              \cmidrule(lr){1-18}
               C_{1,1}   & 1      & \;     & \;     & \;     & \;     & \;     & \;      & \;     & \;     & \;     \\
               \vdots    & \vdots & \;     & \;     & \;     & \;     & \;     & \;      & \;     & \;     & \;     \\
               C_{1,e-1} & 1      & \;     & \;     & \;     & \;     & \;     & \;      & \;     & \;     & \;     \\  
               \cmidrule(lr){1-18}
               C_{1,e}   & \;     & 1      & \cdots & 1      & *      & \cdots & *       & 1      & \cdots & 1      & & & & 1      & \cdots & 1      \\
               \vdots    & \;     & \vdots & \ddots & \vdots & \vdots & \ddots & \vdots  & \vdots & \ddots & \vdots & & & & \vdots & \ddots & \vdots  \\
               C_{1,m_1}   & \;     & 1      & \cdots & 1      & *      & \cdots & *     & 1      & \cdots & 1      & & & & 1      & \cdots & 1       \\
               \cmidrule(lr){1-18}

               C_{2,1}   & 1      & \;     & \;     & \;     & *      & \cdots & *       & \;     & \;     & \;     \\
               \vdots    & \vdots & \;     & \;     & \;     & \vdots & \ddots & \vdots  & \;     & \;     & \;     \\
               C_{2,m_2} & 1      & \;     & \;     & \;     & *      & \cdots & *       & \;     & \;     & \;     \\
               \cmidrule(lr){1-18}
               C_{3,1}   & 1      & \;     & \;     & \;     & \;     & \;     & \;      & *      & \cdots & *      & *     & \cdots & * & \\
               \vdots    & \vdots & \;     & \;     & \;     & \;     & \;     & \;      & \vdots & \ddots & \vdots &\vdots & \ddots & * & \\
               C_{3,m_3} & 1      & \;     & \;     & \;     & \;     & \;     & \;      & *      & \cdots & *      & *     & \cdots & * & \\
               \cmidrule(lr){1-18}
               \vdots    & \vdots & \;     & \;     & \;     & \;     & \;     & \;      &        &        &        & & & & *      & \cdots & *      \\
               \vdots    & \vdots & \;     & \;     & \;     & \;     & \;     & \;      &        &        &        & & & & \vdots & \ddots & \vdots \\
               \vdots    & \vdots & \;     & \;     & \;     & \;     & \;     & \;      &        &        &        & & & & *      & \cdots & *      \\
               \cmidrule(lr){1-18}
               D_1       & 1      & 1      & 0      & 0      & \;     & \;     & \;      & \;     & \;     & \;     \\
               \vdots    & \vdots & 0      & \ddots & 0      & \;     & \;     & \;      & \;     & \;     & \;     \\
               D_{d-t-1} & 1      & 0      & 0      & 1      & \;     & \;     & \;      & \;     & \;     & \;     \\
\end{bmatrix}$
}
\label{mat:caseB-2}
\end{equation}
The columns $q_{2,1}, \dots, q_{2,g'}$ are columns of type 2-1 containing a stair-shaped sub-matrix that we mentioned in Lemma~\ref{lem:configurationoftype2-1andtype2-2}. The column $q_{i,g_i}(i = 3, \dots, t)$ is the column of type 2-2 that we mentioned in the same lemma.
\end{theorem}
\begin{lemma}
The block sub-matrix $M_i'$ of the matrix~\ref{mat:caseA} is an incidence matrix of a cyclic quotient surface singularity $[a_{i,n_i}, \dots, a_{i,1}]$.
\label{lem:onebranch}
\end{lemma}
\begin{proof}
Similar to the proof of lemma~\ref{lem:twobranches}.
\end{proof}
\begin{proof}[proof of theorem~\ref{thm:caseofcombinatorialincidencematrix}]
$\mathbf{Case~A}$. Assume that all entries of the column $p_0$ of the matrix $M$ are $1$. Then the intersection relations $C_{i,j}.C_{i',j'} = 1$ where $i \neq i'$ in Equation~\ref{eq:combinatorial equations of wighted homogeneous surface singularities} are satisfied at $p_0$. Therefore there are no additional intersection points between any two decorated curves $C_{i,j}$ and $C_{i',j'}$. In the aspect of combinatorial incidence matrices, there are no columns $p$ such that $M(C_{i,j}, p) = M(C_{i',j'},p) = 1$ for $i \neq i'$ except the column $p_0$. Therefore by proper column exchanging, we can make the block sub-matrices $M_1', \dots, M_t'$ as in the matrix~\ref{mat:caseA}.  

$\mathbf{Case~B}$. Assume that some entries of the column $p_0$ are $0$. Let $C_{i,j}$ be such rows that $M(C_{i,j},p_0) = 0$. It is equivalent to that every decorated curve passes through $p_0$ except curves $C_{i,j}$. To satisfy the intersection relation $C_{i,j}.D_k = 1$ for $k = 1, \dots, d-t-1$, the curves $C_{i,j}$ must pass through $p_1, \dots, p_{d-t-1}$. It induces the columns $p_1, \dots, p_{d-t-1}$ of the matrices~\ref{mat:caseB-1} and ~\ref{mat:caseB-2}.

To show that the rows containing $0$-entries of the column $p_0$ come from only one branch, assume that $M(C_{1,1},p_0) = M(C_{2,1},p_0) = 0$. This is equivalent to that $C_{1,1}$ and $C_{2,1}$ does not pass through $p_0$. Since $d \geq t + 3$, there are at least two points $p_1$ and $p_2$. To satisfy the intersection relation $C_{1,1}.D_k = 1$ and $C_{2,1}.D_k = 1$ for $k = 1, 2$, the curves $C_{1.1}$ and $C_{2,1}$ pass through $p_1$ and $p_2$. Then $C_{1,1}.C_{2,1} \geq 2$. But $C_{1,1}.C_{2,1} = 1$ by the intersection relation. Therefore decorated curves that do not pass through $p_0$ come from only one branch. We may assume that this branch is the first branch.

To show the rest of the lemma, we suppose that two sub-matrices $[M_1, M_2, D]$ and $[M_1, M_3, D]$ are the type 2-1. By Lemma~\ref{lem:twobranches}, we consider $[M_1,M_2,D]$ and $[M_1,M_3,D]$ are incidence matrices of cyclic quotient surface singularities. By Lemma~\ref{lem:configurationoftype2-1andtype2-2}, the stair-shaped sub-matrix exists only one in $M_1$. Therefore some decorated curves in $M_2$ and $M_3$ intersect at some $q_j$. But the intersection relations $C_{2,j}.C_{3,j'} = 1$ are already satisfied at the column $p_0$. This contradiction means that there is at most one sub-matrix of type 2-1.
\end{proof}
\begin{theorem}
The map $\phi_{PI} : \mathscr{P}(X) \rightarrow C\mathscr{I}(X)$ from the set of $P$-resolutions of $X$ to the set of combinatorial incidence matrices of $X$ is surjective.
\label{thm:main}
\end{theorem}
\begin{proof} 
We construct a P-resolution of $X$ from a given combinatorial incidence matrix $M$ of $X$ for each case in theorem~\ref{thm:caseofcombinatorialincidencematrix}. We then show that the constructed P-resolution induces the given combinatorial incidence matrix by using the MMP algorithm.

$\mathbf{Case~A}$ All entries of the column $p_0$ are $1$.(Matrix~\ref{mat:caseA}).\\
By eliminating the column $p_0$ of the matrix $M$, we obtain block sub-matrices $M_1', \cdots, M_t'$ and $D' = I_{(d-t-1) \times (d-t-1)}$.

By lemma~\ref{lem:onebranch}, we consider the sub-matrices $M_i$ as incidence matrices of cyclic quotient surface singularities $[a_{i,n_i}, \dots, a_{i,1}]$. We can find the P-resolution for each cyclic quotient surface singularity that induces the incidence matrix $M_i'$ respectively. This means that we know where the T-singularities are located, that is, which exceptional curves in the minimal resolution are contracted to be a T-singularity. Since each dual resolution graph of the cyclic quotient surface singularity is an branch of the dual graph of $X$, we can locate T-singularities on each branch. Therefore we obtain a P-resolution of $X$.

We apply the MMP algorithm to the P-resolution that we construct now(refer Figure~\ref{mmp1}). Note that the MMP algorithm on each P-resolution of $[a_{i,n_i}, \dots, a_{i,1}]$ is the same with the MMP algorithm on each branch of the P-resolution of $X$. Therefore the MMP algorithm induces the same matrices $M_i'$. The $(-1)$-curves that appears on each branch do not connect two decorated curves in different branches. A $(-1)$-curve connecting them is only the central curve. The column $p_0$ comes from this $(-1)$-curve. The sub-matrix $D'$ comes from the $(-1)$-curves attached to the decorated curves $D_i$ respectively.
\begin{figure}
\begin{minipage}{\textwidth}
\begin{center}
\begin{tikzpicture}
  [inner sep=1mm,
R/.style={circle,draw=black!255,fill=white!20,thick},
T/.style={circle,draw=red!255,fill=white!255,thick}]
    \node[R] (e) [label={[yshift=-.5cm]:$-d$}]{};
    \node    (d) [below left=.7cm of e]{$\udots$};
    \node[R] (c) [below left=.7cm of d, label=right:$-a_{1,j}$]{};
    \node    (b) [below left=.7cm of c]{$\udots$};
    \node[R] (a) [below left=.7cm of b, label=below:$-a_{1,n_1}$]{};
    \node    (f) [below right=.7cm of e]{$\ddots$};
    \node[R] (g) [below right=.7cm of f, label=left:$-a_{t,j'}$]{};
    \node    (h) [below right=.7cm of g]{$\ddots$};
    \node[R] (i) [below right=.7cm of h, label=below:$-a_{t,n_t}$]{};

    \node[T] (j) [above left=.5cm and 0.8cm of a]{};
    \node    (k) [right=.1cm of j]{$\cdots$};
    \node[T] (l) [above= 0.5cm of a]{};
    \node[R] (m) [above=.4cm of j, label=above:$C_{1,m_1}$]{};
    \node    (k) [right=.1cm of m]{$\cdots$};
    \node[R] (n) [above=.4cm of l]{};

    \node[T] (o) [above left=.5cm and 0.8cm of c]{};
    \node    (p) [right=.1cm of o]{$\cdots$};
    \node[T] (q) [above=.5cm of c]{};
    \node[R] (r) [above=.4cm of o]{};
    \node    (k) [right=.1cm of r]{$\cdots$};
    \node[R] (s) [above=.4cm of q]{};

    \node[T] (t) [above=.5cm of g]{};
    \node    (u) [right=.1cm of t]{$\cdots$};
    \node[T] (v) [above right=.5cm and 0.8cm of g]{};
    \node[R] (w) [above=.4cm of t]{};
    \node    (k) [right=.1cm of w]{$\cdots$};
    \node[R] (x) [above=.4cm of v]{};

    \node[T] (y) [above=.5cm of i]{};
    \node    (z) [right=.1cm of y]{$\cdots$};
    \node[T] (aa) [above right=.5cm and 0.8cm of i]{};
    \node[R] (ab) [above=.4cm of y]{};
    \node    (k) [right=.1cm of ab]{$\cdots$};
    \node[R] (ac) [above=.4cm of aa, label=above:$C_{t,m_2}$]{};

    \node[T] (ad) [above left=.5cm of e]{};
    \node    (ae) [right=.1cm of ad]{$\cdots$};
    \node[T] (af) [above right=.5cm of e]{};
    \node[R] (ag) [above=.4cm of ad, label=above:$D_{1}$]{};
    \node    (k) [right=.1cm of ag]{$\cdots$};
    \node[R] (ah) [above=.4cm of af, label=above:$D_{d-t-1}$]{};

    \node    (CE) [below=1.5cm of e] {$\cdots$};

    \draw[-] (a)--(b);
    \draw[-] (b)--(c);
    \draw[-] (c)--(d);
    \draw[-] (d)--(e);
    \draw[-] (e)--(f);
    \draw[-] (f)--(g);
    \draw[-] (g)--(h);
    \draw[-] (h)--(i);
    
    \draw[-] (a)--(j);
    \draw[-] (a)--(l);
    \draw[-] (j)--(m);
    \draw[-] (l)--(n);

    \draw[-] (c)--(o);
    \draw[-] (c)--(q);
    \draw[-] (r)--(o);
    \draw[-] (s)--(q);

    \draw[-] (t)--(g);
    \draw[-] (v)--(g);
    \draw[-] (w)--(t);
    \draw[-] (x)--(v);

    \draw[-] (y)--(i);
    \draw[-] (aa)--(i);
    \draw[-] (ab)--(y);
    \draw[-] (ac)--(aa);

    \draw[-] (ad)--(e);
    \draw[-] (af)--(e);
    \draw[-] (ag)--(ad);
    \draw[-] (ah)--(af);

\draw [decorate,decoration={brace,amplitude=5pt,raise=4ex}]
  (-4,-1.7) -- (-3,-1.7) node[midway,yshift=3em]{$a_{1,1}-1$};

\draw [decorate,decoration={brace,amplitude=5pt,raise=4ex}]
  (-2.5,-0.6) -- (-1.5,-0.6) node[midway,yshift=3em]{$a_{1,j}-2$};

\draw [decorate,decoration={brace,amplitude=5pt,raise=4ex}]
  (1.5,.-0.6) -- (2.5,.-0.6) node[midway,yshift=3em]{$a_{t,j'}-2$};

\draw [decorate,decoration={brace,amplitude=5pt,raise=4ex}]
  (3,-1.7) -- (4,-1.7) node[midway,yshift=3em]{$a_{t,n_t}-1$};

\draw [decorate,decoration={brace,mirror,amplitude=5pt,raise=4ex}]
  (-2.7,-2.7) -- (-0.3,-0.3) node[midway,xshift=4em,yshift=-2em]{$T$-singularities};

\draw [decorate,decoration={brace,amplitude=5pt,raise=4ex}]
  (2.7,-2.7) -- (0.3,-0.3) node[midway,xshift=4em,yshift=-2em]{};
\end{tikzpicture}
\end{center}
\end{minipage}
\\
\begin{minipage}{\textwidth}
\begin{center}
$\downarrow$
\end{center}
\end{minipage}
\\
\begin{center}
\begin{minipage}{.3\textwidth}
\begin{center}
\begin{tikzpicture}
[inner sep=1mm,
R/.style={circle,draw=black!255,fill=white!20,thick},
T/.style={circle,draw=red!255,fill=white!255,thick}]
    \node[R] (a) [label={[yshift=.3cm]:$-d+t$}]{};
    \node[R] (b) [above left=.5cm and .7cm of a, label=left:$C_{1,1}$]{};
    \node    (c) [above left=.1cm and .9cm of a]{$\vdots$};
    \node[R] (d) [left=.7cm of a, label=left:$C_{1,m_1}$]{};
    \node[T] (e) [below left=.5cm and .3cm of a]{};
    \node    (f) [below=.6cm of a]{$\cdots$};
    \node[T] (g) [below right=.5cm and .3cm of a]{};
    \node[R] (h) [below =.5cm of e, label=below:$D_1$]{};
    \node[R] (i) [below =.5cm of g, label=below:$D_{d-t-1}$]{};
    \node[R] (j) [above right=.5cm and .7cm of a, label=right:$C_{2,1}$]{};
    \node    (k) [above right=.1cm and .9cm of a]{$\vdots$};
    \node[R] (l) [right=.7cm of a, label=right:$C_{2,m_2}$]{};

    \draw[-] (a) -- (b);
    \draw[-] (a) -- (d);
    \draw[-] (a) -- (e);
    \draw[-] (a) -- (g);
    \draw[-] (e) -- (h);
    \draw[-] (i) -- (g);
    \draw[-] (a) -- (j);
    \draw[-] (a) -- (l);
\end{tikzpicture}
\end{center}
\end{minipage}
\begin{minipage}{.1\textwidth}
\begin{center}
$\to$
\end{center}
\end{minipage}
\begin{minipage}{.3\textwidth}
\begin{center}
$\begin{bmatrix}[ccc]
M_1' &        &      \\
     & \ddots &      \\
     &        & M_t'
\end{bmatrix}$
\end{center}
\end{minipage}
\end{center}

\begin{minipage}{\textwidth}
\begin{center}
$\downarrow$
\end{center}
\end{minipage}
\\
\begin{center}
\begin{minipage}{.3\textwidth}
\begin{center}
\begin{tikzpicture}
[inner sep=1mm,
R/.style={circle,draw=black!255,fill=white!20,thick},
T/.style={circle,draw=red!255,fill=white!255,thick}]
    \node[T] (a) [label={[yshift=.3cm]:$-1$}]{};
    \node[R] (b) [above left=.5cm and .7cm of a, label=left:$C_{1,1}$]{};
    \node    (c) [above left=.1cm and .9cm of a]{$\vdots$};
    \node[R] (d) [left=.7cm of a, label=left:$C_{1,m_1}$]{};
    \node[R] (e) [below left=.5cm and .3cm of a, label=below:$D_1$]{};
    \node    (f) [below=.3cm of a]{$\cdots$};
    \node[R] (g) [below right=.5cm and .3cm of a, label=below:$D_{d-t-1}$]{};
    
    \node[R] (j) [above right=.5cm and .7cm of a, label=right:$C_{2,1}$]{};
    \node    (k) [above right=.1cm and .9cm of a]{$\vdots$};
    \node[R] (l) [right=.7cm of a, label=right:$C_{2,m_2}$]{};

    \draw[-] (a) -- (b);
    \draw[-] (a) -- (d);
    \draw[-] (a) -- (e);
    \draw[-] (a) -- (g);

    \draw[-] (a) -- (j);
    \draw[-] (a) -- (l);
\end{tikzpicture}
\end{center}
\end{minipage}
\begin{minipage}{.1\textwidth}
\begin{center}
$\to$
\end{center}
\end{minipage}
\begin{minipage}{.3\textwidth}
\begin{center}
$\begin{bmatrix}[cccc]
     &M_1' &        &      \\
     &     & \ddots &      \\
     &     &        & M_t' \\
D'   &     &        &      \\
\end{bmatrix}$
\end{center}
\end{minipage}
\end{center}

\begin{minipage}{\textwidth}
\begin{center}
$\downarrow$
\end{center}
\end{minipage}
\\
\begin{center}
\begin{minipage}{.3\textwidth}
\begin{center}
$C=\bigcup C_i$
\end{center}
\end{minipage}
\begin{minipage}{.1\textwidth}
\begin{center}
$\to$
\end{center}
\end{minipage}
\begin{minipage}{.3\textwidth}
\begin{center}
$\begin{bmatrix}[ccccc]
\mathbf{1} &      &M_1' &        &      \\
\mathbf{1} &      &     & \ddots &      \\
\mathbf{1} &      &     &        & M_t' \\
\mathbf{1} & D'   &     &        &      \\
\end{bmatrix}$
\end{center}
\end{minipage}
\end{center}
\caption{MMP on $\mathbf{Case A}$}
\label{mmp1}
\end{figure}

$\mathbf{Case~B}$ $M(C_{1,e},p_0)=\cdots=M(C_{1,m_1},p_0)=0$.\\
$\mathbf{Case~B-1})$ We assume that the sub-matrix $[M_1, M_i, D]$ is type 2-2 for all $i = 2, \cdots, t$(Matrix~\ref{mat:caseB-1}).
We define $t-1+s = \sum\limits_{i=2}^t g_i$. That is, $s$ is the sum of $g_i$ such that $g_i \geq 2$. We consecutively blow up the intersection of the central curve and $A_{i,1}$ in the minimal resolution to make $(g_i-2)$ $(-2)$-curves for $i = 2, \cdots, t$.
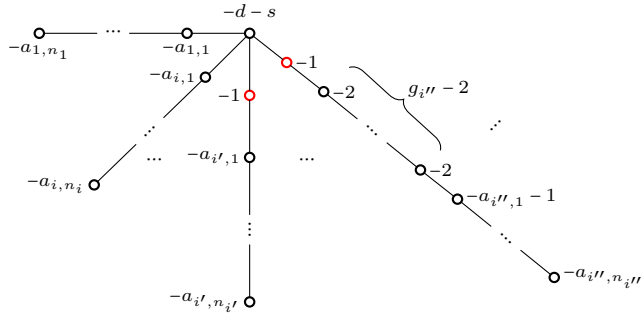
\begin{center}
\begin{tikzpicture}
  [inner sep=1mm,
R/.style={circle,draw=black!255,fill=white!20,thick},
S/.style={rectangle,draw=black!255,fill=white!20,thick},
T/.style={circle,draw=red!255,fill=white!255,thick}]
    \node[R] (c) [label={[yshift=.1cm]:$-d-s$}]{};

    \node[R] (a11) [left=.7cm of c, label=below:$-a_{1,1}$]{};
    \node    (a12) [left=.7cm of a11]{$\cdots$};
    \node[R] (a13) [left=.7cm of a12, label=below:$-a_{1,n_1}$]{};

    \node[T] (at1) [below right=.3cm and .4cm of c, label=right:$-1$]{};
    \node[R] (at2) [below right=.3cm and .4cm of at1, label=right:$-2$]{};
    \node    (at3) [below right=.3cm and .4cm of at2]{$\ddots$};
    \node[R] (at4) [below right=.3cm and .4cm of at3, label=right:$-2$]{};
    \node[R] (at5) [below right=.3cm and .4cm of at4, label=right:$-a_{i'',1}-1$]{};
    \node    (at6) [below right=.3cm and .4cm of at5]{$\ddots$};
    \node[R] (at7) [below right=.3cm and .4cm of at6, label=right:$-a_{i'',n_{i''}}$]{};

    \node[R] (a21) [below left=.7cm of c, label=left:$-a_{i,1}$]{};    
    \node    (a22) [below left=.7cm of a21]{$\udots$}; 
    \node[R] (a23) [below left=.7cm of a22, label=left:$-a_{i,n_i}$]{}; 

    \node[T] (a31) [below=.7cm of c, label=left:$-1$]{};
    \node[R] (a32) [below=.7cm of a31, label=left:$-a_{i',1}$]{};    
    \node    (a33) [below=.7cm of a32]{$\vdots$}; 
    \node[R] (a34) [below=.7cm of a33, label=left:$-a_{i',n_{i'}}$]{};

    \node    (e1) [below left=1.5cm and 1cm of c]{$\cdots$};
    \node    (e2) [below right=1.5cm and .5cm of c]{$\cdots$};
    \node    (e3) [below right=1cm and 3cm of c]{$\udots$};

    \draw[-] (c)--(a11);
    \draw[-] (a11)--(a12);
    \draw[-] (a12)--(a13);

    \draw[-] (c)--(at1);
    \draw[-] (at1)--(at2);
    \draw[-] (at2)--(at3);
    \draw[-] (at3)--(at4);
    \draw[-] (at4)--(at5);
    \draw[-] (at5)--(at6);
    \draw[-] (at6)--(at7);

    \draw[-] (c)--(a21);
    \draw[-] (a21)--(a22);
    \draw[-] (a22)--(a23);

    \draw[-] (c)--(a31);
    \draw[-] (a31)--(a32);
    \draw[-] (a32)--(a33);
    \draw[-] (a33)--(a34);

\draw [decorate,decoration={brace,amplitude=5pt,raise=4ex}]
  (.9,-.9) -- (2.1,-2) node[midway,xshift=2.8em, yshift=2em]{$g_{i''} - 2$};
\end{tikzpicture}
\end{center}
We will locate T-singularities on $[a_{1,n_1},\dots,a_{1,1},d+s]$ and $[a_{i,n_i},\dots,a_{i,1}+1, 2, \dots, 2]$ for $i = 2, \cdots, t$. 

(1) Consider a sub-matrix $[M_1, D]$(Figure~\ref{fig.c1d}). 
\begin{figure}
\centering
$\begin{bmatrix}[c|ccccccccccccccccccc]
                       & p_0 & p_1 & \cdots & p_{d-t} & q_{2,1} & \cdots & q_{t,g_t} &   &   &   &   \\
              \cmidrule(lr){1-12}
               C_{1,1} & 1      & \;     & \;     & \;     & \;     & \;     & \;      & \;     & \;     & \;     \\
               \vdots    & \vdots & \;     & \;     & \;     & \;     & \;     & \;      & \;     & \;     & \;     \\
               C_{1,e-1} & 1      & \;     & \;     & \;     & \;     & \;     & \;      & \;     & \;     & \;     \\  
               \cmidrule(lr){1-12}
               C_{1,e}   & \;     & 1      & \cdots & 1      & 1      & \cdots & 1       & *      & \cdots & *      \\
               \vdots    & \;     & \vdots & \ddots & \vdots & \vdots & \ddots & \vdots  & \vdots & \ddots & \vdots \\
               C_{1,m_1}   & \;     & 1      & \cdots & 1      & 1      & \cdots & 1     & *      & \cdots & *      \\
               \cmidrule(lr){1-12}
               D_1       & 1      & 1      &        &       & \;     & \;     & \;      & \;      & \;      & \;      & \;         \\
               \vdots    & \vdots &        & \ddots &       & \;     & \;     & \;      & \;      & \;      & \;      & \;         \\
               D_{d-t-1} & 1      &        &        & 1      & \;     & \;     & \;      & \;      & \;      & \;      & \;         \\
\end{bmatrix}$
\caption{$[M_1, D]$}
\label{fig.c1d}
\end{figure}
We add decorated curves $E_1, \dots, E_{t-1+s}$ such that combinatorial equations are 
\begin{equation*}
\begin{split}
l(E_i) = 2 \\
E_i.E_{i'} = 1 \\
C_{1,j}.E_i = 1 
\end{split}
\end{equation*}
for $i, i' = 1, \dots, t-1+s$(Figure~\ref{fig:c1de})
\begin{figure}
\centering
$\begin{bmatrix}[c|ccccccccccccccccccc]
                         & p_0    & p_1    & \cdots & p_{d-t} & q_{2,1} & \cdots & q_{t,g_t} &        &        &         &      \\
              \cmidrule(lr){1-12}
               C_{1,1}   & 1      & \;     & \;     & \;      & \;      & \;     & \;        & \;     & \;     & \;      &      \\
               \vdots    & \vdots & \;     & \;     & \;      & \;      & \;     & \;        & \;     & \;     & \;      &      \\
               C_{1,e-1} & 1      & \;     & \;     & \;      & \;      & \;     & \;        & \;     & \;     & \;      &      \\  
                \cmidrule(lr){1-12}
               C_{1,e}   & \;     & 1      & \cdots & 1       & 1       & \cdots & 1         & *      & \cdots & *       &      \\
               \vdots    & \;     & \vdots & \ddots & \vdots  & \vdots  & \ddots & \vdots    & \vdots & \ddots & \vdots  &      \\
               C_{1,m_1} & \;     & 1      & \cdots & 1       & 1       & \cdots & 1         & *      & \cdots & *       &      \\
               \cmidrule(lr){1-12}
               D_1       & 1      & 1      &        &         & \;      & \;     & \;        & \;     & \;     & \;      & \;   \\
               \vdots    & \vdots &        & \ddots &         & \;      & \;     & \;        & \;     & \;     & \;      & \;   \\
               D_{d-t-1} & 1      &        &        & 1       & \;      & \;     & \;        & \;     & \;     & \;      & \;   \\
               \cmidrule(lr){1-12}
               E_1       & 1      & \;     & \;     & \;      & 1       & \cdots & 0         &        &        &         &      \\
               \vdots    & \vdots & \;     & \;     & \;      & \vdots  & \ddots & \vdots    &        &        &         &      \\
               E_{t-1+s} & 1      & \;     & \;     &         & 0       & \cdots & 1         &        &        &         &          
\end{bmatrix}$
\caption{$[M_1, D, E]$}
\label{fig:c1de}
\end{figure}

The combinatorial equations of $E_i$ are the same with those of $D_k$. Therefore we can consider this matrix as an incidence matrix of a cyclic quotient surface singularity $[a_{1,n_1},\dots,a_{1,1},d+s]$(Ref Lemma~\ref{lem:twobranches}). We can find the $P$-resolution of the cyclic quotient surface singularity that induces the matrix $[C_1, D, E]$. Since every decorated curve $D_k$ and $E_l$ has no free point, we know that the $(-d-s)$-curve is an exceptional curve of a Wahl singularity by Lemma~\ref{lem:No free point ensure that the last curve is an exceptional of Wahl singularity}. Moreover, the decorated curves $C_{1,e}, \cdots, C_{1,m_1}$ degenerate to the $(-d-s)$-curve after flips because of the columns $p_1, \cdots, p_{d-k}$.

Let $M_1'$ be the matrix obtained from $M_1$ by deleting the columns $q_{2,1}, \dots, q_{t,g_t}$. If we progress the MMP algorithm until the last $-(d+s)$-curve becomes $-(t+s)$-curve, then we obtain the matrix $[M_1',D]$(See Figure~\ref{fig:partialmmpofc1c2de}). In fact, the additional curves $E_1,\dots,E_{t-1+s}$ correspond to the $(t-1)$ branches.
\begin{figure}
\begin{center}
\begin{minipage}{.75\textwidth}
\begin{center}
\begin{tikzpicture}
  [inner sep=1mm,
R/.style={circle,draw=black!255,fill=white!20,thick},
S/.style={rectangle,draw=black!255,fill=white!20,thick},
T/.style={circle,draw=red!255,fill=white!255,thick}]
    \node[R] (a) [label=below:$-a_{1,n_1}$]{};
    \node    (b) [right=.7cm of a]{$\cdots$};
    \node[R] (c) [right=.7cm of b, label=below:$-a_{1,j}$]{};
    \node    (d) [right=.7cm of c]{$\cdots$};
    \node[S] (e) [right=.7cm of d, label=right:$-d-s$]{};

    \node[T] (j) [above left=.5cm of a]{};
    \node    (k) [above=.2cm of a]{$\cdots$};
    \node[T] (l) [above right=.5cm of a]{};
    \node[R] (m) [above=.4cm of j, label=above:$C_{1,m_1}$]{};
    \node    (k) [above=.75cm of a]{$\cdots$};
    \node[R] (n) [above=.4cm of l]{};

    \node[T] (o) [above left=.5cm of c]{};
    \node    (p) [above=.2cm of c]{$\cdots$};
    \node[T] (q) [above right=.5cm of c]{};
    \node[R] (r) [above=.4cm of o]{};
    \node    (k) [above=.75cm of c]{$\cdots$};
    \node[R] (s) [above=.4cm of q]{};

    \node[T] (ad) [above left=.5cm of e]{};
    \node    (ae) [above=.2cm of e]{$\cdots$};
    \node[T] (af) [above right=.5cm of e]{};
    \node[R] (ag) [above=.4cm of ad, label=above:$D_{1}$]{};
    \node    (k)  [above=.75cm of e]{$\cdots$};
    \node[R] (ah) [above=.4cm of af, label=above:$D_{d-t-1}$]{};

    \node[T] (ai) [below left=.5cm of e]{};
    \node    (aj) [below=.3cm of e]{$\cdots$};
    \node[T] (ak) [below right=.5cm of e]{};
    \node[R] (al) [below=.5cm of ai, label=below:$E_{1}$]{};
    \node    (am) [below=.4cm of aj]{$\cdots$};
    \node[R] (an) [below=.5cm of ak, label=below:$E_{t-1+s}$]{};

    \draw[-] (a)--(b);
    \draw[-] (b)--(c);
    \draw[-] (c)--(d);
    \draw[-] (d)--(e);
    
    \draw[-] (a)--(j);
    \draw[-] (a)--(l);
    \draw[-] (j)--(m);
    \draw[-] (l)--(n);

    \draw[-] (c)--(o);
    \draw[-] (c)--(q);
    \draw[-] (r)--(o);
    \draw[-] (s)--(q);

    \draw[-] (ad)--(e);
    \draw[-] (af)--(e);
    \draw[-] (ag)--(ad);
    \draw[-] (ah)--(af);

    \draw[-] (e)--(ai);
    \draw[-] (e)--(ak);
    \draw[-] (ai)--(al);
    \draw[-] (ak)--(an);

\draw [decorate,decoration={brace,amplitude=5pt,raise=4ex}]
  (-0.5,.8) -- (0.5,.8) node[midway,yshift=3em]{$a_{1,1}-1$};

\draw [decorate,decoration={brace,amplitude=5pt,raise=4ex}]
  (1.4,.8) -- (2.4,.8) node[midway,yshift=3em]{$a_{1,j}-2$};

\end{tikzpicture}
\end{center}
\end{minipage}
\end{center}
\begin{center}
\begin{minipage}{.3\textwidth}
\begin{center}
$\downarrow$
\end{center}
\end{minipage}
\end{center}

\begin{center}
\begin{minipage}{.4\textwidth}
\begin{tikzpicture}
[inner sep=1mm,
R/.style={circle,draw=black!255,fill=white!20,thick},
T/.style={circle,draw=red!255,fill=white!255,thick}]
    \node[R] (a) [label=right:$-t-s$]{};
    \node[R] (b) [above left=.1cm and 1cm of a, label=left:$C_{1,e}$]{};
    \node    (c) [below left=.05cm and .9cm of a]{$\vdots$};
    \node[R] (d) [below left=.5cm and 1cm of a, label=left:$C_{1,m_1}$]{};
    \node[T] (e) [below left=.7cm and .5cm of a]{};
    \node    (f) [below=.6cm of a]{$\cdots$};
    \node[T] (g) [below right=.7cm and .5cm of a]{};
    \node[R] (h) [below =.5cm of e, label=below:$E_1$]{};
    \node[R] (i) [below =.5cm of g, label=below:$E_{t-1+s}$]{};

    \node[R] (m) [above left=.7cm and .7cm of a, label=above:$D_{1}$]{};
    \node[R] (n) [above right=.7cm and .7cm of a, label=above:$D_{d-t-1}$]{};
    \node    (o) [above=.6cm of a]{$\cdots$};

    \node[R] (p) [above right=.1cm and 1cm of a, label=right:$C_{1,1}$]{};
    \node    (q) [below right=.05cm and .9cm of a]{$\vdots$};
    \node[R] (r) [below right=.5cm and 1cm of a, label=right:$C_{1,e-1}$]{};

    \draw[-] (a) -- (b);
    \draw[-] (a) -- (d);
    \draw[-] (a) -- (e);
    \draw[-] (a) -- (g);
    \draw[-] (e) -- (h);
    \draw[-] (i) -- (g);
    \draw[-] (a) -- (j);
    \draw[-] (a) -- (l);
    \draw[-] (a) -- (m);
    \draw[-] (a) -- (n);
    \draw[-] (a) -- (p);
    \draw[-] (a) -- (r);

    \draw[edge,line width=10pt] (a) -- (b);
    \draw[edge,line width=10pt] (a) -- (d);
\end{tikzpicture}
\end{minipage}
\begin{minipage}{.1\textwidth}
$\to$
\end{minipage}
\begin{minipage}{.1\textwidth}
$[M_1', D]$
\end{minipage}
\end{center}
\caption{partial MMP on $[M_1,M_2,D,E]$}
\label{fig:partialmmpofc1c2de}
\end{figure}
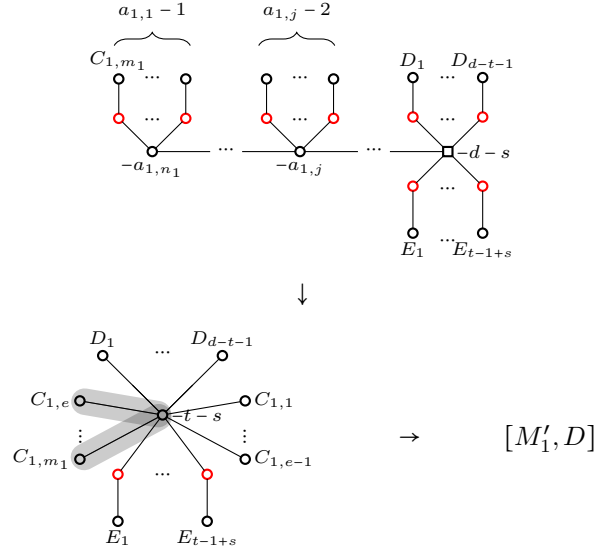

(2) Consider the sub-matrix $M_i$ for $i = 2, \cdots, t$. 
\begin{center}
$M_i = \begin{bmatrix}[c|ccccccccccccccccccc]
                       & p_0 & q_{i,1} & \cdots & q_{i,g_i} \\
              \cmidrule(lr){1-8}
               C_{i,1}   & 1      & *      & \cdots & *      & *      & \cdots & *      \\
               \vdots    & \vdots & \vdots & \ddots & \vdots & \vdots & \ddots & \vdots \\
               C_{i,m_i} & 1      & *      & \cdots & *      & *      & \cdots & *      \\
\end{bmatrix}$
\end{center}
We consider three cases whether $g_i = 1$, $g_i = 2$ or $g_i \geq 3$.

Assume $g_i = 1$. Then 
\begin{center}
$M_i = \begin{bmatrix}[c|ccccccccccccccccccc]
                       & p_0 & q_{i,1} \\
              \cmidrule(lr){1-6}
               C_{i,1}   & 1      & 1      & *      & \cdots & *      \\
               \vdots    & \vdots & \vdots & \vdots & \ddots & \vdots \\
               C_{i,m_i} & 1      & 1      & *      & \cdots & *      \\
\end{bmatrix}$
\end{center}
We delete the column $p_0$ and we denote this matrix as $M_i'$.
\begin{center}
\begin{minipage}{.4\textwidth}
\resizebox{\textwidth}{!}
{
$M_i' = \begin{bmatrix}[c|ccccccccccccccccccc]
                       & q_{i,1} \\
              \cmidrule(lr){1-5}
               C_{i,1}   & 1      & *      & \cdots & *      \\
               \vdots    & \vdots & \vdots & \ddots & \vdots \\
               C_{i,m_i} & 1      & *      & \cdots & *      \\
\end{bmatrix}$
}
\end{minipage}
\end{center}
This matrix satisfies the combinatorial constraints of the quotient surface singularity $[a_{i,n_i}, \dots, a_{i,1}]$. Therefore the matrix is an incidence matrix of the cyclic quotient surface singularity. We can find the corresponding P-resolution. The exceptional curve $A_{i,1}$ is not an exceptional curve of a Wahl singularity because of the column $q_{i,1}$(Lemma~\ref{lem:The 1-column ensures that the last curve is not an exceptional of Wahl singularity}).

Assume that $g_i = 2$.
\begin{center}
$M_i = \begin{bmatrix}[c|ccccccccccccccccccc]
                       & p_0 & q_{i,1} & q_{i,2} \\
              \hline
               C_{i,1}   & 1      & *      & *      & *      & \cdots & *      \\
               \vdots    & \vdots & \vdots & \vdots & \vdots & \ddots & \vdots \\
               C_{i,m_i} & 1      & *      & *      & *      & \cdots & *      \\
\end{bmatrix}$
\end{center}
Define a matrix $M_i'$ by deleting the column $p_0$ from $M_i$ and $M_i''$ by deleting the columns $q_{i,1}$ and $q_{i,2}$ from $M_i'$.
\begin{center}
$M_i' = \begin{bmatrix}[c|ccccccccccccccccccc]
                       & q_{i,1} & q_{i,2} \\
              \hline
               C_{i,1}   & *      & *      & *      & \cdots & *      \\
               \vdots    & \vdots & \vdots & \vdots & \ddots & \vdots \\
               C_{i,m_i} & *      & *      & *      & \cdots & *      \\
\end{bmatrix}$
,
$M_i'' = \begin{bmatrix}[c|ccccccccccccccccccc]
                         &        &  \\
              \hline
               C_{i,1}   & *      & \cdots & *      \\
               \vdots    & \vdots & \ddots & \vdots \\
               C_{i,m_i} & *      & \cdots & *      \\
\end{bmatrix}$
\end{center}
We add a row $F$ to $M_i'$.
\begin{center}
\begin{minipage}{.35\textwidth}
\resizebox{\textwidth}{!}
{
$\begin{bmatrix}[c|ccccccccccccccccccc]
                       & q_{i,1} & q_{i,2} \\
              \hline
               C_{i,1}   & *      & *      & *      & \cdots & *      \\
               \vdots    & \vdots & \vdots & \vdots & \ddots & \vdots \\
               C_{i,m_i} & *      & *      & *      & \cdots & *      \\
               F         & 1      & 1      &        &        &
\end{bmatrix}$
}
\end{minipage}
\end{center}
Note that the intersection relation between $F$ and $C_{i,1}, \cdots, C_{i,m_i}$ is the same with the relation between $C_{1,e}$ and $C_{i,1}, \cdots, C_{i,m_i}$. That is, $F.M_{i,j} = 1$ for $j = 1, \cdots, m_i$ and $l(F) = 2$. Therefore this matrix satisfies the combinatorial equations of the following sandwiched structure.
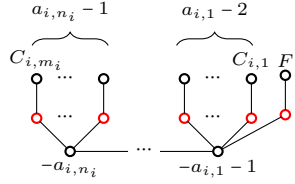
\begin{center}
\begin{tikzpicture}
  [inner sep=1mm,
R/.style={circle,draw=black!255,fill=white!20,thick},
T/.style={circle,draw=red!255,fill=white!255,thick}]
    \node[R] (a) [label=below:$-a_{i,n_i}$]{};
    \node    (b) [right=.7cm of a]{$\cdots$};
    \node[R] (c) [right=.7cm of b, label=below:$-a_{i,1}-1$]{};

    \node[T] (j) [above left=.5cm of a]{};
    \node    (k) [right=.1cm of j]{$\cdots$};
    \node[T] (l) [above right=.5cm of a]{};
    \node[R] (m) [above=.4cm of j, label=above:$C_{i,m_i}$]{};
    \node    (k) [right=.1cm of m]{$\cdots$};
    \node[R] (n) [above=.4cm of l]{};

    \node[T] (o) [above left=.5cm of c]{};
    \node    (p) [right=.1cm of o]{$\cdots$};
    \node[T] (q) [above right=.5cm of c]{};
    \node[R] (r) [above=.4cm of o]{};
    \node    (k) [right=.1cm of r]{$\cdots$};
    \node[R] (s) [above=.4cm of q, label=above:$C_{i,1}$]{};

    \node[T] (af) [above right=.4cm and .8cm of c]{};
    \node[R] (ag) [above=.4cm of af, label=above:$F$]{};

    \draw[-] (a)--(b);
    \draw[-] (b)--(c);

    \draw[-] (o)--(c);
    \draw[-] (q)--(c);
    \draw[-] (o)--(r);
    \draw[-] (s)--(q);

    \draw[-] (a)--(j);
    \draw[-] (a)--(l);
    \draw[-] (j)--(m);
    \draw[-] (l)--(n);

    \draw[-] (c)--(af);
    \draw[-] (af)--(ag);

\draw [decorate,decoration={brace,amplitude=5pt,raise=4ex}]
  (-0.5,.8) -- (0.5,.8) node[midway,yshift=3em]{$a_{i,n_i}-1$};

\draw [decorate,decoration={brace,amplitude=5pt,raise=4ex}]
  (1.4,.8) -- (2.4,.8) node[midway,yshift=3em]{$a_{i,1}-2$};
\end{tikzpicture}
\end{center}
Since the curve $F$ has no free point, the curve $A_{i,1}$ is an exceptional curve of a Wahl singularity(Lemma~\ref{lem:No free point ensure that the last curve is an exceptional of Wahl singularity}). Note that if we continue the MMP-algorithm until the $(-a_{i,1}-1)$-curve becomes a $(-2)$-curve, then we obtain the matrix $M_i''$.

Assume that $g_i \geq 3$.
\begin{center}
$M_i = \begin{bmatrix}[c|ccccccccccccccccccc]
                       & p_0 & q_{i,1} & \cdots & q_{i,g_i} \\
              \hline
               C_{i,1}   & 1      & *      & \cdots & *      & *      & \cdots & *      \\
               \vdots    & \vdots & \vdots & \ddots & \vdots & \vdots & \ddots & \vdots \\
               C_{i,m_i} & 1      & *      & \cdots & *      & *      & \cdots & *      \\
\end{bmatrix}$
\end{center}
Define a matrix $M_i'$ by deleting the column $p_0$ and $M_i''$ by deleting the columns $q_{i,1}, \cdots, q_{i,g_i}$.
\begin{center}
$M_i' = \begin{bmatrix}[c|ccccccccccccccccccc]
                       & q_{i,1} & \cdots & q_{i,g_i} \\
              \hline
               C_{i,1}   & *      & \cdots & *      & *      & \cdots & *      \\
               \vdots    & \vdots & \ddots & \vdots & \vdots & \ddots & \vdots \\
               C_{i,m_i} & *      & \cdots & *      & *      & \cdots & *      \\
\end{bmatrix}$
,  
$M_i'' = \begin{bmatrix}[c|ccccccccccccccccccc]
                         &        &        &         \\
              \hline
               C_{i,1}   & *      & \cdots & *       \\
               \vdots    & \vdots & \ddots & \vdots  \\
               C_{i,m_i} & *      & \cdots & *       \\
\end{bmatrix}$
\end{center}
This matrix satisfies the combinatorial equations of the quotient surface singularity $[a_{i,n_i}, \dots, a_{i,1}]$. Then we add a row $F$ to $M_i'$.
\begin{center}
\begin{minipage}{.4\textwidth}
\resizebox{\textwidth}{!}
{
$\begin{bmatrix}[c|ccccccccccccccccccc]
                       & q_{i,1} & \cdots & q_{i,g_i} \\
              \hline
               C_{i,1}   & *      & \cdots & *      & *      & \cdots & *      \\
               \vdots    & \vdots & \ddots & \vdots & \vdots & \ddots & \vdots \\
               C_{i,m_i} & *      & \cdots & *      & *      & \cdots & *      \\
               F         & 1      & \cdots & 1      &        &        &
\end{bmatrix}$
}
\end{minipage}
\end{center}
Note that the intersection relation between $F$ and $C_{i,1}, \cdots, C_{i,m_i}$ is the same with the relation between $C_{1,e}$ and $C_{i,1}, \cdots, C_{i,m_i}$. That is, $F.C_{i,j} = 1$ for $j = 1, \cdots, m_i$ and $l(F) = g_i$. Therefore this matrix satisfies the combinatorial equations of the following sandwiched structure where $(g_i-2)$ $(-2)$-curves are added.
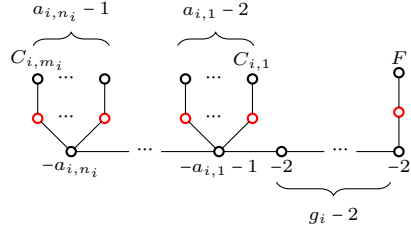
\begin{center}
\begin{tikzpicture}
  [inner sep=1mm,
R/.style={circle,draw=black!255,fill=white!20,thick},
T/.style={circle,draw=red!255,fill=white!255,thick}]
    \node[R] (a) [label=below:$-a_{i,n_i}$]{};
    \node    (b) [right=.7cm of a]{$\cdots$};
    \node[R] (c) [right=.7cm of b, label=below:$-a_{i,1}-1$]{};
    \node[R] (c1) [right=.7cm of c, label=below:$-2$]{};
    \node    (c2) [right=.5cm of c1]{$\cdots$};

    \node[R] (e) [right=.5cm of c2, label=below:$-2$]{};

    \node[T] (j) [above left=.5cm of a]{};
    \node    (k) [right=.1cm of j]{$\cdots$};
    \node[T] (l) [above right=.5cm of a]{};
    \node[R] (m) [above=.4cm of j, label=above:$C_{i,m_i}$]{};
    \node    (k) [right=.1cm of m]{$\cdots$};
    \node[R] (n) [above=.4cm of l]{};

    \node[T] (o) [above left=.5cm of c]{};
    \node    (p) [right=.1cm of o]{$\cdots$};
    \node[T] (q) [above right=.5cm of c]{};
    \node[R] (r) [above=.4cm of o]{};
    \node    (k) [right=.1cm of r]{$\cdots$};
    \node[R] (s) [above=.4cm of q, label=above:$C_{i,1}$]{};

    \node[T] (ad) [above=.4cm of e]{};
    \node[R] (ag) [above=.4cm of ad, label=above:$F$]{};

    \draw[-] (a)--(b);
    \draw[-] (b)--(c);
    \draw[-] (c)--(c1);
    \draw[-] (c1)--(c2);
    \draw[-] (c2)--(e);

    \draw[-] (o)--(c);
    \draw[-] (q)--(c);
    \draw[-] (o)--(r);
    \draw[-] (s)--(q);

    \draw[-] (a)--(j);
    \draw[-] (a)--(l);
    \draw[-] (j)--(m);
    \draw[-] (l)--(n);

    \draw[-] (ad)--(e);
    \draw[-] (ag)--(ad);

\draw [decorate,decoration={brace,amplitude=5pt,raise=4ex}]
  (-0.5,.8) -- (0.5,.8) node[midway,yshift=3em]{$a_{i,n_i}-1$};

\draw [decorate,decoration={brace,amplitude=5pt,raise=4ex}]
  (1.4,.8) -- (2.4,.8) node[midway,yshift=3em]{$a_{i,1}-2$};

\draw [decorate,decoration={brace,mirror,amplitude=5pt,raise=4ex}]
  (2.7,.2) -- (4.2,.2) node[midway,yshift=-3em]{$g_i-2$};
\end{tikzpicture}
\end{center}
The matrix is an incidence matrix of the cyclic quotient surface singularity with the sandwiched structure. Therefore we can find the corresponding P-resolution. Since the decorated curve $F$ has no free points, the last $(-2)$-curve is an exceptional curve of a Wahl singularity(Lemma~\ref{lem:No free point ensure that the last curve is an exceptional of Wahl singularity}). Note that if we progress the MMP-algorithm until the ($-a_{i,1}-1$)-curve becomes $(-2)$-curve, then we obtain the matrix $M_i''$.

We apply flips and divisorial contractions to the $P$-resolution of $X$. We progress until the $(-d-s)$-curve becomes $(-t-s)$-curve, the $-(a_{i,1})$-curve becomes a $(-1)$-curve if $g_i = 1$, $-(a_{i,1})$-curve becomes a $(-2)$-curve if $g_i \geq 1$, then we arrive at the followings.

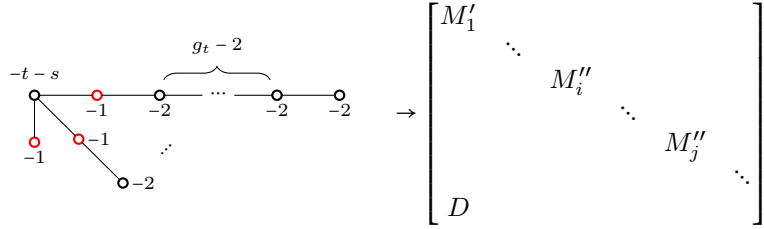
\begin{center}
\begin{minipage}{.4\textwidth}
\begin{tikzpicture}
  [inner sep=1mm,
R/.style={circle,draw=black!255,fill=white!20,thick},
T/.style={circle,draw=red!255,fill=white!255,thick}]
    \node[R] (c) [label={[yshift=.1cm]:$-t-s$}]{};

    \node[T] (a) [below=.5cm of c, label=below:$-1$]{};

    \node[T] (at1) [right=.7cm of c, label=below:$-1$]{};
    \node[R] (at2) [right=.7cm of at1, label=below:$-2$]{};
    \node    (at3) [right=.5cm of at2]{$\cdots$};
    \node[R] (at4) [right=.5cm of at3, label=below:$-2$]{};
    \node[R] (at5) [right=.7cm of at4, label=below:$-2$]{};

    \node[T] (a41) [below right=.7cm of c, label=right:$-1$]{}; 
    \node[R] (a42) [below right=.7cm of a41, label=right:$-2$]{};    

    \node    (e3) [below right=.5cm and 1.5cm of c]{$\udots$};

    \draw[-] (c)--(at1);
    \draw[-] (at1)--(at2);
    \draw[-] (at2)--(at3);
    \draw[-] (at3)--(at4);
    \draw[-] (at4)--(at5);

    \draw[-] (c)--(a41);
    \draw[-] (a41)--(a42);

    \draw[-] (c)--(a);

\draw [decorate,decoration={brace,amplitude=5pt,raise=4ex}]
  (1.7,-.4) -- (3.1,-.4) node[midway,yshift=3em]{$g_t - 2$};
\end{tikzpicture}
\end{minipage}
$\rightarrow$
\begin{minipage}{.4\textwidth}
$\begin{bmatrix}
  M_1' &         &       &        &       &        \\
       & \ddots  &       &        &       &        \\
       &         & M_i'' &        &       &        \\
       &         &       & \ddots &       &        \\
       &         &       &        & M_j'' &        \\
       &         &       &        &       & \ddots \\
  D    &         &       &        &       &        
\end{bmatrix}$
\end{minipage}
\end{center}

As we have noted at each step, we obtain the sub-matrices $[M_1', D], M_2'', \cdots, M_t''$. If $g_i = 1$, then all decorated curves $C_{i,1}, \cdots, C_{i,m_i}$ are connected to the central curve through the $(-1)$-curve. Since the decorated curves $C_{1,e}, \dots, C_{1,m_1}$ degenerate to the central curve, the $(-1)$-curve connects the decorated curves $C_{1,e}, \cdots, C_{1,m_1}$ and $C_{i,1}, \cdots, C_{i,m_i}$. The $(-1)$-curve corresponds to the column $q_{i,1}$.

If $g_i \geq 2$, then we can consider that the additional curve $F$ is replaced by the degenerated curves $C_{1, e}, \cdots, C_{1,m_1}$. This means that in the combinatorial incidence matrix, $F$ is replaced by $C_{1, e}, \cdots, C_{1,m_1}$ and we obtain the columns $q_{i,1}, \cdots, q_{i,g_i}$.

$\mathbf{Case~B-2})$ We assume that $[M_1, M_2, D]$ is the type 2-1 and $[M_1, M_i, D]$ is type 2-2 for $i = 3, \cdots, t$(Matrix~\ref{mat:caseB-2}). 

We define $t-2+s = \sum\limits_{i=3}^t g_i$. That is, $s$ is the sum of $g_i$ such that $g_i \geq 2$. We consecutively blow up the intersection of the central curve and $A_{i,1}$ in the minimal resolution of $(X, 0)$ to make $(g_i-2)$ $(-2)$-curves for $i = 3, \cdots, t$.
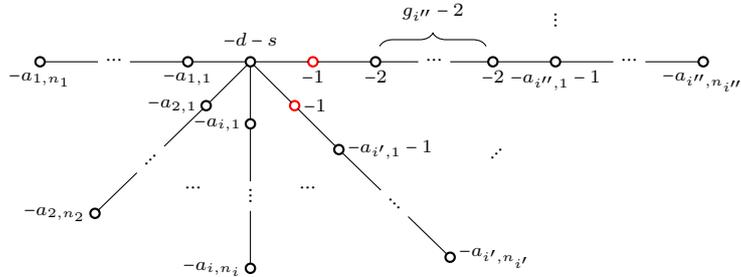
\begin{center}
\begin{tikzpicture}
  [inner sep=1mm,
R/.style={circle,draw=black!255,fill=white!20,thick},
T/.style={circle,draw=red!255,fill=white!255,thick}]
    \node[R] (c) [label={[yshift=.1cm]:$-d-s$}]{};

    \node[R] (a11) [left=.7cm of c, label=below:$-a_{1,1}$]{};
    \node    (a12) [left=.7cm of a11]{$\cdots$};
    \node[R] (a13) [left=.7cm of a12, label=below:$-a_{1,n_1}$]{};

    \node[T] (at1) [right=.7cm of c, label=below:$-1$]{};
    \node[R] (at2) [right=.7cm of at1, label=below:$-2$]{};
    \node    (at3) [right=.5cm of at2]{$\cdots$};
    \node[R] (at4) [right=.5cm of at3, label=below:$-2$]{};
    \node[R] (at5) [right=.7cm of at4, label=below:$-a_{i'',1}-1$]{};
    \node    (at6) [right=.7cm of at5]{$\cdots$};
    \node[R] (at7) [right=.7cm of at6, label=below:$-a_{i'',n_{i''}}$]{};

    \node[R] (a21) [below left=.7cm of c, label=left:$-a_{2,1}$]{};    
    \node    (a22) [below left=.7cm of a21]{$\udots$}; 
    \node[R] (a23) [below left=.7cm of a22, label=left:$-a_{2,n_2}$]{}; 

    \node[R] (a31) [below=.7cm of c, label=left:$-a_{i,1}$]{};    
    \node    (a32) [below=.7cm of a31]{$\vdots$}; 
    \node[R] (a33) [below=.7cm of a32, label=left:$-a_{i,n_i}$]{};      

    \node[T] (a41) [below right=.7cm of c, label=right:$-1$]{}; 
    \node[R] (a42) [below right=.7cm of a41, label=right:$-a_{i',1}-1$]{};    
    \node    (a43) [below right=.7cm of a42]{$\ddots$}; 
    \node[R] (a44) [below right=.7cm of a43, label=right:$-a_{i',n_{i'}}$]{};      

    \node    (e1) [below left=1.5cm and .5cm of c]{$\cdots$}; 
    \node    (e2) [below right=1.5cm and .5cm of c]{$\cdots$}; 
    \node    (e3) [below right=1cm and 3cm of c]{$\udots$}; 
    \node    (e4) [above=.3cm of at5]{$\vdots$}; 

    \draw[-] (c)--(a11);
    \draw[-] (a11)--(a12);
    \draw[-] (a12)--(a13);

    \draw[-] (c)--(at1);
    \draw[-] (at1)--(at2);
    \draw[-] (at2)--(at3);
    \draw[-] (at3)--(at4);
    \draw[-] (at4)--(at5);
    \draw[-] (at5)--(at6);
    \draw[-] (at6)--(at7);

    \draw[-] (c)--(a21);
    \draw[-] (a21)--(a22);
    \draw[-] (a22)--(a23);

    \draw[-] (c)--(a31);
    \draw[-] (a31)--(a32);
    \draw[-] (a32)--(a33);

    \draw[-] (c)--(a41);
    \draw[-] (a41)--(a42);
    \draw[-] (a42)--(a43);
    \draw[-] (a43)--(a44);

\draw [decorate,decoration={brace,amplitude=5pt,raise=4ex}]
  (1.7,-.4) -- (3.1,-.4) node[midway,yshift=3em]{$g_{i''} - 2$};
\end{tikzpicture}

\end{center}
We will locate T-singularities on $[a_{1,n_1},\dots,a_{1,1},-d-s, a_{2,1}, \dots, a_{2,n_2}]$ and $[a_{i,n_i},\dots,a_{i,1}+1, 2, \dots, 2]$ for $i = 3, \cdots, t$.

(1) Consider a matrix $[M_1, M_2, D]$(Figure~\ref{fig.c1c2d}). 

\begin{figure}
\centering
$\begin{bmatrix}[c|ccccccccccccccccccc]
                       & p_0 & p_1 & \cdots & p_{d-t} & q_1 & \cdots & q_{g'} & q_{3,1} & \cdots & q_{t,g_t} & \cdots \\
              \cmidrule(lr){1-12}
               C_{1,1} & 1      & \;     & \;     & \;     & \;     & \;     & \;      & \;     & \;     & \;     \\
               \vdots    & \vdots & \;     & \;     & \;     & \;     & \;     & \;      & \;     & \;     & \;     \\
               C_{1,e-1} & 1      & \;     & \;     & \;     & \;     & \;     & \;      & \;     & \;     & \;     \\  
               \cmidrule(lr){1-12}
               C_{1,e}   & \;     & 1      & \cdots & 1      & *      & \cdots & *       & 1      & \cdots & 1      \\
               \vdots    & \;     & \vdots & \ddots & \vdots & \vdots & \ddots & \vdots  & \vdots & \ddots & \vdots \\
               C_{1,m_1}   & \;     & 1      & \cdots & 1      & *      & \cdots & *     & 1      & \cdots & 1      \\
               \cmidrule(lr){1-12}
               C_{2,1}   & 1      & \;     & \;     & \;     & *      & \cdots & *       & \;      & \;      & \;      & \;         \\
               \vdots    & \vdots & \;     & \;     & \;     & \vdots & \ddots & \vdots  & \;      & \;      & \;      & \;         \\
               C_{2,m_2} & 1      & \;     & \;     & \;     & *      & \cdots & *       & \;      & \;      & \;      & \;         \\
               \cmidrule(lr){1-12}
               D_1       & 1      & 1      & 0      & 0      & \;     & \;     & \;      & \;      & \;      & \;      & \;         \\
               \vdots    & \vdots & 0      & \ddots & 0      & \;     & \;     & \;      & \;      & \;      & \;      & \;         \\
               D_{d-t-1} & 1      & 0      & 0      & 1      & \;     & \;     & \;      & \;      & \;      & \;      & \;         \\
\end{bmatrix}$
\caption{$[M_1, M_2, D]$}
\label{fig.c1c2d}
\end{figure}
We add decorated curves $E_1, \dots, E_{t-2+s}$ such that combinatorial equations are 
$$l(E_i) = 2$$
$$E_i.E_{i'} = 1$$
$$C_{1,j}.E_i = 1$$
$$C_{2,j'}.E_i = 1$$
for $i, i' = 1, \dots, t-2+s$.
\begin{figure}
\centering
$\begin{bmatrix}[c|ccccccccccccccccccc]
                       & p_0 & p_1 & \cdots & p_{d-t} & q_1 & \cdots & q_{g'}& q_{3,1} & \cdots & q_{t,g_t} & \cdots\\
              \cmidrule(lr){1-12}
               C_{1,1} & 1      & \;     & \;     & \;     & \;     & \;     & \;      & \;     & \;     & \;     \\
               \vdots    & \vdots & \;     & \;     & \;     & \;     & \;     & \;      & \;     & \;     & \;     \\
               C_{1,e-1} & 1      & \;     & \;     & \;     & \;     & \;     & \;      & \;     & \;     & \;     \\  
                \cmidrule(lr){1-12}
               C_{1,e}   & \;     & 1      & \cdots & 1      & *      & \cdots & *       & 1      & \cdots & 1      \\
               \vdots    & \;     & \vdots & \ddots & \vdots & \vdots & \ddots & \vdots  & \vdots & \ddots & \vdots \\
               C_{1,m_1}   & \;     & 1      & \cdots & 1      & *      & \cdots & *       & 1      & \cdots & 1      \\
               \cmidrule(lr){1-12}
               C_{2,1}   & 1      & \;     & \;     & \;     & *      & \cdots & *       & \;      & \;      & \;      & \;        \\
               \vdots    & \vdots & \;     & \;     & \;     & \vdots & \ddots & \vdots  & \;      & \;      & \;      & \;        \\
               C_{2,m_2} & 1      & \;     & \;     & \;     & *      & \cdots & *       & \;      & \;      & \;      & \;     \\
               \cmidrule(lr){1-12}
               D_1       & 1      & 1      & 0      & 0      & \;     & \;     & \;      & \;      & \;      & \;      & \;        \\
               \vdots    & \vdots & 0      & \ddots & 0      & \;     & \;     & \;      & \;      & \;      & \;      & \;         \\
               D_{d-t-1} & 1      & 0      & 0      & 1      & \;     & \;     & \;      & \;      & \;      & \;      & \;          \\
               \cmidrule(lr){1-12}
               E_1       & 1      & \;     & \;     & \;     & \;     & \;     & \;     & 1      & \cdots & 0      \\
               \vdots    & \vdots & \;     & \;     & \;     & \;     & \;     & \;     & \vdots & \ddots & \vdots \\
               E_{t-2+s}       & 1      & \;     & \;     & \;     & \;     & \;     & \;     & 0      & \cdots & 1       
\end{bmatrix}$
\caption{$[M_1,M_2,D,E]$}
\label{fig.c1c2de}
\end{figure}

The combinatorial equations of $E_i$ are the same with those of $D_k$. Therefore we can consider this matrix as an incidence matrix of a cyclic quotient surface singularity $[a_{1,n_1},\dots,a_{1,1},d+s, a_{2,1}, \dots, a_{2,n_2}]$(Ref Lemma~\ref{lem:twobranches}). We can find the $P$-resolution of the cyclic quotient surface singularity that induces the matrix $[M_1, M_2, D, E]$. 

Let $M_1'$ be a matrix obtained from $M_1$ by deleting the columns $q_{3,1}, \cdots, q_{t,g_t}$. If we progress the MMP-algorithm until the central curve becomes $(t - 2 + s - 1)$-curve, then we obtain the matrix $[M_1',M_2,D]$(See Figure~\ref{fig:partialmmp1}).
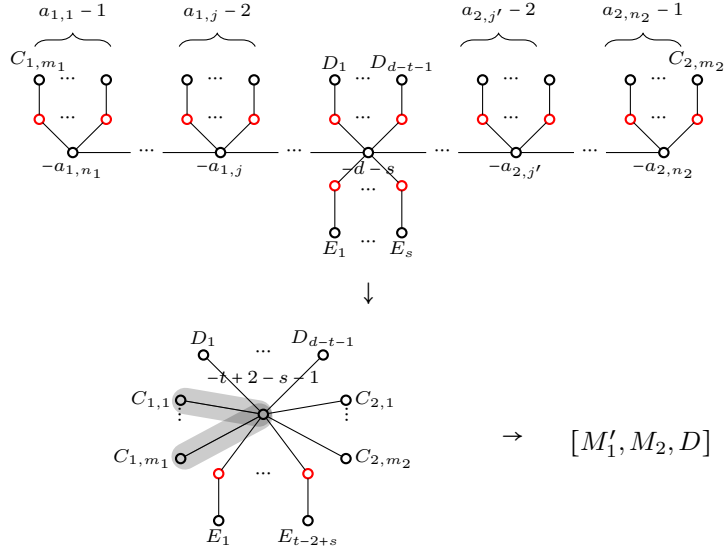
\begin{figure}
\begin{center}
\begin{minipage}{\textwidth}
\begin{center}
\begin{tikzpicture}
  [inner sep=1mm,
R/.style={circle,draw=black!255,fill=white!20,thick},
T/.style={circle,draw=red!255,fill=white!255,thick}]
    \node[R] (a) [label=below:$-a_{1,n_1}$]{};
    \node    (b) [right=.7cm of a]{$\cdots$};
    \node[R] (c) [right=.7cm of b, label=below:$-a_{1,j}$]{};
    \node    (d) [right=.7cm of c]{$\cdots$};
    \node[R] (e) [right=.7cm of d, label=below:$-d-s$]{};
    \node    (f) [right=.7cm of e]{$\cdots$};
    \node[R] (g) [right=.7cm of f, label=below:$-a_{2,j'}$]{};
    \node    (h) [right=.7cm of g]{$\cdots$};
    \node[R] (i) [right=.7cm of h, label=below:$-a_{2,n_2}$]{};

    \node[T] (j) [above left=.5cm of a]{};
    \node    (k) [right=.1cm of j]{$\cdots$};
    \node[T] (l) [above right=.5cm of a]{};
    \node[R] (m) [above=.4cm of j, label=above:$C_{1,m_1}$]{};
    \node    (k) [right=.1cm of m]{$\cdots$};
    \node[R] (n) [above=.4cm of l]{};

    \node[T] (o) [above left=.5cm of c]{};
    \node    (p) [right=.1cm of o]{$\cdots$};
    \node[T] (q) [above right=.5cm of c]{};
    \node[R] (r) [above=.4cm of o]{};
    \node    (k) [right=.1cm of r]{$\cdots$};
    \node[R] (s) [above=.4cm of q]{};

    \node[T] (t) [above left=.5cm of g]{};
    \node    (u) [right=.1cm of t]{$\cdots$};
    \node[T] (v) [above right=.5cm of g]{};
    \node[R] (w) [above=.4cm of t]{};
    \node    (k) [right=.1cm of w]{$\cdots$};
    \node[R] (x) [above=.4cm of v]{};

    \node[T] (y) [above left=.5cm of i]{};
    \node    (z) [right=.1cm of y]{$\cdots$};
    \node[T] (aa) [above right=.5cm of i]{};
    \node[R] (ab) [above=.4cm of y]{};
    \node    (k) [right=.1cm of ab]{$\cdots$};
    \node[R] (ac) [above=.4cm of aa, label=above:$C_{2,m_2}$]{};

    \node[T] (ad) [above left=.5cm of e]{};
    \node    (ae) [right=.1cm of ad]{$\cdots$};
    \node[T] (af) [above right=.5cm of e]{};
    \node[R] (ag) [above=.4cm of ad, label=above:$D_{1}$]{};
    \node    (k) [right=.1cm of ag]{$\cdots$};
    \node[R] (ah) [above=.4cm of af, label=above:$D_{d-t-1}$]{};

    \node[T] (ai) [below left=.5cm of e]{};
    \node    (aj) [below=.3cm of e]{$\cdots$};
    \node[T] (ak) [below right=.5cm of e]{};
    \node[R] (al) [below=.5cm of ai, label=below:$E_{1}$]{};
    \node    (am) [below=.4cm of aj]{$\cdots$};
    \node[R] (an) [below=.5cm of ak, label=below:$E_{s}$]{};

    \draw[-] (a)--(b);
    \draw[-] (b)--(c);
    \draw[-] (c)--(d);
    \draw[-] (d)--(e);
    \draw[-] (e)--(f);
    \draw[-] (f)--(g);
    \draw[-] (g)--(h);
    \draw[-] (h)--(i);
    
    \draw[-] (a)--(j);
    \draw[-] (a)--(l);
    \draw[-] (j)--(m);
    \draw[-] (l)--(n);

    \draw[-] (c)--(o);
    \draw[-] (c)--(q);
    \draw[-] (r)--(o);
    \draw[-] (s)--(q);

    \draw[-] (t)--(g);
    \draw[-] (v)--(g);
    \draw[-] (w)--(t);
    \draw[-] (x)--(v);

    \draw[-] (y)--(i);
    \draw[-] (aa)--(i);
    \draw[-] (ab)--(y);
    \draw[-] (ac)--(aa);

    \draw[-] (ad)--(e);
    \draw[-] (af)--(e);
    \draw[-] (ag)--(ad);
    \draw[-] (ah)--(af);

    \draw[-] (e)--(ai);
    \draw[-] (e)--(ak);
    \draw[-] (ai)--(al);
    \draw[-] (ak)--(an);

\draw [decorate,decoration={brace,amplitude=5pt,raise=4ex}]
  (-0.5,.8) -- (0.5,.8) node[midway,yshift=3em]{$a_{1,1}-1$};

\draw [decorate,decoration={brace,amplitude=5pt,raise=4ex}]
  (1.4,.8) -- (2.4,.8) node[midway,yshift=3em]{$a_{1,j}-2$};

\draw [decorate,decoration={brace,amplitude=5pt,raise=4ex}]
  (5.1,.8) -- (6.1,.8) node[midway,yshift=3em]{$a_{2,j'}-2$};

\draw [decorate,decoration={brace,amplitude=5pt,raise=4ex}]
  (7,.8) -- (8,.8) node[midway,yshift=3em]{$a_{2,n_2}-1$};
\end{tikzpicture}
\end{center}
\end{minipage}
\end{center}
\begin{center}
\begin{minipage}{.3\textwidth}
\begin{center}
$\downarrow$
\end{center}
\end{minipage}
\end{center}
\begin{center}
\begin{minipage}{.4\textwidth}
\begin{center}
\begin{tikzpicture}
[inner sep=1mm,
R/.style={circle,draw=black!255,fill=white!20,thick},
T/.style={circle,draw=red!255,fill=white!255,thick}]
    \node[R] (a) [label={[yshift=.3cm]:$-t+2-s-1$}]{};
    \node[R] (b) [above left=.1cm and 1cm of a, label=left:$C_{1,1}$]{};
    \node    (c) [left=.9cm of a]{$\vdots$};
    \node[R] (d) [below left=.5cm and 1cm of a, label=left:$C_{1,m_1}$]{};
    \node[T] (e) [below left=.7cm and .5cm of a]{};
    \node    (f) [below=.6cm of a]{$\cdots$};
    \node[T] (g) [below right=.7cm and .5cm of a]{};
    \node[R] (h) [below =.5cm of e, label=below:$E_1$]{};
    \node[R] (i) [below =.5cm of g, label=below:$E_{t-2+s}$]{};
    \node[R] (j) [above right=.1cm and 1cm of a, label=right:$C_{2,1}$]{};
    \node    (k) [right=.9cm of a]{$\vdots$};
    \node[R] (l) [below right=.5cm and 1cm of a, label=right:$C_{2,m_2}$]{};

    \node[R] (m) [above left=.7cm and .7cm of a, label=above:$D_{1}$]{};
    \node[R] (n) [above right=.7cm and .7cm of a, label=above:$D_{d-t-1}$]{};
    \node    (o) [above=.6cm of a]{$\cdots$};

    \draw[-] (a) -- (b);
    \draw[-] (a) -- (d);
    \draw[-] (a) -- (e);
    \draw[-] (a) -- (g);
    \draw[-] (e) -- (h);
    \draw[-] (i) -- (g);
    \draw[-] (a) -- (j);
    \draw[-] (a) -- (l);
    \draw[-] (a) -- (m);
    \draw[-] (a) -- (n);

    \draw[edge,line width=10pt] (a) -- (b);
    \draw[edge,line width=10pt] (a) -- (d);
\end{tikzpicture}
\end{center}
\end{minipage}
\begin{minipage}{.1\textwidth}
\begin{center}
$\to$
\end{center}
\end{minipage}
\begin{minipage}{.1\textwidth}
\begin{center}
$[M_1', M_2, D]$
\end{center}
\end{minipage}
\end{center}
\caption{partial MMP on $[M_1, M_2, D]$}
\label{fig:partialmmp1}
\end{figure}

(2) For the sub-matrices $M_i$ for $i = 3, \cdots, t$, we obtain $P$-resolutions as in the $\mathbf{Case~B-1}$.

(3) We found $P$-resolutions of $[a_{1,n_1},\dots, a_{1,1}, d+s, a_{2,1}, \dots, a_{2,n_2}]$ and $[a_{i,n_i}, \dots, a_{i,1}]$, we can locate $T$-singularities on the minimal resolution of $(X, 0)$.

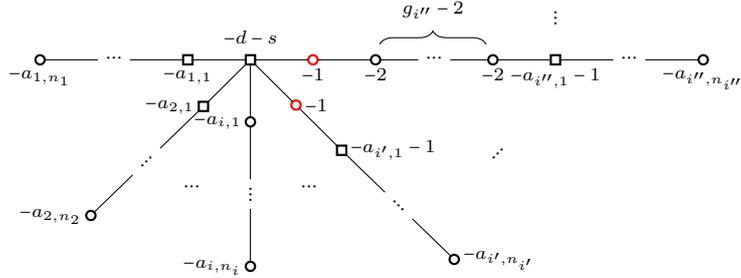
\begin{center}
\begin{tikzpicture}
  [inner sep=1mm,
R/.style={circle,draw=black!255,fill=white!20,thick},
S/.style={rectangle,draw=black!255,fill=white!20,thick},
T/.style={circle,draw=red!255,fill=white!255,thick}]
    \node[S] (c) [label={[yshift=.1cm]:$-d-s$}]{};

    \node[S] (a11) [left=.7cm of c, label=below:$-a_{1,1}$]{};
    \node    (a12) [left=.7cm of a11]{$\cdots$};
    \node[R] (a13) [left=.7cm of a12, label=below:$-a_{1,n_1}$]{};

    \node[T] (at1) [right=.7cm of c, label=below:$-1$]{};
    \node[R] (at2) [right=.7cm of at1, label=below:$-2$]{};
    \node    (at3) [right=.5cm of at2]{$\cdots$};
    \node[R] (at4) [right=.5cm of at3, label=below:$-2$]{};
    \node[S] (at5) [right=.7cm of at4, label=below:$-a_{i'',1}-1$]{};
    \node    (at6) [right=.7cm of at5]{$\cdots$};
    \node[R] (at7) [right=.7cm of at6, label=below:$-a_{i'',n_{i''}}$]{};

    \node[S] (a21) [below left=.7cm of c, label=left:$-a_{2,1}$]{};    
    \node    (a22) [below left=.7cm of a21]{$\udots$}; 
    \node[R] (a23) [below left=.7cm of a22, label=left:$-a_{2,n_2}$]{}; 

    \node[R] (a31) [below=.7cm of c, label=left:$-a_{i,1}$]{};    
    \node    (a32) [below=.7cm of a31]{$\vdots$}; 
    \node[R] (a33) [below=.7cm of a32, label=left:$-a_{i,n_i}$]{};      

    \node[T] (a41) [below right=.7cm of c, label=right:$-1$]{}; 
    \node[S] (a42) [below right=.7cm of a41, label=right:$-a_{i',1}-1$]{};    
    \node    (a43) [below right=.7cm of a42]{$\ddots$}; 
    \node[R] (a44) [below right=.7cm of a43, label=right:$-a_{i',n_{i'}}$]{};      

    \node    (e1) [below left=1.5cm and .5cm of c]{$\cdots$}; 
    \node    (e2) [below right=1.5cm and .5cm of c]{$\cdots$}; 
    \node    (e3) [below right=1cm and 3cm of c]{$\udots$}; 
    \node    (e4) [above=.3cm of at5]{$\vdots$}; 

    \draw[-] (c)--(a11);
    \draw[-] (a11)--(a12);
    \draw[-] (a12)--(a13);

    \draw[-] (c)--(at1);
    \draw[-] (at1)--(at2);
    \draw[-] (at2)--(at3);
    \draw[-] (at3)--(at4);
    \draw[-] (at4)--(at5);
    \draw[-] (at5)--(at6);
    \draw[-] (at6)--(at7);

    \draw[-] (c)--(a21);
    \draw[-] (a21)--(a22);
    \draw[-] (a22)--(a23);

    \draw[-] (c)--(a31);
    \draw[-] (a31)--(a32);
    \draw[-] (a32)--(a33);

    \draw[-] (c)--(a41);
    \draw[-] (a41)--(a42);
    \draw[-] (a42)--(a43);
    \draw[-] (a43)--(a44);

\draw [decorate,decoration={brace,amplitude=5pt,raise=4ex}]
  (1.7,-.4) -- (3.1,-.4) node[midway,yshift=3em]{$g_{i''} - 2$};
\end{tikzpicture}

\end{center}

(4) We run the MMP algorithm on the constructed $P$-resolution until the $(-d-s)$-curve becomes a $(-t-s+1)$-curve, the $-a_{i,1}$-curve becomes a $(-1)$-curve when $g_i = 2$ and the $-(a_{i,1}+1)$-curve becomes a $(-2)$-curve when $g_i \geq 2$, then we obtain the sub-matrices $[M_1', M_2, D], M_3'', \dots, M_t''$.
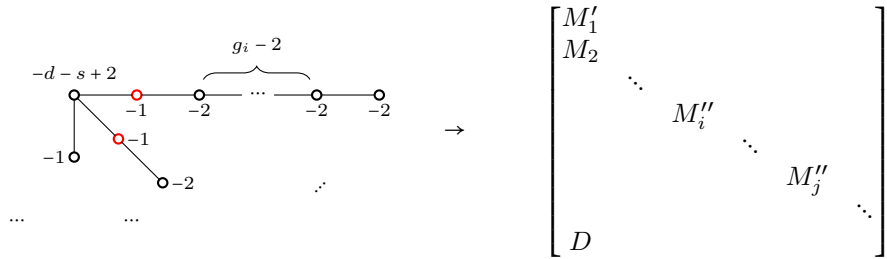
\begin{center}
\begin{minipage}{.45\textwidth}
\begin{tikzpicture}
  [inner sep=1mm,
R/.style={circle,draw=black!255,fill=white!20,thick},
T/.style={circle,draw=red!255,fill=white!255,thick}]
    \node[R] (c) [label={[yshift=.1cm]:$-d-s+2$}]{};

    \node[T] (at1) [right=.7cm of c, label=below:$-1$]{};
    \node[R] (at2) [right=.7cm of at1, label=below:$-2$]{};
    \node    (at3) [right=.5cm of at2]{$\cdots$};
    \node[R] (at4) [right=.5cm of at3, label=below:$-2$]{};
    \node[R] (at5) [right=.7cm of at4, label=below:$-2$]{};

    \node[R] (a31) [below=.7cm of c, label=left:$-1$]{};        

    \node[T] (a41) [below right=.7cm of c, label=right:$-1$]{}; 
    \node[R] (a42) [below right=.7cm of a41, label=right:$-2$]{};    

    \node    (e1) [below left=1.5cm and .5cm of c]{$\cdots$}; 
    \node    (e2) [below right=1.5cm and .5cm of c]{$\cdots$}; 
    \node    (e3) [below right=1cm and 3cm of c]{$\udots$}; 

    \draw[-] (c)--(at1);
    \draw[-] (at1)--(at2);
    \draw[-] (at2)--(at3);
    \draw[-] (at3)--(at4);
    \draw[-] (at4)--(at5);

    \draw[-] (c)--(a31);

    \draw[-] (c)--(a41);
    \draw[-] (a41)--(a42);

\draw [decorate,decoration={brace,amplitude=5pt,raise=4ex}]
  (1.7,-.4) -- (3.1,-.4) node[midway,yshift=3em]{$g_i - 2$};
\end{tikzpicture}
\end{minipage}
\begin{minipage}{.1\textwidth}
  $\to$
\end{minipage}
\begin{minipage}{.3\textwidth}
$\begin{bmatrix}
  M_1' &         &       &        &       &        \\
  M_2  &         &       &        &       &        \\
       & \ddots  &       &        &       &        \\
       &         & M_i'' &        &       &        \\
       &         &       & \ddots &       &        \\
       &         &       &        & M_j'' &        \\
       &         &       &        &       & \ddots \\
  D    &         &       &        &       &        
\end{bmatrix}$
\end{minipage}
\end{center}

Similar to the $\mathbf{Case~B-1}$, we can check that the $P$-resolution induces the given combinatorial incidence matrix.

In general, it is not guaranteed that the ampleness still holds after blow-ups. Therefore, in case~B, we have to check the ampleness on the $(-1)$-curves near the central curve. The ampleness is equivalent to that the sum of the discrepancies of two curves connected through a $(-1)$-curve being equal to or less than $-1$. 

Simplifying cases, we have
\\
(1) One blow up
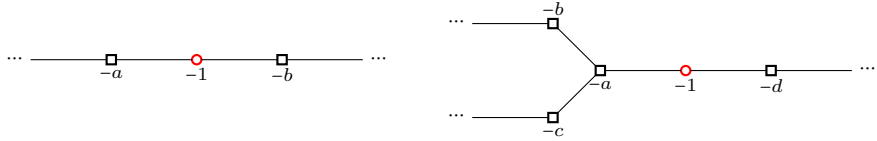
\begin{center}
\begin{minipage}{.45\textwidth}
  \begin{tikzpicture}
  [inner sep=1mm,
    R/.style={circle,draw=black!255,fill=white!20,thick},
    S/.style={rectangle,draw=black!255,fill=white!20,thick},
    T/.style={circle,draw=red!255,fill=white!255,thick}]
  
      \node    (a) []{$\cdots$};
      \node[S] (b) [right= of a, label=below:$-a$]{};
      \node[T] (c) [right= of b, label=below:$-1$]{};
      \node[S] (d) [right= of c, label=below:$-b$]{};
      \node    (e) [right= of d]{$\cdots$};

      \draw[-] (a) -- (b);
      \draw[-] (b) -- (c);
      \draw[-] (c) -- (d);
      \draw[-] (d) -- (e);

  \end{tikzpicture}
  \end{minipage}
  \begin{minipage}{.45\textwidth}
  \begin{tikzpicture}
  [inner sep=1mm,
    R/.style={circle,draw=black!255,fill=white!20,thick},
    S/.style={rectangle,draw=black!255,fill=white!20,thick},
    T/.style={circle,draw=red!255,fill=white!255,thick}]
  
      \node[T] (a) [label=below:$-1$]{};
      \node[S] (b) [left= of a, label=below:$-a$]{};
      \node[S] (c) [above left=.5cm and .5cm of b, label=above:$-b$]{};
      \node[S] (d) [below left=.5cm and .5cm of b, label=below:$-c$]{};
      \node    (f) [left= of c]{$\cdots$};
      \node    (g) [left= of d]{$\cdots$};
      \node[S] (h) [right= of a, label=below:$-d$]{};
      \node    (i) [right= of h]{$\cdots$};

      \draw[-] (a) -- (b);
      \draw[-] (b) -- (c);
      \draw[-] (b) -- (d);
      \draw[-] (c) -- (f);
      \draw[-] (d) -- (g);
      \draw[-] (a) -- (h);
      \draw[-] (h) -- (i);
    \end{tikzpicture}
  \end{minipage}
\end{center}
where $a \geq 5$ and $d \geq 3$.
\\
(2)Blow up $n$ times
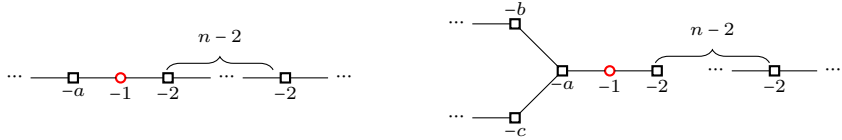
\begin{center}
\begin{minipage}{.45\textwidth}
  \begin{tikzpicture}
  [inner sep=1mm,
    R/.style={circle,draw=black!255,fill=white!20,thick},
    S/.style={rectangle,draw=black!255,fill=white!20,thick},
    T/.style={circle,draw=red!255,fill=white!255,thick}]
  
      \node    (a) []{$\cdots$};
      \node[S] (b) [right=.5cm of a, label=below:$-a$]{};
      \node[T] (c) [right=.5cm of b, label=below:$-1$]{};
      \node[S] (d) [right=.5cm of c, label=below:$-2$]{};
      \node    (e) [right=.5cm of d]{$\cdots$};
      \node[S] (f) [right=.5cm of e, label=below:$-2$]{};
      \node    (g) [right=.5cm of f]{$\cdots$};

      \draw[-] (a) -- (b);
      \draw[-] (b) -- (c);
      \draw[-] (c) -- (d);
      \draw[-] (d) -- (e);
      \draw[-] (e) -- (f);
      \draw[-] (f) -- (g);

\draw [decorate,decoration={brace,amplitude=5pt,raise=4ex}]
  (2,-.5) -- (3.4,-.5) node[midway,yshift=3em]{$n-2$};
  \end{tikzpicture}
  \end{minipage}
  \begin{minipage}{.45\textwidth}
  \begin{tikzpicture}
  [inner sep=1mm,
    R/.style={circle,draw=black!255,fill=white!20,thick},
    S/.style={rectangle,draw=black!255,fill=white!20,thick},
    T/.style={circle,draw=red!255,fill=white!255,thick}]
  
      \node[T] (a) [label=below:$-1$]{};
      \node[S] (b) [left=.5cm of a, label=below:$-a$]{};
      \node[S] (c) [above left=.5cm and .5cm of b, label=above:$-b$]{};
      \node[S] (d) [below left=.5cm and .5cm of b, label=below:$-c$]{};
      \node    (f) [left=.5cm of c]{$\cdots$};
      \node    (g) [left=.5cm of d]{$\cdots$};
      \node[S] (h) [right=.5cm of a, label=below:$-2$]{};
      \node    (i) [right=.5cm of h]{$\cdots$};
      \node[S] (j) [right=.5cm of i, label=below:$-2$]{};
      \node    (l) [right=.5cm of j]{$\cdots$};

      \draw[-] (a) -- (b);
      \draw[-] (b) -- (c);
      \draw[-] (b) -- (d);
      \draw[-] (c) -- (f);
      \draw[-] (d) -- (g);
      \draw[-] (a) -- (h);
      \draw[-] (i) -- (j);
      \draw[-] (j) -- (l);
      \draw [decorate,decoration={brace,amplitude=5pt,raise=4ex}]
  (.6,-.5) -- (2.1,-.5) node[midway,yshift=3em]{$n-2$};
    \end{tikzpicture}
  \end{minipage}
\end{center}
where $a \geq 4+n$.
Let the cases above be case $1, 2, 3$ and $4$.
We need some upper bounds for discrepancies of Wahl singularities. The following can be found in the appendix of \cite{choi2022symplectic}. 

We use the description of discrepancies of Urz\'{u}a-Vilches(\cite{UV}). Let $Y = \frac{1}{n^2}(1, na - 1)$ be a Wahl singularity and $f : \widetilde{Y} \rightarrow Y$ be the minimal resolution of $Y$. Then the canonical divisor $K_{\widetilde{Y}}$ of $\widetilde{Y}$ is represented as $K_{\widetilde{Y}} = f^* K_Y + \sum m_iE_i$ for exceptional curves $E_i$ of $f$. The $m_i$ is called the discrepancy of $E_i$. It is well known that $-1 < m_i < 0$ because $Y$ is a terminal singularity.

Let $[a_1, \cdots, a_r]$ be the Hirzebruch-Jung continued fraction of a Wahl singularity(Wahl continued fraction, for short) and $m_i$ be the discrepancy corresponding to $a_i$. We define an integer sequence $\delta_1, \cdots, \delta_r$ in the following inductive way. 

For $r = 1$, that is, for $[4]$, we assign an integer $\delta_1 = 1$ to $[4]$. If an integer sequence $\delta_1, \cdots, \delta_r$ is assigned to a Wahl singularity $[a_1, \cdots, a_r]$, then we assign 
$$\delta_1, \cdots, \delta_r, \delta_1 + \delta_r \text{ to } [a_1+1, a_2, \cdots, a_r, 2],$$
$$\delta_1 + \delta_r, \delta_1, \cdots, \delta_r \text{ to } [2, a_1, \cdots, a_{r-1}, a_r+1].$$
Then the discrepancy $m_i$ is $\left(-1 + \frac{\delta_i}{\delta_1 + \delta_r}\right)$.

\begin{lemma}[Urz\'{u}a-Vilches \cite{UV}*{Lemma 4.4}]
\label{lem:lemma of UV}
Let $[a_1, \cdots, a_t]$ be a Wahl singularity, assume $t \geq 2$ and $a_t = 2$, and let us denote its discrepancies by $m_1, \cdots, m_t$. Then we have the following bounds:\\
(Type M) If $a_2 = a_3 = \cdots = a_t$, then $m_1 = -1 + 1/(a_1-2)$ and $m_t = -1/(a_1-2)$.\\
(Type B) Otherwise, $m_1 = -1 + \mu$ and $m_t = -\mu$, where $1/a_1 < \mu < 1/(a_1 - 1)$.
\end{lemma}
\begin{lemma}
Let $[a_1, \cdots, a_r]$ be a Wahl continued fraction with $a_1 \geq 3$. Then the discrepancy of $a_1$ is less than $\frac{2-a_1}{a_1-1}$.
\end{lemma}
\begin{proof}
The last number $a_r$ must be $2$ because of the inductive construction of Wahl singularities. Then it is direct from Lemma~\ref{lem:lemma of UV}.
\end{proof}
In the case 1, let the discrepancies for the $(-a)$-curve be $m_a$ and the discrepancies for the $(-b)$-curve be $m_b$. Then $m_a + m_b < \frac{2-a}{a-1} + \frac{2-b}{b-1} < -\frac{1}{2} - \frac{1}{2} = -1$.
\begin{lemma}
Let $[a_1, \cdots, a_r]$ be a Wahl continued fraction with $a_1 = a_2 = \cdots = a_n = 2$ and $a_{n+1} \geq 3$ for $1 \leq n < r$. Then the discrepancy of $a_1$ is less than $-1/(n+2)$.
\end{lemma}
\begin{proof}
Consider the inverse of $[a_1, \cdots, a_r]$ : $[a_r, \cdots, a_1]$. Then $a_r$ must be $n+2$ because of the inductive construction of Wahl singularities. Then it follows directly from Lemma~\ref{lem:lemma of UV}.
\end{proof}
In the case 3, let the discrepancies for the $(-a)$-curve be $m_a$ and the discrepancies for the $(-b)$-curve be $m_b$. Then $m_a + m_b < \frac{2-a}{a-1} - \frac{1}{n+1} < \frac{-n-2}{n+3} - \frac{1}{n+1}  < \frac{-n^2-4n-5}{n^2 + 4n + 3} = -1 -\frac{2}{n^2+4n+3}  < -1$.
\begin{lemma}
Let $[a_1, \cdots, a_t, \cdots, a_r]$ be a Wahl continued fraction with $a_t \geq 5$. Then the discrepancy $m_t$ of $a_t$ is less than or equal to $(-a_t + 1)/a_t$.
\end{lemma}
\begin{proof}
Let $Y$ be the Wahl singularity corresponding to the given fraction. Then $K_{\widetilde{Y}} = f^*K_Y + \sum\limits_{i=1}^{r} m_i E_i$ where $E_i^2 = -a_i$. By multiplying $E_t$, we obtain $-2 + a_t = m_{t-1} + m_{t+1} -m_t a_t$. Therefore, $m_t = (2 - a_t + m_{t-1} + m_{t+1})/a_t$. If we show that $m_{t-1} + m_{t+1} \leq -1$, then we conclude that $m_t  \leq(-a_t + 1)/a_t$. 

We consider two cases. First, assume that $E_t$ is the initial curve of $Y$. Note that $[a_1, \cdots, a_t, \cdots, a_r]$ must be constructed from $[3, 5, 2]$ or $[2, 5, 3]$. Without loss of generality, assume that it is constructed from $[3, 5, 2]$. Then the $\delta$ sequence assigned to $[3, 5, 2]$ is $(2, 1, 3)$. If the sequence $(\delta_1, \cdots, \delta_r)$ is assigned to $[a_1, \cdots, a_t, \cdots, a_r]$, then $\delta_{t-1} = 2$ and $\delta_{t+1} = 3$. Note also that $\delta_1 + \delta_r \geq 2 + 3 = 5$. From the $\delta$ sequence, we obtain a bound $m_{t-1} + m_{t+1} = \left(-1 + \frac{\delta_{t-1}}{\delta_1 + \delta_r} \right) + \left(-1 + \frac{\delta_{t+1}}{\delta_1 + \delta_r} \right) = \left(-2 + \frac{\delta_{t-1} + \delta_{t+1}}{\delta_1 + \delta_r} \right) \leq -2 + \frac{5}{5} = -1$.

Second, assume that $E_t$ is not the initial curve and that the initial curve is left side of $E_t$. We track the inductive process to obtain $[a_1, \cdots, a_r]$. Starting from 
$$[4] \leftrightarrow (1),$$ 
we obtain 
$$[a_s-1, \cdots, a_{t-1}] \leftrightarrow (\delta_s, \cdots, \delta_{t-1}).$$ 
By adding a $2$ to the right side, we obtain 
$$[a_s, \cdots, a_{t-1}, 2] \leftrightarrow (\delta_s, \cdots, \delta_{t-1}, \delta_s + \delta_{t-1}).$$ 
To make the number $2$ to be $a_t$, we add $(a_t -2)$ $2$ to the left and we get 
$$[2, \cdots, 2, \cdots, a_t] \leftrightarrow ((a_t-1)\delta_s + (a_t-2)\delta_{t-1}, \cdots, 2\delta_u + \delta_{t-1}, \delta_s, \cdots, \delta_{t-1}, \delta_s + \delta_{t-1}).$$To fix the number $a_t$, we must add a $2$ to the right and we get 
$$[3, \cdots, 2, \cdots, a_t, 2]  \leftrightarrow ((a_t-1)\delta_s + (a_t-2)\delta_{t-1}, \cdots, 2\delta_u + \delta_{t-1}, \delta_s, \cdots, \delta_{t-1}, \delta_s + \delta_{t-1}, a_t\delta_s + (a_t-1)\delta_{t-1}).$$
Finally, we obtain
$$[a_1, \cdots, a_r] \leftrightarrow (\delta_1, \cdots, \delta_r).$$
Therefore we have 
\begin{equation*}
\begin{split}
  m_{t-1} + m_{t+1} & =  \left(-1+\frac{\delta_{t-1}}{\delta_1 + \delta_r}\right) + \left(-1+\frac{a_t\delta_s+(a_t-1)\delta_{t-1}}{\delta_1 + \delta_r}\right) \\
                    & = \left(-2+\frac{a_t\delta_s+a_t\delta_{t-1}}{\delta_1+\delta_t}\right) \\
                    & < \left(-2 + \frac{a_t\delta_s+a_t\delta_{t-1}}{(a_t-1)\delta_s + (a_t-2)\delta_{t-1} + a_t\delta_s+(a_t-1)\delta_{t-1})} \right) \\
                    & = \left(-2 + \frac{a_t\delta_s+a_t\delta_{t-1}}{(2a_t-1)\delta_s + (2a_t-3)\delta_{t-1}} \right) \\
                    & < -1 
\end{split}
\end{equation*}
\end{proof}
In the cases 2 and 4, let the discrepancies for the $(-a)$-curve be $m_a$, the discrepancies for the $(-d)$-curve be $m_d$ and the discrepancies for the $(-2)$-curve be $m_2$. Then $m_a + m_d < \frac{1-a_t}{a_t} + \frac{2-d}{d-1} < -\frac{4}{5} - \frac{1}{2}  < -1$. And $m_a + m_2 < \frac{1-a_t}{a_t} - \frac{1}{n+1} < -\frac{-3-n}{4+n}-\frac{1}{n+1} = \frac{-n^2 -5n -7}{n^2 + 5n + 4} = -1 - \frac{3}{n^2 + 5n + 4} < -1$. The ampleness of each case is confirmed.

\end{proof}
\begin{remark}
In the definition of case B, the condition `decorated curves that degenerate to the central curve come from only one branch ' is essential for finding the corresponding $P$-resolutions. Even if $d = t + 2$, if the condition is still satisfied, then we can construct $P$-resolutions in similar way. But there exist combinatorial incidence matrices that the condition is not satisfied. If $d = t + 2$, then there are incidence matrices that do not correspond to P-resolutions. 

For example, for a weighted homogeneous surface singularity of type $(6, (2,1), (2,1), (2,1), (2,1))$, we have the following combinatorial incidence matrix.
$$\begin{bmatrix}[c|cccccc]
C_1 &  1 &   & 1 & 1 &   &   \\
C_2 &  1 &   &   &   & 1 & 1 \\
C_3 &    & 1 & 1 &   & 1 &   \\
C_4 &    & 1 &   & 1 &   &   \\
D_1 &  1 & 1 &   &   &   & 1 \\
\end{bmatrix}$$
We expect that it is a non-cyclic normal singularity admitting a $\mathbb{Q}$-Gorenstein smoothing.
\end{remark}

\subsection{An example}
We consider a weighted homogeneous surface singularity of type $(6, (3,5), (9,13), (7,10))$. Then its dual resolution graph is Figure~\ref{fig:a WHSS of type (6, (3,5), (9,13), (7,10))}.
\begin{figure}
\label{fig:a WHSS of type (6, (3,5), (9,13), (7,10))}
\centering
  \begin{tikzpicture}
    [inner sep=1mm,
  R/.style={circle,draw=black!255,fill=white!20,thick},
  T/.style={rectangle,draw=black!255,fill=white!255,thick}]
      \node[R] (a) [label=above:$-3$]{};
      \node[R] (b) [right=of a, label=above:$-2$] {};
      \node[R] (c) [right=of b, label=above:$-6$] {};
      \node[R] (d) [right=of c, label=above:$-2$] {};
      \node[R] (e) [right=of d, label=above:$-2$] {};
      \node[R] (f) [right=of e, label=above:$-4$] {};
      \node[R] (g) [below=of c, label=right:$-2$] {};
      \node[R] (h) [below=of g, label=right:$-2$] {};
      \node[R] (i) [below=of h, label=right:$-5$] {};
      \draw[-] (a)--(b);
      \draw[-] (b)--(c);
      \draw[-] (c)--(d);
      \draw[-] (d)--(e);
      \draw[-] (e)--(f);
      \draw[-] (c)--(g);
      \draw[-] (g)--(h);
      \draw[-] (h)--(i);
  \end{tikzpicture}
  \caption{Dual resolution graph of a WHSS of type $(6, (3,5), (9,13), (7,10))$}
\end{figure}

$\mathbf{Case~A)}$ We have a following combinatorial incidence matrix of $\mathbf{Case~A}$
\begin{equation*}
\begin{bmatrix}[c|ccccccccccccccccc]
          & p_0  & \; & \; & \; & \; & \; & \; & \; & \; & \; & \; & \; & \; & \; & \; & \; & \; \\
\cmidrule(lr){1-18}
C_{1,1}   & 1    & 1  & 1  & 1  & 0  & \; & \; & \; & \; & \; & \; & \; & \; & \; & \; & \; & \; \\
C_{1,2}   & 1    & 1  & 1  & 0  & 1  & \; & \; & \; & \; & \; & \; & \; & \; & \; & \; & \; & \; \\
\cmidrule(lr){1-18}
C_{2,1}   & 1    &    &    &    &    & 1  & 1  & 1  & 1  & 0  & \; & \; & \; & \; & \; & \; & \; \\
C_{2,2}   & 1    &    &    &    &    & 1  & 1  & 1  & 0  & 1  & \; & \; & \; & \; & \; & \; & \; \\
C_{2,3}   & 1    &    &    &    &    & 1  & 1  & 0  & 1  & 1  & \; & \; & \; & \; & \; & \; & \; \\
C_{2,4}   & 1    &    &    &    &    & 1  & 0  & 1  & 1  & 1  & \; & \; & \; & \; & \; & \; & \; \\
\cmidrule(lr){1-18}
C_{3,1}   & 1    &    &    &    &    & \; & \; & \; & \; & \; & 1  & 1  & 1  & 1  & 0  & \; & \; \\
C_{3,2}   & 1    &    &    &    &    & \; & \; & \; & \; & \; & 1  & 1  & 1  & 0  & 1  & \; & \; \\
C_{3,3}   & 1    &    &    &    &    & \; & \; & \; & \; & \; & 1  & 1  & 0  & 1  & 1  & \; & \; \\
\cmidrule(lr){1-18}
D_1       & 1    &    &    &    &    & \; & \; & \; & \; & \; & \; & \; & \; & \; & \; & 1  & 0  \\
D_2       & 1    &    &    &    &    & \; & \; & \; & \; & \; & \; & \; & \; & \; & \; & 0  & 1  \\
\end{bmatrix}
\end{equation*}
Then we obtain three sub-matrices\\
\begin{equation*}
\begin{bmatrix}[c|cccc]
C_{1,1}   & 1  & 1  & 1  & 0  \\
C_{1,2}   & 1  & 1  & 0  & 1  \\
\end{bmatrix}
\;\;
\begin{bmatrix}[c|ccccc]
C_{2,1}   & 1  & 1  & 1  & 1  & 0   \\
C_{2,2}   & 1  & 1  & 1  & 0  & 1   \\
C_{2,3}   & 1  & 1  & 0  & 1  & 1   \\
C_{2,4}   & 1  & 0  & 1  & 1  & 1   \\
\end{bmatrix}
\;\;
\begin{bmatrix}[c|ccccc]
C_{3,1}   & 1  & 1  & 1  & 1  & 0 \\
C_{3,2}   & 1  & 1  & 1  & 0  & 1 \\
C_{3,3}   & 1  & 1  & 0  & 1  & 1 \\
\end{bmatrix}
\end{equation*}
We can find corresponding $P$-resolutions.\\
\begin{center}
\begin{tikzpicture}
  [inner sep=1mm,
R/.style={circle,draw=black!255,fill=white!20,thick},
T/.style={rectangle,draw=black!255,fill=white!255,thick}]
    \node[R] (a) [label=below:$-3$]{};
    \node[R] (b) [right=of a, label=below:$-2$] {};
    \draw[-] (a)--(b);
\end{tikzpicture}
\begin{tikzpicture}
  [inner sep=1mm,
R/.style={circle,draw=black!255,fill=white!20,thick},
T/.style={rectangle,draw=black!255,fill=white!255,thick}]
    \node[T] (a) [label=below:$-5$]{};
    \node[T] (b) [right=of a, label=below:$-2$] {};
    \node[R] (c) [right=of b, label=below:$-2$] {};
    \draw[-] (a)--(b);
    \draw[-] (b)--(c);
\end{tikzpicture}
  \begin{tikzpicture}
    [inner sep=1mm,
  R/.style={circle,draw=black!255,fill=white!20,thick},
  T/.style={rectangle,draw=black!255,fill=white!255,thick}]
      \node[T] (a) [label=below:$-4$]{};
      \node[R] (b) [right=of a, label=below:$-2$] {};
      \node[R] (c) [right=of b, label=below:$-2$] {};
      \draw[-] (a)--(b);
      \draw[-] (b)--(c);
  \end{tikzpicture}
  \end{center}
  From these $P$-resolutions, we get a $P$-resolution of $X$.
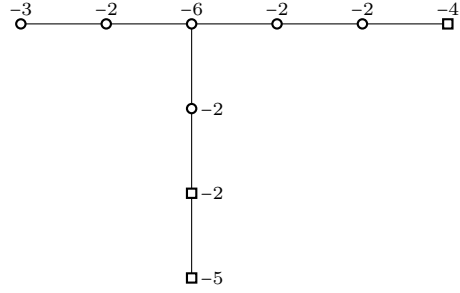
  \begin{center}
  \begin{tikzpicture}
    [inner sep=1mm,
  R/.style={circle,draw=black!255,fill=white!20,thick},
  T/.style={rectangle,draw=black!255,fill=white!255,thick}]
      \node[R] (a) [label=above:$-3$]{};
      \node[R] (b) [right=of a, label=above:$-2$] {};
      \node[R] (c) [right=of b, label=above:$-6$] {};
      \node[R] (d) [right=of c, label=above:$-2$] {};
      \node[R] (e) [right=of d, label=above:$-2$] {};
      \node[T] (f) [right=of e, label=above:$-4$] {};
      \node[R] (g) [below=of c, label=right:$-2$] {};
      \node[T] (h) [below=of g, label=right:$-2$] {};
      \node[T] (i) [below=of h, label=right:$-5$] {};
      \draw[-] (a)--(b);
      \draw[-] (b)--(c);
      \draw[-] (c)--(d);
      \draw[-] (d)--(e);
      \draw[-] (e)--(f);
      \draw[-] (c)--(g);
      \draw[-] (g)--(h);
      \draw[-] (h)--(i);
  \end{tikzpicture} 
  \end{center}
$\mathbf{Case~B)}$
We have an incidence matrix of case B.
\begin{center}
$\begin{bmatrix}[c|cccccccccccccc]
              C_{1,1} & \; &  1 &  1 & 1  &  1 &  0 & \; & \; & \; & \; & \; & \; & \; \\
              C_{1,2} & \; &  1 &  1 & 1  &  0 &  1 & \; & \; & \; & \; & \; & \; & \; \\
              C_{2,1} & 1  & \; & \; & 1  & \; & \; & 1  & 1  & 1  & 0  & \; & \; & \; \\
              C_{2,2} & 1  & \; & \; & 1  & \; & \; & 1  & 1  & 0  & 1  & \; & \; & \; \\
              C_{2,3} & 1  & \; & \; & 1  & \; & \; & 1  & 0  & 1  & 1  & \; & \; & \; \\
              C_{2,4} & 1  & \; & \; & 1  & \; & \; & 0  & 1  & 1  & 1  & \; & \; & \; \\  
              C_{3,1} & 1  & \; & \; & \; & 1  & 1  & \; & \; & \; & \; & 1  & 1  & 0  \\
              C_{3,2} & 1  & \; & \; & \; & 1  & 1  & \; & \; & \; & \; & 1  & 0  & 1  \\
              C_{3,3} & 1  & \; & \; & \; & 1  & 1  & \; & \; & \; & \; & 0  & 1  & 1  \\
              D_1     & 1  & 1  & 0  & \; & \; & \; & \; & \; & \; & \; & \; & \; & \; \\
              D_2     & 1  & 0  & 1  & \; & \; & \; & \; & \; & \; & \; & \; & \; & \; \\
\end{bmatrix}$
\end{center}
Then $[M_1,M_3,D]$ is of type 2-1 and $g' = 2$. And $[M_1, M_2,D]$ is of type 2-2 and $g_2 = 1$. Therefore we obtain two sub-matrices.
\begin{center}
$\begin{bmatrix}[c|cccccccccccccc]
              C_{1,1} & \; &  1 &  1 & 1  &  1 &  0 & \; & \; & \; \\
              C_{1,2} & \; &  1 &  1 & 1  &  0 &  1 & \; & \; & \; \\
              C_{3,1} & 1  & \; & \; & \; & 1  & 1  & 1  & 1  & 0  \\
              C_{3,2} & 1  & \; & \; & \; & 1  & 1  & 1  & 0  & 1  \\
              C_{3,3} & 1  & \; & \; & \; & 1  & 1  & 0  & 1  & 1  \\
              D_1     & 1  & 1  & 0  & 0  & \; & \; & \; & \; & \; \\
              D_2     & 1  & 0  & 1  & 0  & \; & \; & \; & \; & \; \\
              F_1     & 1  & 0  & 0  & 1  & \; & \; & \; & \; & \; \\              
\end{bmatrix}$
$\begin{bmatrix}[c|cccccccccccccc]
              C_{2,1} & 1  & 1  & 1  & 1  & 0  \\
              C_{2,2} & 1  & 1  & 1  & 0  & 1  \\
              C_{2,3} & 1  & 1  & 0  & 1  & 1  \\
              C_{2,4} & 1  & 0  & 1  & 1  & 1  \\  
\end{bmatrix}$
\end{center}
From these incidence matrices, we obtain corresponding $P$-resolutions.
\begin{center}
\begin{tikzpicture}
  [inner sep=1mm,
R/.style={circle,draw=black!255,fill=white!20,thick},
T/.style={rectangle,draw=black!255,fill=white!255,thick}]
    \node[T] (a) [label=below:$-3$]{};
    \node[T] (b) [right=of a, label=below:$-2$] {};
    \node[T] (c) [right=of b, label=below:$-6$] {};
    \node[T] (d) [right=of c, label=below:$-2$] {};
    \node[R] (e) [right=of d, label=below:$-2$] {};
    \node[R] (f) [right=of e, label=below:$-4$] {};

    \draw[-] (a)--(b);
    \draw[-] (b)--(c);
    \draw[-] (c)--(d);
    \draw[-] (d)--(e);
    \draw[-] (e)--(f);
\end{tikzpicture}
\begin{tikzpicture}
  [inner sep=1mm,
R/.style={circle,draw=black!255,fill=white!20,thick},
T/.style={rectangle,draw=black!255,fill=white!255,thick}]
    \node[T] (a) [label=below:$-5$]{};
    \node[T] (b) [right=of a, label=below:$-2$] {};
    \node[R] (c) [right=of b, label=below:$-2$] {};

    \draw[-] (a)--(b);
    \draw[-] (b)--(c);
\end{tikzpicture}
\end{center}
By combining them, we get a $P$-resolution of $X$.
\begin{center}
\begin{tikzpicture}
  [inner sep=1mm,
R/.style={circle,draw=black!255,fill=white!20,thick},
T/.style={rectangle,draw=black!255,fill=white!255,thick}]
    \node[T] (a) [label=below:$-3$]{};
    \node[T] (b) [right=of a, label=below:$-2$] {};
    \node[T] (c) [right=of b, label=below:$-6$] {};
    \node[T] (d) [right=of c, label=below:$-2$] {};
    \node[R] (e) [right=of d, label=below:$-2$] {};
    \node[T] (f) [right=of e, label=below:$-4$] {};
    \node[R] (g) [above=of c, label=right:$-2$] {};
    \node[T] (h) [above=of g, label=right:$-2$] {};
    \node[T] (i) [above=of h, label=right:$-5$] {};
    \draw[-] (a)--(b);
    \draw[-] (b)--(c);
    \draw[-] (c)--(d);
    \draw[-] (d)--(e);
    \draw[-] (e)--(f);
    \draw[-] (c)--(g);
    \draw[-] (g)--(h);
    \draw[-] (h)--(i);
\end{tikzpicture}
\end{center}

\end{document}